\documentclass[aip,jmp]{revtex4-1}
\usepackage{}
\usepackage{amsfonts}
\usepackage{graphicx}
\usepackage{dcolumn}
\usepackage{bm}
\usepackage{amssymb}
\usepackage{amsmath,amsthm}
\makeatletter

\newcommand{\Rmnum}[1]{\expandafter\@slowromancap\romannumeral #1@}
\makeatother
\newfam\msbfam
\font\twelmsb=msbm10 at 12pt 
\font\sevenmsb=msbm10 at 7pt \font\fivemsb=msbm10 at 5pt
\textfont\msbfam=\twelmsb
\scriptfont\msbfam=\sevenmsb \scriptscriptfont\msbfam=\fivemsb

  \def \beginproof{\par\noindent {\bf Proof.}\ \ }
  \def \endproof{\hskip .5cm $\Box$ \vskip .5cm}

 \newtheorem{proposition}{Proposition}[section]

 \newtheorem{lemma}{Lemma}[section]

 \newtheorem{theorem}{Theorem}[section]
 
\allowdisplaybreaks
\numberwithin{equation}{section}

\begin{document}


\title { The leading coefficients for an affine Weyl group of type $\tilde{G_{2}}$ I:\\
The lowest two-sided cell case}

\author{Peng-Fei Guo,  Hai-Tao Ma, Zhu-Jun Zheng}
\thanks{Supported in part by the Fundamental Research Funds for the central Universities with grant Number 2013ZM0109} \email{zhengzj@scut.edu.cn}
\affiliation{Department of Mathematics, South China University of
 Technology, Guangzhou 510641, P. R. China.}





\begin{abstract} \label{abstract}
This paper studies the leading coefficients $\mu(u,w)$ of the
Kazhdan-Lusztig polynomials $P_{u,w}$ for the lowest two-sided cell ${c_{0}}$ of an affine Weyl group $W$ of type
$\widetilde{G_{2}}$. We gives an estimation $\mu(u,w)\leq 3$ for $u,w\in c_{0}$ through calculations of multiplications
of Kazhdan-Lusztig basis and some conclusions of representation of algebraic group.
\end{abstract}

\maketitle

\maketitle

\section{Introduction}
The Hecke algebra $\mathcal {H}$ of Coxeter group $W$ over $\mathbb{Z}[q^{\frac{1}{2}},q^{-\frac{1}{2}}]$ is a free $\mathcal {A}$-module with a basis $\{T_{w}\}_{w\in\widetilde{W}}$. Its multiplication is defined by the relations $(T_{s}-q)(T_{s}+1)=0$ if $s\in S$ and $T_{w}T_{u} = T_{wu}$ if $l(wu)=l(w)+l(u)$. David Kazhdan and George Lusztig proved that there exist an unique baisis $C_{w}$ in $\mathcal {H}$ such that $C_{w} = q^{-\frac{l(w)}{2}}\sum_{u\leq w}P_{u, w}T_{u}$ and $\overline{C_{w}}=C_{w}$ (see [1]) for each $w\in W$, where $P_{u,w}\in \mathbb{Z}[q]$ of degree less than or equal to $\frac{1}{2}(l(w)-l(u)-1)$ for $u \leq w$ and $P_{w,w} = 1$ in [1]. It is called Kazhdan-Lusztig basis. The polynomial $P_{u, w}$ is called Kazhdan-Lusztig polynomial.

Many related research and results have shown the Kazhdan-Lusztig basis and Kazhdan-Lusztig polynomials are very important. Kazhdan-Lusztig cells of $W$, which are defined by Kazhdan-Lusztig polynomial, gives rise to a W-graph hence to a representation of Hecke algebra $\mathcal {H}$. Considering $W$ as a Weyl group and $X$ as a left cell of W, an integral representation of $W$ is obtained from $X$ by taking $q^{\frac{1}{2}}$ to 1. The correspondence representation over $\mathbb{Q}$ is not irreducible in general. However, it contains a unique irreducible component in the class $\mathcal{S}_{W}$ (see [25]) which is a distinguished class of irreducible representations of a Weyl group $W$. Two left cells give rise to the same representation in $\mathcal{S}_{W}$ if and only if they are contained in the same two-sided cell.

Let $G$ be a semisimple algebric group over an algebraically closed field with Weyl group $W$. For each $w\in W$, let $\mathcal {B}_{w}$ be the Schubert cell associated with $w$ and $\overline{\mathcal {B}_{w}}$ be the Schubert variety associated with $w$. Then $P_{y,w}$ can be regarded as a measure of the failure of local $Poincar\acute{e}$ duality on the Schubert variety $\overline{\mathcal {B}_{w}}$ in a neighborhood of a point in $\mathcal {B}_{y}$.

Let $W$ be Weyl group of  a semisimple complex Lie algebra $\mathbf{g}$. For each $w\in W$, let $M_{w}$ be the Verma module with the highest weight $-w(\rho)-\rho$ and let $L_{w}$ be its unique irreducible quotient. There are relations between characters of $L_{w}$ and $ M_{y}$. $chL_{w}=\sum_{y\leqslant w}(-1)^{l(y)-l(w)}P_{y,w}chM_{y}$ and $chM_{w}=\sum_{y\leqslant w}P_{w_{0}y,w_{0}w}chL_{y}$. They are the well known Kazhdan-Lusztig conjectures which have been proved independently by Beilinson Alexandre, Bernstein Joseph(see[26]) and by Brylinski Jean-Luc, Kashiwara Masaki(see[26]), respectively.

The leading coefficient $\mu(u, w)$ of $P_{u,w}$ plays important roles in Kazhdan-Lusztig theory, Lie theory and representation theory. It can be used to estimate the dimensions of the first extension groups for irreducible modules of algebraic groups and finite groups of Lie type(see[20]). Some $\mu(u, w)$ play crucial roles in the recursive formula of Kazhdan-Lusztig polynomials in [1].

However, there are few results for the leading coefficients in general. In [8], in order to understand the local finite problem of $\widetilde{B_{2}}$ type, George Lusztig computed some special the leading coefficients for an affine Weyl group of this type. In [19], for an affine Weyl group of type $A_{5}$, some non-trivial the leading coefficients are found. Timothy McLarnan and Gregory Warrington have shown that $\mu(u, w)$ can be greater than 1 for a symmetric group which is a counter example for $(0, 1)$-conjecture(see [9]). In [17], Nanhua Xi showed that when $u\leq w$ and $a(u)\leq a(w)$, then $\mu(u, w)\leq1$ if $W$ is a symmetric group or an affine Weyl group of type $A_{n}$. In [18], Leonard Scott and Nanhua Xi showed that the leading coefficient $\mu(u, w)$ of some Kazhdan-Lusztig polynomial $P_{u,w}$ for an affine Weyl group of type $A_{n}$ is $n+2$ if $n\geq4$. Recently, Brant Jones and Richard Green (see [21, 22]) showed that the coefficients are $\leq1$ for some special pairs of elements in a symmetric group and other Weyl groups. In [13, 14], Liping Wang gave some results on the leading coefficients for an affine Weyl group of type $\widetilde{B_{2}}$ and $\widetilde{A_{2}}$.

It's tedious to computed the leading coefficients for an affine Weyl group of $\widetilde{G_{2}}$ type. In [23], Nanhua Xi gave the explicit description of all five two-sided cells of affine Weyl group of type $\widetilde{G_{2}}$, $c_{1} = \{w \in W | a(w) = 0\} = \{e\}, c_{2} = \{w \in W | a(w) = 1\}, c_{3} = \{w \in W | a(w) = 2\}, c_{4}= {w \in W | a(w) = 3}, c_{0} = \{w \in W | a(w) = 6\}$, where $a(\cdot)$ is $a$-function (see[5]). We compute $\mu(u, w)$ through following two cases.

(1) Compute $\mu(u, w)$ for $u,v$  which belong to the same cell.

(i) For $c_{1}=\{e\}$, $\mu(e, e)=0$;

(ii) For $c_{2}$, this two-sided cell is called the second lowest cell, which means that $w$ has the unique expression form for any $w\in c_{2}$. In [13], for any $u, w\in c_{2}$, Liping Wang showed that $\mu(u, w)=1$ if $l(w)-l(u)=1$; otherwise, $\mu(u, w)=0$;

(iii) For $c_{3}$ and $c_{4}$, we use some semi-linear equations(see[8]) to compute $\mu(u, w)$, and this part is given in another article which we have done;

(iv) For $c_{0}$ which is the lowest two-sided cell, we  use some conclusions of representation of algebraic groups (see[18]) to compute all $\mu(u, w)$ for $u, w\in c_{0}$.

(2) Compute $\mu(u, w)$ for $u,v$ which belong to different cells.

In this paper, we calculate the leading coefficients $\mu(x, y)$ for $x, y\in c_{0}$. And the coming papers (see[29]) we give the results of the leading coefficients for $c_{3}$, $c_{4}$ and different cells.

In this paper, we use the conclusions in [18] to calculate $\mu(u, w)$ for $u, w\in c_{0}$. The calculation is so complicated that it can not be done manually. Although there are some programs to compute $P_{u,v}$, such as Mark Goresky's Pascal, Fokko du Cloux's Coxeter 3.0, these algorithms are not efficient for our problem. Our main work is find a new algorithm to compute $P_{u,v}$ and $C_{s}C_{y}$ for affine Weyl group of type $\widetilde{G_{2}}$ and realize it by visual C++. Moreover, by some long and tedious calculations, we obtain all $\delta_{x,y,z}$ defined in [3] which play an important role in the proof of our main result theorem IV.1.

We now outline the structure of this paper. In section II, we recall the definitions and properties Kazhdan-Lusztig basis, Kazhdan-Lusztig polynomial, Kazhdan-Lusztig cell and Kazhdan-Lusztig coefficient. We compute all elements of $c_{0}$ in Section II.C by the method mentioned in [20]. In Section 3, we give an explicit description of the algorithm which is used to compute Kazhdan-Lusztig polynomials and multiplications of Kazhdan-Lusztig basis. In Section 4, we compute all multiplications of Kazhdan-Lusztig basis which we need. From these multiplications, we get the main theorem of this paper that $\mu(u, v)\leqslant 3$ for $u, w\in c_{0}$.

  \section{Preliminaries}\label{sec2}

In this section, we review some basic facts for Kazhdan-Lusztig theory, properties of Kazhdan-Lusztig polynomial and Kazhdan-Lusztig coefficient, lowest two-sided cell and the formula of Kazhdan-Lusztig coefficient for the lowest two-sided cell.

\subsection{Basic definitions and conventions}

Let $R$ be a root system and $W_{0}$ be its Weyl group. The semi-direct product $W=W_{0}\ltimes\Lambda_{r}$
is an affine Weyl group and $\widetilde{W}=W_{0}\ltimes\Lambda$ is an extended affine Weyl group where $\Lambda$ is the weight lattice
of $R$ and $\Lambda_{r}$ = $\mathbb{Z}R$ is the root lattice.

There is an abelian
subgroup $\Omega$ of $\widetilde{W}$ such that $\omega S = S\omega$ for any $\omega\in\Omega$ and $\widetilde{W}=W \ltimes\Omega$, where $S$ is the set of simple reflections of $W$. The length function $l$ of $W$ and the Bruhat order $\leq$ on $W$ can be extended to $\widetilde{W}$
by setting $l(\omega w) = l(w)$ and $\omega w \leq \omega^{\prime} u$ if and only if $\omega=\omega^{\prime}$ and $w\leq u$,
where $\omega, \omega^{\prime}$ are in $\Omega$ and $w, u$ are in $W$.

Let $\widetilde{\mathcal {H}}$ be the Hecke algebra of $(\widetilde{W}, S)$ over Lurant polynomial ring $\mathcal
{A}=\mathbb{Z}[q^{{\frac{1}{2}}},q^{-{\frac{1}{2}}}]$ ($q$ is an indeterminate). $\widetilde{\mathcal {H}}$ is a
free $\mathcal {A}$-module and has a basis $\{T_{w}\}_{w\in\widetilde{W}}$. Its multiplication is defined by the relations $(T_{s} -
q)(T_{s} + 1) = 0$ if $s\in S$ and $T_{w}T_{u} = T_{wu}$ if $l(wu)=l(w)+l(u)$. $C_{w} = q^{-\frac{l(w)}{2}}\sum_{u\leq w}P_{u, w}T_{u}, w\in\widetilde{W}$ be its Kazhdan-Lusztig basis, where
$P_{u,w}\in \mathbb{Z}[q]$ $(u\leq w)$ are the Kazhdan-Lusztig polynomials(Here the $C_{w}$ is just the
$C^{\prime}_{w}$ in [1].). The degree of $P_{u, w}$ is less than or equal to
$\frac{1}{2}(l(w)-l(u)-1)$ if $u<w$ and $P_{w, w} = 1$. The subalgebra $\mathcal {H}$ of $\widetilde{\mathcal {H}}$
generated by all $T_{s}$ $(s\in S)$ is the Hecke algebra of $(W, S)$. For $\omega$ in $\Omega$
and $w$ in $W$, $C_{\omega w}=T_{\omega}C_{w}$, $P_{\omega u, \omega w}=P_{u, w}$ and $P_{\omega^{\prime} u,\omega w}=0$
for any different $\omega^{\prime},\omega$ in $\Omega$ and $u,w$ in $W$. It is known that the coefficients of these polynomials are all non-negative (see[2]).

The Kazhdan-Lusztig polynomials can be write as $P_{u,w}=\mu(u, w)q^{\frac{1}{2}(l(w)-l(u)-1)}+$ lower degree terms. The coefficient
$\mu(u, w)$ is very interesting, which can be seen even from the recursive formula
for Kazhdan-Lusztig polynomials [1, (2.2.c)]. $\mu(u,w)$ is called the
leading coefficient of $P_{u, w}$. Write $u\prec w$ if $u\leq w$ and $\mu(u, w)\neq0$.
Define $\widetilde{\mu}(u,w) = \mu(u,w)$ if $u\leq w$, $(-1)^{l(w)-l(u)}<0$ and $\widetilde{\mu}(u, w) = \mu(w, u)$ if $w\leq u$.

For $\mu(u, w)$, there are following conclusions (see[1]).

(1) Let $u, w\in W, s\in S$ be such that $u<w, su>u, sw<w$. Then $u\prec w$
if and only if $w=su$. Moreover, $\mu(u, w)=1$ in this case.

(2) Let $u,w\in W, s\in S$ be such that $u<w, us>u, ws<w$. Then $u\prec w$
if and only if $w=us$. Moreover, $\mu(u, w)=1$ in this case.

The elements $C_{w}$ have the following properties (see[1]). For $s\in S$ we have
$$C_{s}C_{w}=\left\{
             \begin{array}{ll}
               (q^{\frac{1}{2}}+q^{-\frac{1}{2}})C_{w}, & \hbox{$sw\leq w$;} \\
               C_{sw}+\sum_{y\prec w, sy\leq y}C_{y},   & \hbox{$sw\geq w$.}
             \end{array}
           \right.$$

$$C_{w}C_{s}=\left\{
             \begin{array}{ll}
               (q^{\frac{1}{2}}+q^{-\frac{1}{2}})C_{w}, & \hbox{$ws\leq w$;} \\
               C_{ws}+\sum_{y\prec w, ys\leq y}C_{y},   & \hbox{$ws\geq w$.}
             \end{array}
           \right.$$

Given $x,y\in W$, $C_{x}C_{y}=\sum_{z\in W}h_{x,y,z}C_{z}$, $h_{x,y,z}\in \mathcal {A}$. For any $z\in W$, define $a(z)$=min\{$i\in \mathbb{N}\mid q^{-\frac{1}{2}}h_{x,y,z}\in \mathbb{Z}[q^{-\frac{1}{2}}]$, for any $x, y\in W$\}. Denote $h_{x,y,z}=\gamma_{x,y,z}q^{\frac{a(z)}{2}}+\delta_{x,y,z}q^{\frac{a(z)-1}{2}}+$ lower degree terms, where $\gamma_{x,y,z}$ and $\delta_{x,y,z}$ are all non-negative.

\subsection{The lowest two-sided cell}

It is known that (see[11,12]) $c_{0}=\{w\in \widetilde{W}\mid a(w)=l(w_{0})\}$ is a two-sided cell, which is the
lowest one with respect to the partial order $\leq_{LR}$, where $w_{0}$ is the longest element of $W_{0}$. $c_{0}$ is called the lowest two-sided cell of $\widetilde{W}$.

In [20], Leonard Scott and Nanhua Xi gave a description of $c_{0}$. Let $R^{+}$ (resp. $R^{+}$, $\Delta$) be the set of positive (resp. negtive, simple) roots in the root system $R$ of $W_{0}$. The set of dominant weights $\Lambda^{+}$ is the set $\{z\in\Lambda\mid l(zw_{0})=l(z)+l(w_{0})\}$. For each simple root $\alpha$, $s_{\alpha}$ is the corresponding simple reflection in $W_{0}$ and $x_{\alpha}$ is the corresponding fundamental weight. For each $w\in W_{0}$, set $d_{w}=w\prod_{\{\alpha\in\Delta,w(\alpha)\in R^{-}\}}x_{\alpha}$. Then $c_{0}=\{d_{w}zw_{0}d_{u}^{-1}\mid w,u\in W_{0}, z\in\Lambda^{+}\}$.

Moreover, the set $c^{\prime}_{0,w}=\{d_{w}zw_{0}d_{u}^{-1}\mid u\in W_{0}, z\in\Lambda^{+}\}$ is a right cell of $\widetilde{W}$ and
$c_{0,u}=\{d_{w}zw_{0}d_{u}^{-1}\mid w\in W_{0}, z\in\Lambda^{+}\}$ is a left cell of $\widetilde{W}$. The distinguished involutions of $c_{0}$ are $d_{w}w_{0}d^{-1}_{w}$, $w\in W_{0}$.

For $x\in \Lambda^{+}$, let $V(x)$ be a rational irreducible $G$-module of highest weight $x$ and $S_{x}$ be the corresponding
element defined in [4]. For any $x\in \Lambda$, Bernstein Joseph gave an element $\theta_{x}\in\mathcal {H}$, defined as $$\theta_{x}=(q^{-l(x^{\prime})}T_{x^{\prime}})(q^{-l(x^{\prime\prime})}T_{x^{\prime\prime}})^{-1},$$
where $x^{\prime},x^{\prime\prime}\in \Lambda^{+}$ satisfies $x=x^{\prime}-x^{\prime\prime}$.For $x\in \Lambda^{+}$, $(x)W_{0}$ is $W_{0}$-orbit in $\Lambda$ which contains $x$. Set $z_{x}=\sum_{x^{\prime}\in(x)W_{0}}\theta_{x^{\prime}}$.
For $x\in \Lambda^{+}$, define
$$S_{x}=\sum_{\{x^{\prime}\in \Lambda^{+},x^{\prime}\leq x\}}d_{x^{\prime}}(x)z_{x^{\prime}},$$
where $d_{x^{\prime}}(x)$ is the multiplicity of $x^{\prime}$ in $V(x)$. Then, $\{S_{x}\mid x\in \Lambda^{+}\}$ is a $\mathcal {A}$-base of the
center of $\mathcal {H}$.

For every $w\in W_{0}$, define $E_{d_{w}}=q^{-\frac{l(d_{w})}{2}}\sum_{\{y\leq d_{w}, l(yw_{0})=l(y)+l(w_{0})\}}P_{yw_{0},d_{w}w_{0}}T_{y}$
and $F_{d_{w}}=q^{-\frac{l(d_{w})}{2}}\sum_{\{y\leq d_{w}, l(yw_{0})=l(y)+l(w_{0})\}}P_{yw_{0},d_{w}w_{0}}T_{y^{-1}}$. Then

(a) $E_{d_{w}}S_{x}C_{w_{0}}F_{d_{u}} = C_{d_{w}xw_{0}d^{-1}_{u}}$ for any $w, u\in W$ and $x\in \Lambda^{+}$;

(b) $S_{x}S_{y} =\sum _{z\in \Lambda^{+}}m_{x, y, z}S_{z}$ for any $x, y\in\Lambda^{+}$. Here $m_{x, y, z}$ is defined to be the multiplicity of $V(z)$ in the tensor product $V (x)\bigotimes V (y)$;

(c) For any $y,w\in c_{0}$, if $\mu(y,w)\neq 0$, then we have $y\thicksim_{L}w$ and $y\nsim_{R}w$, or $y\thicksim_{R}w$ and $y\nsim_{L}w$;

(d) Let $w,w^{\prime},u \in W_{0}$, $y,z\in \Lambda^{+}$ and let $d_{w}yw_{0}d^{-1}_{u} , d_{w^{\prime}}zw_{0}d^{-1}_{u}$ be elements of the left cell $c_{0,u}$. Then $\mu(d_{w}yw_{0}d^{-1}_{u} , d_{w^{\prime}}zw_{0}d^{-1}_{u})=\mu(d_{w}yw_{0}, d_{w^{\prime}}zw_{0})$;

(e) Let $x=d_{u}xw_{0},y=d_{u^{\prime}}x^{\prime}w_{0}$,  where $u, u^{\prime}, v \in W_{0}$ and $x,x^{\prime} \in \Lambda^{+}$. Then the set $\{z\in W\mid h_{w_{0}d_{u}^{-1},d_{u^{\prime}}w_{0},z}\neq 0\}$ is finte, and by $h_{w_{0}d_{u}^{-1},d_{u^{\prime}}w_{0},z}\neq 0$ we know that $z=z_{1}w_{0}$ for some $z_{1}\in \Lambda^{+}$. Moreover, $\mu(y,w)=\sum_{z\in \Lambda^{+}}m_{x^{\ast},x^{\prime},z^{\ast}}\delta_{w_{0}d_{u}^{-1},d_{u^{\prime}}w_{0},zw_{0}}$, where $x^{\ast}=w_{0}x^{-1}w_{0}$.

\subsection{The lowest two-sided cell of $\widetilde{G}_{2}$}

For an affine Weyl group $(W, S)$ of type $\widetilde{G}_{2}$, $W$ is generated by $S=\{r, s, t\}$ and the generated relations are $r^{2}=s^{2}=t^{2}=e$ and $(st)^{6}=(rt)^{2}=(rs)^{3}=e$. $W_{0}=\{e,s,t,st,ts,sts,tst,stst,tsts,ststs,tstst,ststst\}$ and $\Delta=\{\alpha,\beta\}$.

The fundamental weights of $\widetilde{G}_{2}$ are $x_{\alpha}=2\alpha+\beta=ststsr$ and $x_{\beta}=3\alpha+2\beta=(tstsr)^{2}$ in $W$.

As above subsection mentioned, for each $w\in W_{0}$, $d_{w}=w\prod_{\{\alpha\in\Delta,w(\alpha)\in R^{-}\}}x_{\alpha}.$
Then we get the following results.
\begin{eqnarray*}
&&d_{e}=e,\\
&&d_{s}=sx_{\alpha}=s(ststsr)=tstsr,\\
&&d_{t}=tx_{\beta}=s(tstsr)^{2}=ststrstsr,\\
&&d_{st}=stx_{\beta}=st(tstsr)^{2}=tstrstsr,\\
&&d_{ts}=nntsx_{\alpha}=ts(ststsr)=stsr,\\
&&d_{sts}=nnstsx_{\alpha}=sts(ststsr)=tsr,\\
&&d_{tst}=ntstx_{\beta}=tst(tstsr)^{2}=strstsr,\\
&&d_{stst}=nststx_{\beta}=stst(tstsr)^{2}=trstsr,\\
&&d_{tsts}=nntstsx_{\alpha}=tsts(ststsr)=sr,\\
&&d_{ststs}=nststsx_{\alpha}=ststs(ststsr)=r,\\
&&d_{tstst}=tststx_{\beta}=tstst(tstsr)^{2}=rstsr,\\
&&d_{ststst}=stststx_{\alpha}x_{\beta}=stst(ststsr)(tstsr)^{2}=rstsrtstsr.\\
\end{eqnarray*}

From these results and $c_{0}=\{d_{w}zw_{0}d_{u}^{-1}\mid w,u\in W_{0}, z\in\Lambda^{+}\}$, the two-sided cell $c_{0}$ consists of the
following elements.
$$v(i,j,a,b)=v_{i}ststst(x)^{a}(y)^{b}v_{j}^{-1}, 1\leq i,j\leq12, a,b\in\mathbb{N},$$
where
$$v_{1}=e, v_{2}=r, v_{3}=sr, v_{4}=tsr, v_{5}=stsr, v_{6}=tstsr,v_{7}=rstsr,$$
$$v_{8}=trstsr, v_{9}=strstsr, v_{10}=tstrstsr,v_{11}=ststrstsr, v_{12}=rststrstsr,$$
$$x=rststs, y=(rstst)^{2}.$$

  \section{Algorithm for computing $C_{s}C_{w}$}\label{sec3}
In this section, we describe our algorithm for an affine Weyl group $(W, S)$ of type $G_{2}$. Our idea is useful for any affine Wyel group.
\begin{lemma} \label{lemma3.1}
 Given a point $v_{0}\in A_{0}^{+}$ where  $A_{0}^{+} = \{x \in \mathfrak{H}^{*} | \langle x,\check{\alpha}\rangle > 0  ,  \alpha \in \Pi, \langle x, \check{\alpha_{0}}\rangle < 1\}$ and $w_{1}, w_{2}\in W$, then $w_{1}v_{0}=w_{2}v_{0}$ if and only if $w_{1}=w_{2}$, where $w_{1}$ and $w_{1}$ are the element of the affine Weyl group.
\end{lemma}
\beginproof
 If $w_{1}=w_{2}$, it obvious that $w_{1}v_{0}=w_{2}v_{0}$.

 If $w_{1}v_{0}=w_{2}v_{0}$, then $w_{1}v_{0}=w_{2}v_{0}\in w_{1}A^{+}\cap w_{2}A^{+}\neq\emptyset$ and $w_{1}A^{+}=w_{2}A^{+}$. Since $W_{a}$ acts
 on all alcoves simply transitively, we get $w_{1}= w_{2}$.
\endproof

\begin{lemma} \label{lemma3.2}(see[24])
 Let $(W,S)$ be a Coxeter group and $w\in W$. Fix any reduced expression $w=s_{1}\cdots s_{n}$. Then $v\leq w$ if and only if there exist $1\leq j_{1}< j_{2}<\cdots<j_{p}\leq n$ such that $v=s_{1}\cdots \hat{s}_{j_{1}}\cdots \hat{s}_{j_{p}}\cdots s_{n}$. In fact, for $v\leq w$, $j=$($j_{1}< j_{2}<\cdots<j_{p}$) can be choosed such that the above expression is a reduced expression of $v$.
\end{lemma}

\begin{lemma} \label{lemma3.3}(see[26])
Assume that $w=s_{1}s_{2}...s_{q}$ with $s_{i}\in S$. We can find a subsequence $i_{1}<i_{2}<...<i_{r}$ of $1,2,...,r$ such that
$w=s_{i_{1}}s_{i_{2}}...s_{i_{r}}$ is a reduced expression.
\end{lemma}

\begin{lemma} \label{lemma3.4}
Let $x, y\in W$ and $l(x)-l(y)\geq 3$. If there is a $s\in S$, such that $sx<x$ and $sy>y$ (resp.$xs<x$ and $ys>y$), then $\mu(x, y)=0$
\end{lemma}
\beginproof
If there is a $s\in S$, such that $sx<x$ and $sy>y$, then $P_{y, x}=P_{sy, x}$. since $l(x)-l(y)\geq 3$, the degree of $P_{y, x}$ equal to the degree of $P_{sy, x}\leq\frac{l(x)-l(sy)-1}{2}=\frac{l(x)-l(y)-2}{2}<\frac{l(x)-l(y)-1}{2}$, then $\mu(x, y)=0$.
\endproof

\textbf{First Step}. The subwords of $w\in W$ is obtained by following two cases.

\textbf{Case 1}. If the length of $w$ is small ($l(w)\leq 10$), the subwords of $w$ is obtained by the algorithm of binary tree. The set of the results is denoted by $M$. For example, if $w=rstsr$, then $M=\{e, r, s, t, rs, rt, rr, sr, st, ss, ts, rst, rss, rsr, rts, rtr, tr, sts, str, ssr, tsr, rsts,\\ rstr, rssr, rtsr, stsr, rstsr\}$.

\textbf{Case 2}. If the length of $w$ is large enough ($l(w)> 10$),  the subwords of $w$ is obtained by the following procedure.

(i). Split $w$ into two smaller subwords $w_{1}, w_{2}$;

(ii). Get set $A$ (resp. $B$) of all subwords of $w_{1}$ (resp. $w_{2}$) by the algorithm of binary tree;

(iii). Use the generated relations of affine Weyl group $W$ mentioned in subsection II.C to do some preliminary reduce for $A$ and $B$.  More precisely, delete $u$ if $u=v\in A$ (reps. $B$) by the following relations and $|u|\geq|v|$ where $u$ has $|u|$ letters. The results are denoted by $A^{'}$ and $B^{'}$;

(iv). Glue the elements of $A^{'}$ and $B^{'}$. The set of the results is denoted by $M$. That is, $M=\{uv|u\in A^{'}$, $v\in B^{'}\}$.

For example, $w=rtstsrtstsrststst$. If we use the method mentioned in case (1), it takes 20 seconds to get 21764 subwords of $w$. But if $w$ is split into $w_{1}=rtstsrts$ and $w_{2}=tsrststst$, we use the method mentioned in case (2). The 1706 subwords of $w$ are founded less than a second. It is obvious that the method mentioned in case (2) is more efficient than the algorithm of binary tree.

\textbf{Second Step}. Find all the pairwise different reduced subwords of $w$.

(i). For any $m \in M$, $|m|$ is the number of letters of $m$. For example, $m=sts$, then $|m|=3$. The elements of $M$ can be reordered as $w_{1},w_{2},...,w_{r}$, where $|w_{1}|\leq |w_{2}|\leq\cdots\leq |w_{r}|$.

(iii). Fix a rational point $v_{0}$ of $A^{+}_{0}$. Delete $w_{i}$ if there exists $j<i$ such that $w_{i}v_{0}=w_{j}v_{0}$. The the set of the results is denoted by $A$.

The elements of $A$ are the pairwise different reduced subwords of $w$ by the lemma III.2.

\textbf{Third Step}. Compute Kazhdan-Lusztig polynomial $P_{x,y}$ for $x,y\in W$ by following cases.

\textbf{Case I}.  $l(y)<l(x)$, $P_{x,y}=0$.

\textbf{Case II}.  $l(y)=l(x)$, if $yv_{0}=xv_{0}$, then $y=x$, $P_{x,y}=1$. Otherwise, $P_{x,y}=0$.

\textbf{Case III}.  $l(x)<l(y)$, $y=sw$ and $l(y)=l(w)+1$. We use the following recursive formula to compute $P_{x,y}$.

$$P_{x,y}=P_{x,sw}=q^{1-a}P_{sx,w}+q^{a}P_{x,w}+\sum_{x\leq z\prec w,sz<z}q^{\frac{l(w)-l(z)+1}{2}}\mu(z,w)P_{x,z},$$ where a=1 if $sx<x$, a=0 if $sx>x$.

The reduced form of $sx$ is obtained by following procedure. Firstly, find the reordered subwords set $M$ of $sx$ by the second step. Secondly, find the first subword $x^{1}_{i}$ of $sx$ which satisfies $x^{1}_{i}v_{0}=sxv_{0}$. $x^{1}_{i}$ is the reduced form of $sx$ by lemma III.3. If $l(x^{1}_{i})>l(x)$, then $a=0$; otherwise, $a=1$.

All $z$ in the above recursive formula is obtained by following procedure.

(1). The set $A$ of all pairwise different reduced subwords of w is obtained by second step;

(2). Delete the element of $A$ satisfies the condition that $l(w)-l(u)$ is even;

(3). Delete the remaining element which does not satisfy the condition that $su_{1}v_{0}=u^{r}_{1}v_{0}$, where $u^{r}_{1}$ is a reduced subword of $u_{1}$;

(4). Delete the remaining element which does not satisfy the condition that $xv_{0}=u^{r}_{2}v_{0}$, where $u^{r}_{2}$ is a reduced subword of $u_{2}$;

(5). Delete the remaining element which satisfies the conditions of lemma III.4;

(6). Delete the remaining element which does not satisfy the condition that the term $\mu(w_{i},w)q^{\frac{l(w)-l(w_{i})+1}{2}}$ of $P_{w_{i},w}$ is nonzero.

If $P_{sx,w}$, $P_{x,w}$, $P_{z,w}$ and $P_{x,z}$ are known, $P_{x,y}$ are computed by the above recursive formula; otherwise, we compute $P_{x,y}$ through recursive algorithm.

\textbf{Fourth Step}. Compute $C_{s}C_{x}$ for $s\in S, x\in W$ by following formula:
$$C_{sw}=C_{s}C_{w}-\sum_{y\prec w, sy\leq y}\mu(y,w)C_{y}, sw\geq w$$.

  \section{Computation $\mu(y,w)$ for $y,w \in c_{0}$}\label{sec4}

In this section, we use the algorithm mentioned in section III to calculate $C_{s}C_{x}$ for some $s\in S, x\in W$. These multiplications play important role in computing $C_{w_{0}d_{u}^{-1}}C_{d_{u^{'}}w_{0}}$ for $u, u^{'}\in W_{0}$ in order to get $\delta_{w_{0}d_{u}^{-1},d_{u^{'}}w_{0},z_{1}w_{0}}$ which is the coefficient of the term $q^{\frac{1}{2}(a(z_{1}w_{0})-1)}$ of $h_{w_{0}d_{u}^{-1},d_{u^{'}}w_{0},z_{1}w_{0}}$. By II.B(e) in section II, $\mu(y,w)=\sum_{z\in \Lambda^{+}}
m_{x^{\ast},x^{\prime},z^{\ast}}\delta_{w_{0}d_{u}^{-1},d_{u^{\prime}}w_{0},zw_{0}}$, where $x^{\ast}=w_{0}x^{-1}w_{0}$. Our main result theorem IV.1 is obtained by lemma IV.1, proposition IV.1 and lemma IV.2.

We use an algebra antiautomorphism of $\mathcal {H}$ which is defined by George Lusztig in [3] to get the following lemma. The amount of calculation for our problem is reduced by this lemma.

\begin{lemma} \label{lemma4.1}
Let $x=d_{u}xw_{0},y=d_{u^{\prime}}x^{\prime}w_{0}$,  where $u, u^{\prime}, v \in W_{0}$ and $x,x^{\prime} \in \Lambda^{+}$. Then
$C_{x}C_{y}=C_{y^{-1}}C_{x^{-1}}$.
\end{lemma}
\beginproof
There is a unique algebra antiautomorphism of $H$ which is denoted by $\varphi$ such that $\varphi(T_{x})=T_{x^{-1}}$, for all $x \in
W$; it maps $C_{x}\longrightarrow C_{x^{-1}}$. For any $x,y \in W$, $\varphi(C_{x}C_{y})=\varphi(C_{y})\varphi(C_{x})=C_{y^{-1}}C_{x^{-1}}$; That is,
$\varphi(C_{x}C_{y})=\varphi(\sum_{z\in W}h_{x,y,z}C_{z})=\sum_{z\in W}h_{x,y,z}C_{z^{-1}}=\sum_{z\in W}h_{y^{-1},x^{-1},z}C_{z}$, we get
$h_{x,y,z}=h_{y^{-1},x^{-1},z^{-1}}$. By 2.2(d), we know that $h_{x,y,z}\neq 0$, $z=z_{1}w_{0}$for some $z_{1}\in \Lambda^{+}$. Then if $h_{x,y,z}\neq
0$, $h_{x,y,z}=h_{x,y,z_{1}w_{0}}=h_{y^{-1},x^{-1},z^{-1}}=h_{y^{-1},x^{-1},w_{0}z^{-1}}$. In the case of $\widetilde{G_{2}}$, we know that
$z_{1}w_{0}=w_{0}z^{-1}$; that is, $h_{x,y,z}=h_{y^{-1},x^{-1},z}$ for all $z\in W$.
\endproof

In the following, all $C_{w_{0}d_{u}^{-1}}C_{d_{u^{'}}w_{0}}$ for $u, u^{'}\in W_{0}$ are computed to obtain $\delta_{w_{0}d_{u}^{-1},d_{u^{'}}w_{0},z_{1}w_{0}}$ which is the coefficient of the term $q^{\frac{1}{2}(a(z_{1}w_{0})-1)}$ of $h_{w_{0}d_{u}^{-1},d_{u^{'}}w_{0},z_{1}w_{0}}$.

\begin{proposition} \label{Prop4.1}
For any $u, u^{'}\in W_{0}$, if $\delta_{w_{0}d_{u}^{-1},d_{u^{'}}w_{0},z_{1}w_{0}}\neq0$, we must have
$\delta_{w_{0}d_{u}^{-1},d_{u^{'}}w_{0},z_{1}w_{0}}=1$ and $z_{1}\in\{0,x_{\alpha},x_{\beta}\}$, where $x_{\alpha},x_{\beta}\in\Lambda^{+}$ are
fundamental weights.
\end{proposition}
\beginproof
We compute the $C_{x}C_{y}$ through the following steps.

(a) We use the formula $C_{w}C_{s}=C_{x}+\sum_{u\prec w, us<u}\mu(u, w)C_{u}$ for $x=ws$ and $l(x)=l(w)+1$ (see[1]) to get a decomposition of $C_{x}$;

(b) For the $s$ in step (a), we compute $C_{s}C_{y}$ by the formula $C_{s}C_{y}=C_{sy}+\sum_{v\prec y, sv<y}\mu(v, y)C_{v}$ for $sy>y$ and $C_{s}C_{y}=(q^{\frac{1}{2}}+q^{-\frac{1}{2}})C_{y}$ for $sy<y$ (see[1]);

(c) We get $C_{x}C_{y}=\sum_{u, v\in W}C_{u}C_{v}$ where $l(u)<l(x)$ by step (a) and step (b). Repeat this proceeding until $C_{x}C_{y}=\sum_{u^{\prime},v^{\prime}\in W}C_{u^{\prime}}C_{v^{\prime}}$ where $l(u^{\prime})=0$, then we obtain $C_{x}C_{y}$.

It is known that $P_{w,v}=1$ for any $w, v\in W_{0}$. By this conclusion and above proceeding, the complete decomposition $C_{w_{0}}=C_{1}C_{2}C_{1}C_{2}C_{1}C_{2}-4C_{1}C_{2}C_{1}C_{2}+3C_{1}C_{2}$ is obtained. This complete decomposition plays a crucial role in the proof of the following proposition. $C_{w_{0}d_{u}^{-1}}C_{d_{u^{'}}w_{0}}$ for $u, u^{'}\in W_{0}$ will be calculated with the increase of $l(w_{0}d_{u}^{-1})$ by following cases.

(i) If $C_{w_{0}d_{u}^{-1}}=C_{w_{0}}$, we get

$$C_{w_{0}d_{u}^{-1}}C_{d_{u^{'}}w_{0}}=C_{w_{0}}C_{d_{u^{'}}w_{0}}=(C_{1}C_{2}C_{1}C_{2}C_{1}C_{2}-4C_{1}C_{2}C_{1}C_{2}+3C_{1}C_{2})C_{d_{u^{'}}w_{0}}.$$
The algorithm is applied to compute $C_{1}C_{2}C_{1}C_{2}C_{1}C_{2}C_{d_{u^{'}}w_{0}}, C_{1}C_{2}C_{1}C_{2}C_{d_{u^{'}}w_{0}},C_{1}C_{2}C_{d_{u^{'}}w_{0}}$. Then the decomposition of $C_{w_{0}}C_{d_{u^{'}}w_{0}}$ is obtained;

(ii) If $l(w_{0}d_{u}^{-1})>l(w_{0})$, there are two types of calculation of $C_{w_{0}d_{u}^{-1}}C_{d_{u^{'}}w_{0}}$.

(1) Step (a) and step (b) mentioned above will be repeatedly used to get $C_{w_{0}d_{u}^{-1}}C_{d_{u^{'}}w_{0}}=\sum_{u,v\in W}C_{u}C_{v}$, where $C_{u}C_{v}$ can be computed by recursion from case (i);

(2) If $C_{w_{0}d_{u}^{-1}}C_{d_{u^{'}}w_{0}}$ can not be computed by (1), step (a) and step (b) mentioned above will be repeatedly used to get $C_{w_{0}d_{u}^{-1}}C_{d_{u^{'}}w_{0}}=\sum_{u,v\in W}C_{u}C_{v}$, where $C_{u}=C_{w_{0}}$. All $C_{w_{0}}C_{v}$ can be compute by (i).

Set $[2]\triangleq q^{\frac{1}{2}}+q^{-\frac{1}{2}}$, $1\triangleq s$, $2\triangleq t$ and $0\triangleq r$.

(1) Computing $C_{w_{0}}C_{w_{0}}$

\begin{eqnarray*}
C_{12}&=&C_{1}C_{2},\\
C_{2}C_{12}&=&C_{121}+C_{2},\\
C_{121}&=&C_{2}C_{1}C_{2}-C_{2},\\
C_{2}C_{121}&=&C_{1212}+C_{12},\\
C_{1212}&=&C_{1}C_{2}C_{1}C_{2}-2C_{1}C_{2},\\
C_{2}C_{1212}&=&C_{21212}+C_{212},\\
C_{21212}&=&C_{2}C_{1}C_{2}C_{1}C_{2}-3C_{2}C_{1}C_{2}+C_{2},\\
C_{1}C_{21212}&=&C_{w_{0}}+C_{1212},\\
C_{w_{0}}&=&C_{1}C_{2}C_{1}C_{2}C_{1}C_{2}-4C_{1}C_{2}C_{1}C_{2}+3C_{1}C_{2}.\\
\end{eqnarray*}

Since $C_{s}C_{w}=(q^{\frac{1}{2}}+q^{-\frac{1}{2}})C_{w}$ if $sw<w$, then

\begin{eqnarray*}
C_{w_{0}}C_{w_{0}}&=&(C_{1}C_{2}C_{1}C_{2}C_{1}C_{2}-4C_{1}C_{2}C_{1}C_{2}+3C_{1}C_{2})C_{w_{0}}\\
&=&[(q^{\frac{1}{2}}+q^{-\frac{1}{2}})^{6}-4(q^{\frac{1}{2}}+q^{-\frac{1}{2}})^{4}+3(q^{\frac{1}{2}}+q^{-\frac{1}{2}})^{2}]C_{w_{0}}\\
&=&([2]^{6}-4[2]^{4}+3[2]^{2})C_{w_{0}}.\\
\end{eqnarray*}

(2) Computing $C_{w_{0}}C_{0121212}$
\begin{eqnarray*}
C_{2}C_{0121212}&=&[2]C_{0121212},\\
C_{1}C_{0121212}&=&C_{10121212}+C_{w_{0}}, C_{1}C_{2}C_{0121212}=[2](C_{10121212}+C_{w_{0}}),\\
C_{2}C_{10121212}&=&C_{210121212}+C_{0121212},\\
C_{2}C_{1}C_{2}C_{0121212}&=&[2][C_{210121212}+C_{0121212}+[2]C_{w_{0}}],\\
C_{1}C_{210121212}&=&C_{1210121212}+C_{10121212},\\
C_{1}C_{2}C_{1}C_{2}C_{0121212}&=&[2]\{C_{1210121212}+2C_{10121212}+[[2]^{2}+1]C_{w_{0}}\},\\
C_{2}C_{1210121212}&=&C_{21210121212}+C_{210121212},\\
C_{2}C_{1}C_{2}C_{1}C_{2}C_{0121212}&=&[2]\{C_{21210121212}+3C_{210121212}+2C_{0121212}\\
& &+[[2]^{3}+[2]]C_{w_{0}}\},\\
C_{1}C_{21210121212}&=&C_{121210121212}+C_{1210121212},\\
C_{1}C_{2}C_{1}C_{2}C_{1}C_{2}C_{0121212}&=&[2]\{C_{121210121212}+4C_{1210121212}+5C_{10121212}\\
& &+([2]^{4}+[2]^{2}+2)C_{121212}\}.\\
\end{eqnarray*}
\begin{eqnarray*}
C_{121212}C_{0121212}&=&(C_{1}C_{2}C_{1}C_{2}C_{1}C_{2}-4C_{1}C_{2}C_{1}C_{2}+3C_{1}C_{2})C_{0121212}\\
&=&[2]\{C_{121210121212}+4C_{1210121212}+5C_{10121212}+([2]^{4}+[2]^{2}+2)C_{121212}\}\\
& &-4[2]\{C_{1210121212}+2C_{10121212}+([2]^{2}+1)C_{121212}\}\\
& &+3[[2](C_{10121212}+C_{121212})]\\
&=&[2]C_{121210121212}+([2]^{5}-3[2]^{2}+[2])C_{121212}\\
&=&[2]C_{x_{\alpha}w_{0}}+([2]^{5}-3[2]^{3}+[2])C_{w_{0}}.\\
\end{eqnarray*}

(3) Computing $C_{121212}C_{10121212}$
\begin{eqnarray*}
C_{2}C_{10121212}&=&C_{210121212}+C_{0121212},\\
C_{1}C_{210121212}&=&C_{1210121212}+C_{10121212},\\
C_{1}C_{2}C_{10121212}&=&C_{1210121212}+2C_{10121212}+C_{121212},\\
C_{2}C_{1210121212}&=&C_{21210121212}+C_{210121212},\\
C_{2}C_{1}C_{2}C_{10121212}&=&C_{21210121212}+3C_{210121212}+2C_{0121212}+[2]C_{121212},\\
C_{1}C_{21210121212}&=&C_{121210121212}+C_{1210121212},\\
C_{1}C_{2}C_{1}C_{2}C_{10121212}&=&C_{121210121212}+4C_{1210121212}+5C_{10121212}+([2]+2)C_{121212},\\
C_{2}C_{121210121212}&=&[2]C_{121210121212},\\
C_{2}C_{1}C_{2}C_{1}C_{2}C_{10121212}&=&[2]C_{121210121212}+4C_{21210121212}+9C_{210121212}\\
& &+5C_{0121212}+([2]^{3}+2[2]+2)C_{121212},\\
C_{1}C_{2}C_{1}C_{2}C_{1}C_{2}C_{10121212}&=&([2]^{2}+4)C_{121210121212}+13C_{1210121212}+14C_{10121212}\\
& &+([2]^{4}+2[2]^{2}+5)C_{121212},\\
\end{eqnarray*}
\begin{eqnarray*}
C_{121212}C_{10121212}&=&(C_{1}C_{2}C_{1}C_{2}C_{1}C_{2}-4C_{1}C_{2}C_{1}C_{2}+3C_{1}C_{2})C_{10121212}\\
&=&([2]^{2}+4)C_{121210121212}+13C_{1210121212}+14C_{10121212}\\
& &+([2]^{4}+2[2]^{2}+5)C_{121212}-4\{C_{121210121212}+4C_{1210121212}\\
& &+5C_{10121212}+([2]+2)C_{121212}\}+3[C_{1210121212}\\
& &+2C_{10121212}+C_{121212}]\\
&=&[2]^{2}C_{121210121212}+([2]^{4}-2[2]^{2})C_{121212}\\
&=&[2]^{2}C_{x_{\alpha}w_{0}}+([2]^{4}-2[2]^{2})C_{w_{0}}.\\
\end{eqnarray*}

(4) Computing $C_{121212}C_{210121212}$
\begin{eqnarray*}
C_{2}C_{210121212}&=&[2]C_{210121212}, \\
C_{1}C_{210121212}&=&C_{1210121212}+C_{10121212},\\
C_{1}C_{2}C_{210121212}&=&[2](C_{1210121212}+C_{10121212}),\\
C_{2}C_{1210121212}&=&C_{21210121212}+C_{210121212},\\
C_{2}C_{1}C_{2}C_{210121212}&=&[2](C_{21210121212}+2C_{210121212}+C_{0121212}),\\
C_{1}C_{21210121212}&=&C_{121210121212}+C_{1210121212},\\
C_{1}C_{2}C_{1}C_{2}C_{210121212}&=&[2](C_{121210121212}+3C_{1210121212}+3C_{10121212}+C_{121212}),\\
C_{2}C_{121210121212}&=&[2]C_{121210121212},\\
C_{1}C_{121210121212}&=&[2]C_{121210121212},\\
C_{2}C_{1}C_{2}C_{1}C_{2}C_{210121212}&=&[2]([2]C_{121210121212}+3C_{21210121212}+6C_{210121212}\\
& &+3C_{0121212}+[2]C_{121212}),\\
C_{1}C_{2}C_{1}C_{2}C_{1}C_{2}C_{210121212}&=&[2]\{([2]^{2}+3)C_{21210121212}+9C_{1210121212}\\
& &+9C_{10121212}+([2]^{2}+3)C_{121212}\}.\\
\end{eqnarray*}
\begin{eqnarray*}
C_{121212}C_{210121212}&=&(C_{1}C_{2}C_{1}C_{2}C_{1}C_{2}-4C_{1}C_{2}C_{1}C_{2}+3C_{1}C_{2})C_{210121212}\\
&=&[2]\{([2]^{2}+3)C_{21210121212}+9C_{1210121212}+9C_{10121212}\\
& &+([2]^{2}+3)C_{121212}\}-4[2]\{C_{121210121212}+3C_{1210121212}\\
& &+3C_{10121212}+C_{121212}\}+3[2][C_{1210121212}+C_{10121212}]\\
&=&([2]^{3}-[2])C_{121210121212}+([2]^{3}-[2])C_{121212}\\
&=&([2]^{3}-[2])C_{x_{\alpha}w_{0}}+([2]^{3}-[2])C_{w_{0}}.\\
\end{eqnarray*}

(5) Computing $C_{121212}C_{1210121212}$
\begin{eqnarray*}
C_{2}C_{1210121212}&=&C_{21210121212}+C_{210121212},\\
C_{1}C_{21210121212}&=&C_{121210121212}+C_{1210121212},\\
C_{1}C_{2}C_{1210121212}&=&C_{121210121212}+2C_{1210121212}+C_{10121212},\\
C_{2}C_{121210121212}&=&[2]C_{121210121212},\\
C_{2}C_{1}C_{2}C_{1210121212}&=&[2]C_{121210121212}+2C_{21210121212}+3C_{210121212}+C_{0121212},\\
C_{1}C_{121210121212}&=&[2]C_{121210121212},\\
C_{1}C_{2}C_{1}C_{2}C_{1210121212}&=&([2]^{2}+2)C_{121210121212}+5C_{1210121212}+4C_{10121212}+C_{121212},\\
C_{2}C_{1}C_{2}C_{1}C_{2}C_{1210121212}&=&([2]^{3}+2[2])C_{121210121212}+5C_{21210121212}\\
& &+9C_{210121212}+4C_{0121212}+[2]C_{121212},\\
C_{1}C_{2}C_{1}C_{2}C_{1}C_{2}C_{1210121212}&=&([2]_{4}+2[2]^{2}+5)C_{121210121212}+14C_{1210121212}\\
& &+13C_{10121212}+(4+[2]^{2})C_{121212}.\\
\end{eqnarray*}
\begin{eqnarray*}
C_{121212}C_{1210121212}&=&(C_{1}C_{2}C_{1}C_{2}C_{1}C_{2}-4C_{1}C_{2}C_{1}C_{2}+3C_{1}C_{2})C_{1210121212}\\
&=&([2]_{4}+2[2]^{2}+5)C_{121210121212}+14C_{1210121212}\\
& &+13C_{10121212}+(4+[2]^{2})C_{121212}\\
& &-4\{([2]_{2}+2)C_{121210121212}+5C_{1210121212}+4C_{10121212}+C_{121212}\}\\
& &+3[C_{121210121212}+2C_{1210121212}+C_{10121212}]\\
&=&([2]^{4}-2[2]^{2})C_{121210121212}+[2]^{2}C_{121212}\\
&=&([2]^{4}-2[2]^{2})C_{x_{\alpha}w_{0}}+[2]^{2}C_{w_{0}}.\\
\end{eqnarray*}

(6) Computing $C_{121212}C_{21210121212}$
\begin{eqnarray*}
C_{2}C_{21210121212}&=&[2]C_{21210121212},\\
C_{1}C_{2}C_{21210121212}&=&[2](C_{121210121212}+C_{1210121212}),\\
C_{2}C_{1}C_{2}C_{21210121212}&=&[2]([2]C_{121210121212}+C_{21210121212}+C_{210121212}),\\
C_{1}C_{2}C_{1}C_{2}C_{21210121212}&=&[2][([2]^{2}+1)C_{121210121212}+2C_{1210121212}+C_{10121212}],\\
C_{2}C_{1}C_{2}C_{1}C_{2}C_{21210121212}&=&[2][([2]^{3}+[2])C_{121210121212}+2C_{21210121212}\\
& &+3C_{210121212}+C_{0121212}],\\
C_{1}C_{2}C_{1}C_{2}C_{1}C_{2}C_{21210121212}&=&([2]^{5}+[2]^{3}+2[2])C_{121210121212}+5[2]C_{1210121212}\\
& &+4[2]C_{10121212}+[2]C_{0121212}).\\
\end{eqnarray*}
\begin{eqnarray*}
C_{121212}C_{21210121212}&=&(C_{1}C_{2}C_{1}C_{2}C_{1}C_{2}-4C_{1}C_{2}C_{1}C_{2}+3C_{1}C_{2})C_{21210121212}\\
&=&([2]^{5}+[2]^{3}+2[2])C_{121210121212}+5[2]C_{1210121212}\\
& &+4[2]C_{10121212}+[2]C_{0121212})-4\{[2][([2]^{2}+1)C_{121210121212}\\
& &+2C_{1210121212}+C_{10121212}]\}+3[[2](C_{121210121212}+C_{1210121212})]\\
&=&([2]^{5}-3[2]^{3}+[2])C_{121210121212}+[2]C_{121212}\\
&=&([2]^{5}-3[2]^{3}+[2])C_{x_{\alpha}w_{0}}+[2]C_{w_{0}}.\\
\end{eqnarray*}

(7) Computing $C_{121212}C_{01210121212}$
\begin{eqnarray*}
C_{121212}&=&C_{2}C_{1}C_{2}C_{1}C_{2}C_{1}-4C_{2}C_{1}C_{2}C_{1}+3C_{2}C_{1},\\
C_{1}C_{01210121212}&=&[2]C_{01210121212},\\
C_{2}C_{01210121212}&=&C_{201210121212},\\
C_{2}C_{1}C_{01210121212}&=&[2]C_{201210121212},\\
C_{1}C_{201210121212}&=&C_{1201210121212}+C_{01210121212},\\
C_{1}C_{2}C_{1}C_{01210121212}&=&[2](C_{1201210121212}+C_{01210121212}),\\
C_{2}C_{1201210121212}&=&C_{21201210121212}+C_{201210121212}+C_{x_{\alpha}w_{0}},\\
C_{2}C_{1}C_{2}C_{1}C_{01210121212}&=&[2](C_{21201210121212}+2C_{201210121212}+C_{x_{\alpha}w_{0}}),\\
C_{1}C_{21201210121212}&=&C_{121201210121212}+C_{1201210121212},\\
C_{1}C_{2}C_{1}C_{2}C_{1}C_{01210121212}&=&[2](C_{121201210121212}+3C_{1201210121212}\\
& &+2C_{01210121212}+[2]C_{121210121212}),\\
C_{2}C_{121201210121212}&=&C_{2121201210121212}+C_{21201210121212}+C_{w_{0}},\\
C_{2}C_{1}C_{2}C_{1}C_{2}C_{1}C_{01210121212}&=&[2][C_{2121201210121212}+4C_{21201210121212}\\
& &5C_{201210121212}+(3+[2]^{2})C_{x_{\alpha}w_{0}}+C_{w_{0}}].\\
\end{eqnarray*}
\begin{eqnarray*}
C_{121212}C_{01210121212}&=&(C_{2}C_{1}C_{2}C_{1}C_{2}C_{1}-4C_{2}C_{1}C_{2}C_{1}+3C_{2}C_{1})C_{01210121212}\\
&=&[2][C_{x_{\beta}w_{0}}+4C_{21201210121212}+5C_{201210121212}\\
& &+(3+[2]^{2})C_{x_{\alpha}w_{0}}+C_{w_{0}}]-4\{[2](C_{21201210121212}\\
& &+2C_{201210121212}+C_{x_{\alpha}w_{0}})\}+3([2]C_{201210121212})\\
&=&[2]C_{2121021210121212}+([2]^{3}-[2])C_{x_{\alpha}w_{0}}+[2]C_{w_{0}}\\
&=&[2]C_{x_{\beta}w_{0}}+([2]^{3}-[2])C_{x_{\alpha}w_{0}}+[2]C_{w_{0}}.\\
\end{eqnarray*}

(8) Computing $C_{121212}C_{201210121212}$
\begin{eqnarray*}
C_{2}C_{01210121212}&=&C_{201210121212},\\
C_{121212}C_{01210121212}&=&C_{121212}C_{2}C_{01210121212}=[2]C_{121212}C_{01210121212}\\
&=&[2][[2]C_{2121021210121212}+([2]^{3}-[2])C_{121210121212}]\\
&=&[2]^{2}C_{x_{\beta}w_{0}}+([2]^{4}-[2]^{2})C_{x_{\alpha}121212}+[2]^{2}C_{w_{0}}.\\
\end{eqnarray*}

(9) Computing $C_{121212}C_{1201210121212}$
\begin{eqnarray*}
C_{1}C_{201210121212}&=&C_{1201210121212}+C_{01210121212},\\
C_{1201210121212}&=&C_{1}C_{201210121212}-C_{01210121212},\\
C_{121212}C_{1201210121212}&=&C_{121212}(C_{1}C_{201210121212}-C_{01210121212})\\
&=&[2]C_{121212}C_{201210121212}-C_{121212}C_{01210121212}\\
&=&[2][[2]^{2}C_{x_{\beta}121212}+([2]^{4}-[2]^{2})C_{121212}+[2]^{2}C_{w_{0}}]\\
& &-([2]C_{x_{\beta}121212}+([2]^{3}-[2])C_{121212}+[2]C_{w_{0}})\\
&=&([2]^{3}-[2])C_{x_{\beta}w_{0}}+([2]^{5}-2[2]^{3}+[2])C_{x_{\alpha}w_{0}}+([2]^{3}-[2])C_{w_{0}}.\\
\end{eqnarray*}

(10) Computing $C_{121212}C_{21201210121212}$
\begin{eqnarray*}
C_{2}C_{1201210121212}&=&C_{21201210121212}+C_{201210121212}+C_{121210121212},\\
C_{21201210121212}&=&C_{2}C_{1201210121212}-C_{201210121212}-C_{121210121212},\\
C_{121212}C_{21201210121212}&=&C_{121212}(C_{2}C_{1201210121212}-C_{201210121212}-C_{121210121212})\\
&=&[2]C_{121212}C_{1201210121212}-C_{121212}C_{201210121212}\\
& &-C_{121212}C_{121210121212}\\
&=&[2][[2]^{3}-[2])C_{x_{\beta}121212}+([2]^{5}-2[2]^{3}+[2])C_{x_{\alpha}121212}\\
& &+([2]^{3}-[2])C_{w_{0}}]-[[2]^{2}C_{x_{\beta}121212}+([2]^{4}-[2]^{2})C_{x_{\alpha}121212}\\
& &+[2]^{2}C_{w_{0}}]-[([2]^{6}-4[2]^{4}+3[2]^{2})C_{x_{\alpha}12121}]\\
&=&([2]^{4}-2[2]^{2})C_{x_{\beta}w_{0}}+([2]^{4}-[2]^{2})C_{x_{\alpha}w_{0}}+([2]^{4}-2[2]^{2})C_{w_{0}}.\\
\end{eqnarray*}

(11) Computing $C_{121212}C_{121201210121212}$
\begin{eqnarray*}
C_{1}C_{21201210121212}&=&C_{121201210121212}+C_{1201210121212},\\
C_{121201210121212}&=&C_{1}C_{21201210121212}-C_{1201210121212},\\
C_{121212}C_{121201210121212}&=&C_{121212}(C_{1}C_{21201210121212}-C_{1201210121212})\\
&=&[2]C_{121212}C_{21201210121212}-C_{121212}C_{1201210121212}\\
&=&[2][([2]^{4}-2[2]^{2})C_{x_{\beta}w_{0}}+([2]^{4}-[2]^{2})C_{x_{\alpha}w_{0}}+([2]^{4}-2[2]^{2})C_{w_{0}}]\\
& &-[([2]^{3}-[2])C_{x_{\beta}w_{0}}+([2]^{5}-2[2]^{3}+[2])C_{x_{\alpha}w_{0}}+([2]^{3}-[2])C_{w_{0}}]\\
&=&([2]^{5}-3[2]^{3}+[2])C_{x_{\beta}w_{0}}+([2]^{3}-[2])C_{x_{\alpha}w_{0}}+([2]^{5}-3[2]^{3}+[2])C_{w_{0}}.\\
\end{eqnarray*}

(12) Computing $C_{121212}C_{0121201210121212}$
\begin{eqnarray*}
C_{2}C_{0121201210121212}&=&C_{20121201210121212}+C_{021201210121212}+C_{0121212}\\
&=&C_{02121201210121212}+C_{201021210121212}+C_{0121212}\\
&=&C_{02121201210121212}+C_{210121210121212}+C_{0121212},\\
C_{1}C_{02121201210121212}&=&C_{102121201210121212}+C_{x_{\beta}w_{0}},\\
C_{1}C_{210121210121212}&=&C_{1210121210121212}+C_{10121210121212},\\
C_{1}C_{0121212}&=&C_{10121212}+C_{w_{0}},\\
C_{2}C_{102121201210121212}&=&C_{2102121201210121212}+C_{02121201210121212},\\
C_{2}C_{10121212}&=&C_{210121212}+C_{0121212},\\
C_{2}C_{1210121210121212}&=&C_{21210121210121212}+C_{210121210121212},\\
C_{2}C_{10121210121212}&=&C_{210121210121212}+C_{0121210121212},\\
C_{1}C_{2102121201210121212}&=&C_{12102121201210121212}+C_{102121201210121212},\\
C_{1}C_{21210121210121212}&=&C_{x_{\alpha}^{2}w_{0}}+C_{w_{0}}+C_{x_{\alpha}w_{0}}+C_{x_{\beta}w_{0}}+C_{1210121210121212},\\
C_{1}C_{0121210121212}&=&C_{10121210121212}+C_{121210121212}=C_{10121210121212}+C_{x_{\alpha}w_{0}},\\
C_{2}C_{12102121201210121212}&=&C_{212102121201210121212}+C_{2102121201210121212},\\
C_{1}C_{212102121201210121212}&=&C_{x_{\alpha}x_{\beta}121212}+C_{12102121201210121212}+C_{x_{\alpha}w_{0}}+C_{x_{\alpha}^{2}w_{0}}.\\
\end{eqnarray*}
\begin{eqnarray*}
C_{1}C_{2}C_{0121201210121212}&=&C_{102121201210121212}+C_{x_{\beta}w_{0}}+C_{10121212}\\
& &+C_{w_{0}}+C_{1210121210121212}+C_{10121210121212},\\
C_{2}C_{1}C_{2}C_{0121201210121212}&=&C_{2102121201210121212}+C_{02121201210121212}+[2]C_{x_{\beta}w_{0}}\\
& &+C_{210121212}+C_{0121212}+[2]C_{w_{0}},\\
& &+C_{21210121210121212}+2C_{210121210121212}+C_{0121210121212}\\
C_{1}C_{2}C_{1}C_{2}C_{0121201210121212}&=&C_{12102121201210121212}+2C_{102121201210121212}+C_{x_{\alpha}^{2}121212}\\
& &+3C_{1210121210121212}+3C_{10121210121212}+C_{1210121212}\\
& &+2C_{10121212}+2C_{x_{\alpha}w_{0}}+([2]^{2}+2)C_{x_{\beta}w_{0}}+([2]^2+2)C_{w_{0}},\\
C_{2}C_{1}C_{2}C_{1}C_{2}C_{0121201210121212}&=&C_{212102121201210121212}+3C_{2102121201210121212}+C_{21210121212}\\
& &+3C_{21210121210121212}+6C_{210121210121212}+3C_{0121210121212}\\
& &+2C_{02121201210121212}+3C_{210121212}+[2]C_{x_{\alpha}^{2}w_{0}}+2[2]C_{x_{\alpha}w_{0}}\\
& &+2C_{0121212}+([2]^{3}+2[2])C_{x_{\beta}w_{0}}+([2]^3+2[2])C_{w_{0}},\\
C_{1}C_{2}C_{1}C_{2}C_{1}C_{2}C_{0121201210121212}&=&C_{x_{\alpha}x_{\beta}121212}+4C_{12102121201210121212}+5C_{102121201210121212}\\
& &+([2]^{2}+4)C_{x_{\alpha}^{2}121212}+9C_{1210121210121212}+9C_{10121210121212}\\
& &+4C_{1210121212}+5C_{10121212}+([2]^{4}+2[2]^{2}+5)C_{w_{0}}\\
& &+(2[2]^{2}+8)C_{x_{\alpha}121212}+([2]^{4}+[2]^{2}+5)C_{x_{\beta}121212}.\\
C_{121212}C_{0121201210121212}&=&(C_{2}C_{1}C_{1}C_{2}C_{1}C_{2}-4C_{2}C_{1}C_{2}C_{1}+3C_{2}C_{1})C_{01210121212}\\
&=&C_{x_{\alpha}x_{\beta}121212}+4C_{12102121201210121212}+5C_{102121201210121212}\\
& &+([2]^{2}+4)C_{x_{\alpha}^{2}121212}+9C_{1210121210121212}+9C_{10121210121212}\\
& &+4C_{1210121212}+5C_{10121212}+([2]^{4}+2[2]^{2}+5)C_{w_{0}}\\
& &+(2[2]^{2}+8)C_{x_{\alpha}121212}+([2]^{4}+[2]^{2}+5)C_{x_{\beta}121212}.\\
& &-4\{C_{12102121201210121212}+2C_{102121201210121212}+C_{x_{\alpha}^{2}121212}\\
& &+3C_{1210121210121212}+3C_{10121210121212}+C_{1210121212}\\
& &+2C_{10121212}+2C_{x_{\alpha}w_{0}}+([2]^{2}+2)C_{x_{\beta}w_{0}}+([2]^2+2)C_{w_{0}}\}\\
& &+3[C_{102121201210121212}+C_{x_{\beta}w_{0}}+C_{10121212}+C_{w_{0}}\\
& &+C_{1210121210121212}+C_{10121210121212}]\\
&=&C_{x_{\alpha}x_{\beta}w_{0}}+[2]^{2}C_{x_{\alpha}^{2}w_{0}}+(2[2]^{2})C_{x_{\alpha}w_{0}}\\
& &+([2]^{4}-2[2]^{2})C_{x_{\beta}w_{0}}+([2]^{4}-2[2]^{2})C_{w_{0}}.\\
\end{eqnarray*}

(12*) $C_{121212}C_{02121021210121212}$
\begin{eqnarray*}
C_{2}C_{02121021210121212}&=&[2]C_{02121021210121212},\\
C_{1}C_{02121021210121212}&=&C_{102121021210121212}+C_{x_{\beta}w_{0}},\\
C_{1}C_{2}C_{02121021210121212}&=&[2](C_{102121021210121212}+C_{x_{\beta}w_{0}}),\\
C_{2}C_{102121021210121212}&=&C_{2102121021210121212}+C_{02121021210121212},\\
C_{2}C_{1}C_{2}C_{02121021210121212}&=&[2](C_{2102121021210121212}+C_{02121021210121212}+[2]C_{x_{\beta}w_{0}}),\\
C_{1}C_{2102121021210121212}&=&C_{12102121021210121212}+C_{102121021210121212},\\
C_{1}C_{2}C_{1}C_{2}C_{02121021210121212}&=&[2][C_{12102121021210121212}+2C_{102121021210121212}\\
& &+([2]^2+1)C_{x_{\beta}w_{0}}],\\
C_{2}C_{12102121021210121212}&=&C_{212102121021210121212}+C_{2102121021210121212},\\
C_{2}C_{1}C_{2}C_{1}C_{2}C_{02121021210121212}&=&[2][C_{212102121021210121212}+3C_{2102121021210121212}\\
& &+C_{02121021210121212}+([2]^3+[2])C_{x_{\beta}w_{0}}],\\
C_{1}C_{212102121021210121212}&=&C_{x_{\alpha}x_{\beta}w_{0}}+C_{x_{\alpha}w_{0}}+C_{x_{\alpha}^{2}w_{0}}+C_{12102121021210121212},\\
C_{1}C_{2}C_{1}C_{2}C_{1}C_{2}C_{02121021210121212}&=&[2]\{C_{x_{\alpha}x_{\beta}w_{0}}+C_{x_{\alpha}w_{0}}+C_{x_{\alpha}^{2}w_{0}}+4C_{12102121021210121212}\\
& &+5C_{102121021210121212}+([2]^4+[2]^{2}+2)C_{x_{\beta}w_{0}}\}.\\
C_{121212}C_{02121021210121212}&=&[2]\{C_{x_{\alpha}x_{\beta}w_{0}}+C_{x_{\alpha}w_{0}}+C_{x_{\alpha}^{2}w_{0}}+4C_{12102121021210121212}\\
& &+5C_{102121021210121212}+([2]^4+[2]^{2}+2)C_{x_{\beta}w_{0}}\}\\
& &-4[2][C_{12102121021210121212}+2C_{102121021210121212}\\
& &+([2]^2+1)C_{x_{\beta}w_{0}}]+3[2](C_{102121021210121212}+C_{x_{\beta}w_{0}})\\
&=&[2]C_{x_{\alpha}x_{\beta}w_{0}}+([2]^{5}-3[2]^{3}+[2])C_{x_{\beta}w_{0}}+[2]C_{x_{\alpha}w_{0}}\\
& &+[2]C_{x_{\alpha}^{2}w_{0}}.\\
\end{eqnarray*}

(12**) $C_{121212}C_{210121210121212}$
\begin{eqnarray*}
C_{2}C_{210121210121212}&=&[2]C_{210121210121212},\\
C_{1}C_{210121210121212}&=&C_{1210121210121212}+C_{10121210121212},\\
C_{2}C_{1210121210121212}&=&C_{21210121210121212}+C_{210121210121212},\\
C_{2}C_{10121210121212}&=&C_{210121210121212}+C_{0121210121212},\\
C_{1}C_{21210121210121212}&=&C_{x_{\alpha}^{2}w_{0}}+C_{x_{\alpha}w_{0}}+C_{x_{\beta}w_{0}}+C_{w_{0}}+C_{1210121210121212}.\\
\end{eqnarray*}
\begin{eqnarray*}
C_{1}C_{2}C_{210121210121212}&=&[2][C_{1210121210121212}+C_{10121210121212}],\\
C_{2}C_{1}C_{2}C_{210121210121212}&=&[2]C_{21210121210121212}+2[2]C_{210121210121212}\\
& &+[2]C_{0121210121212},\\
C_{1}C_{2}C_{1}C_{2}C_{210121210121212}&=&[2][C_{x_{\alpha}^{2}w_{0}}+2C_{x_{\alpha}w_{0}}+C_{x_{\beta}w_{0}}+C_{w_{0}}\\
& &+3C_{1210121210121212}+3C_{10121210121212}],\\
C_{2}C_{1}C_{2}C_{1}C_{2}C_{210121210121212}&=&[2][[2]C_{x_{\alpha}^{2}w_{0}}+2[2]C_{x_{\alpha}w_{0}}+[2]C_{x_{\beta}w_{0}}\\
& &+[2]C_{w_{0}}+3C_{21210121210121212}\\
& &+6C_{210121210121212}+3C_{0121210121212}],\\
C_{1}C_{2}C_{1}C_{2}C_{1}C_{2}C_{210121210121212}&=&[2][([2]^{2}+3)C_{x_{\alpha}^{2}w_{0}}+(2[2]^{2}+6)C_{x_{\alpha}w_{0}}\\
& &+([2]^{2}+3)C_{x_{\beta}w_{0}}+([2]^{2}+3)C_{w_{0}}\\
& &+9C_{1210121210121212}+9C_{10121210121212}].\\
C_{121212}C_{210121210121212}&=&[2][([2]^{2}+3)C_{x_{\alpha}^{2}w_{0}}+(2[2]^{2}+6)C_{x_{\alpha}w_{0}}\\
& &+([2]^{2}+3)C_{x_{\beta}w_{0}}+([2]^{2}+3)C_{w_{0}}\\
& &+9C_{1210121210121212}+9C_{10121210121212}]\\
& &-4[2][C_{x_{\alpha}^{2}w_{0}}+2C_{x_{\alpha}w_{0}}+C_{x_{\beta}w_{0}}+C_{w_{0}}\\
& &+3C_{1210121210121212}+3C_{10121210121212}]\\
& &+3[2][C_{1210121210121212}+C_{10121210121212}]\\
&=&([2]^{3}-[2])C_{x_{\alpha}^{2}w_{0}}+2([2]^{3}-[2])C_{x_{\alpha}w_{0}}\\
& &+([2]^{3}-[2])C_{x_{\beta}w_{0}}+([2]^{3}-[2])C_{w_{0}}.\\
\end{eqnarray*}

(13) Computing $C_{1212120}C_{121212}$
\begin{eqnarray*}
C_{1212120}C_{121212}&=&C_{121212}C_{0121212}=[2]C_{x_{\alpha}w_{0}}+([2]^{5}-3[2]^{3}+[2])C_{w_{0}}.
\end{eqnarray*}

(14) Computing $C_{1212120}C_{0121212}$
\begin{eqnarray*}
C_{121212}C_{0}&=&C_{1212120},\\
C_{0}C_{121212}&=&C_{0121212},\\
C_{1212120}C_{0121212}&=&C_{121212}C_{0}C_{0121212}\\
&=&[2]C_{121212}C_{0121212}\\
&=&[2]^{2}C_{x_{\alpha}w_{0}}+([2]^{6}-3[2]^{4}+[2]^{2})C_{w_{0}}.\\
\end{eqnarray*}

(15) Computing $C_{1212120}C_{10121212}$
\begin{eqnarray*}
C_{0}C_{10121212}&=&[2]C_{10121212},\\
C_{1212120}C_{10121212}&=&C_{121212}C_{0}C_{10121212}\\
&=&[2]C_{121212}C_{10121212}\\
&=&[2]^{3}C_{x_{\alpha}w_{0}}+([2]^{5}-2[2]^{3})C_{w_{0}}.\\
\end{eqnarray*}

(16) Computing $C_{1212120}C_{210121212}$
\begin{eqnarray*}
C_{0}C_{210121212}&=&[2]C_{210121212},\\
C_{1212120}C_{210121212}&=&C_{121212}C_{0}C_{210121212}\\
&=&[2]C_{121212}C_{210121212}\\
&=&([2]^{4}-[2]^{2})C_{x_{\alpha}w_{0}}+([2]^{4}-[2]^{2})C_{w_{0}}.\\
\end{eqnarray*}

(17) Computing $C_{1212120}C_{1210121212}$
\begin{eqnarray*}
C_{0}C_{1210121212}&=&C_{01210121212}+C_{210121212}.\\
C_{1212120}C_{1210121212}&=&C_{121212}C_{0}C_{1210121212}\\
&=&C_{121212}(C_{01210121212}+C_{210121212})\\
&=&[2]C_{x_{\beta}121212}+([2]^{3}-[2])C_{x_{\alpha}121212}+[2]C_{w_{0}}\\
& &+([2]^{3}-[2])C_{x_{\alpha}121212}+([2]^{3}-[2])C_{121212}\\
&=&[2]C_{x_{\beta}w_{0}}+2([2]^{3}-[2])C_{x_{\alpha}w_{0}}+[2]^{3}C_{w_{0}}.\\
\end{eqnarray*}

(18) Computing $C_{1212120}C_{21210121212}$
\begin{eqnarray*}
C_{0}C_{21210121212}&=&C_{201210121212}.\\
C_{1212120}C_{21210121212}&=&C_{121212}C_{0}C_{21210121212}\\
&=&C_{121212}C_{201210121212}\\
&=&[2]^{2}C_{x_{\beta}121212}+([2]^{4}-[2]^{2})C_{x_{\alpha}121212}+[2]^{2}C_{121212}\\
&=&[2]^{2}C_{x_{\beta}w_{0}}+([2]^{4}-[2]^{2})C_{x_{\alpha}w_{0}}+[2]^{2}C_{w_{0}}\\
\end{eqnarray*}

(19) Computing $C_{1212120}C_{01210121212}$
\begin{eqnarray*}
C_{0}C_{01210121212}&=&[2]C_{01210121212}.\\
C_{1212120}C_{01210121212}&=&C_{121212}C_{0}C_{01210121212}\\
&=&[2]C_{121212}C_{01210121212}\\
&=&[2][[2]C_{x_{\beta}121212}+([2]^{3}-[2])C_{x_{\alpha}121212}+[2]C_{121212}]\\
&=&[2]^{2}C_{x_{\beta}w_{0}}+([2]^{4}-[2]^{2})C_{x_{\alpha}w_{0}}+[2]^{2}C_{w_{0}}\\
\end{eqnarray*}

(20) Computing $C_{1212120}C_{201210121212}$
\begin{eqnarray*}
C_{0}C_{201210121212}&=&[2]C_{201210121212}.\\
C_{1212120}C_{201210121212}&=&C_{121212}C_{0}C_{201210121212}\\
&=&[2]C_{121212}C_{201210121212}\\
&=&[2][[2]^{2}C_{x_{\beta}121212}+([2]^{4}-[2]^{2})C_{x_{\alpha}121212}+[2]^{2}C_{121212}]\\
&=&[2]^{3}C_{x_{\beta}w_{0}}+([2]^{5}-[2]^{3})C_{x_{\alpha}w_{0}}+[2]^{3}C_{w_{0}}\\
\end{eqnarray*}

(21) Computing $C_{1212120}C_{1201210121212}$
\begin{eqnarray*}
C_{0}C_{1201210121212}&=&C_{01201210121212}+C_{201210121212}+C_{10121212}\\
&=&C_{10121210121212}+C_{201210121212}+C_{10121212}.\\
\end{eqnarray*}
Now,we need to Computing $C_{121212}C_{10121210121212}$;
\begin{eqnarray*}
C_{210121210121212}&=&C_{2}C_{10121210121212}-C_{0121210121212}\\
&=&C_{2}(C_{1}C_{0121210121212}-C_{x_{\alpha}121212})-C_{0121210121212}\\
&=&C_{2}C_{1}C_{0121210121212}-[2]C_{x_{\alpha}121212})-C_{0121210121212},\\
C_{121212}C_{210121210121212}&=&C_{121212}C_{2}C_{1}C_{0121210121212}-[2]C_{121212}C_{x_{\alpha}121212}\\
& &-C_{121212}C_{0121210121212}\\
&=&([2]^{2}-1)C_{121212}C_{0121210121212}-[2]C_{121212}C_{x_{\alpha}121212},\\
C_{121212}C_{210121210121212}&=&([2]^{2}-1)C_{121212}C_{0121210121212}-[2]C_{121212}C_{x_{\alpha}121212}\\
([2]^{2}-1)C_{121212}C_{0121210121212}&=&C_{121212}C_{210121210121212}+[2]C_{121212}C_{x_{\alpha}121212}\\
&=&([2]^{3}-[2])C_{x_{\alpha}^{2}w_{0}}+2([2]^{3}-[2])C_{x_{\alpha}w_{0}}\\
& &+([2]^{3}-[2])C_{x_{\beta}w_{0}}+([2]^{3}-[2])C_{w_{0}}\\
& &+([2]^{7}-4[2]^{5}+3[2]^{3})C_{x_{\alpha}w_{0}}\\
&=&([2]^{3}-[2])C_{x_{\alpha}^{2}w_{0}}+([2]^{7}-4[2]^{5}+5[2]^{3}-2[2])C_{x_{\alpha}w_{0}}\\
& &+([2]^{3}-[2])C_{x_{\beta}w_{0}}+([2]^{3}-[2])C_{w_{0}},\\
C_{121212}C_{0121210121212}&=&[2]C_{x_{\alpha}^{2}w_{0}}+([2]^{5}-3[2]^{3}+2[2])C_{x_{\alpha}w_{0}}\\
& &+[2]C_{x_{\beta}w_{0}}+[2]C_{w_{0}},\\
C_{210121210121212}&=&C_{2}C_{10121210121212}-C_{0121210121212},\\
C_{121212}C_{210121210121212}&=&[2]C_{121212}C_{10121210121212}-C_{121212}C_{0121210121212},\\
C_{121212}C_{10121210121212}[2]&=&C_{121212}C_{210121210121212}-C_{121212}C_{0121210121212},\\
C_{121212}C_{10121210121212}&=&[2]^{2}C_{x_{\alpha}^{2}w_{0}}+([2]^{4}-[2]^{2})C_{x_{\alpha}w_{0}}\\
& &+[2]^{2}C_{x_{\beta}w_{0}}+[2]^{2}C_{w_{0}}\\
\end{eqnarray*}

\begin{eqnarray*}
C_{1212120}C_{1201210121212}&=&C_{121212}C_{0}C_{1201210121212}\\
&=&C_{121212}(C_{10121210121212}+C_{201210121212}+C_{10121212})\\
&=&[2]^{2}C_{x_{\alpha}^{2}w_{0}}+([2]^{4}-[2]^{2})C_{x_{\alpha}w_{0}}\\
& &+[2]^{2}C_{x_{\beta}w_{0}}+[2]^{2}C_{w_{0}}\\
& &+[2]^{2}C_{x_{\beta}121212}+([2]^{4}-[2]^{2})C_{x_{\alpha}121212}+[2]^{2}C_{w_{0}}\\
& &+[2]^{2}C_{x_{\alpha}121212}+([2]^{4}-2[2]^{2})C_{w_{0}}\\
&=&[2]^{2}C_{x_{\alpha}^{2}121212}+(2[2]^{4}-[2]^{2})C_{x_{\alpha}121212}+2[2]^{2}C_{x_{\beta}121212}+[2]^{4}C_{w_{0}}\\
\end{eqnarray*}

(22) Computing $C_{1212120}C_{21201210121212}$
\begin{eqnarray*}
C_{0}C_{21201210121212}&=&C_{021201210121212}+C_{0121212}+C_{210121212}\\
&=&C_{210121210121212}+C_{0121212}+C_{210121212},\\
C_{1212120}C_{21201210121212}&=&C_{121212}C_{0}C_{21201210121212}\\
&=&C_{121212}(C_{210121210121212}+C_{0121212}+C_{210121212})\\
&=&([2]^{3}-[2])C_{x_{\alpha}^{2}w_{0}}+2([2]^{3}-[2])C_{x_{\alpha}w_{0}}\\
& &+([2]^{3}-[2])C_{x_{\beta}w_{0}}+([2]^{3}-[2])C_{w_{0}}\\
& &+[2]C_{x_{\beta}w_{0}}+([2]^{5}-3[2]^{3}+[2])C_{w_{0}}\\
& &+([2]^{3}-[2])C_{x_{\alpha}w_{0}}+([2]^{3}-[2])C_{w_{0}}\\
&=&([2]^{3}-[2])C_{x_{\alpha}^{2}w_{0}}+(3[2]^{3}-2[2])C_{x_{\alpha}w_{0}}\\
& &+([2]^{3}-[2])C_{x_{\beta}w_{0}}+([2]^{5}-[2]^{3}-[2])C_{w_{0}}\\
\end{eqnarray*}

(23) Computing $C_{1212120}C_{121201210121212}$
\begin{eqnarray*}
C_{0}C_{121201210121212}&=&C_{0121201210121212}+C_{10121212},\\
C_{1212120}C_{121201210121212}&=&C_{121212}C_{0}C_{201210121212}\\
&=&C_{121212}(C_{0121201210121212}+C_{10121212})\\
&=&C_{x_{\alpha}x_{\beta}w_{0}}+[2]^{2}C_{x_{\alpha}^{2}w_{0}}+(2[2]^{2})C_{x_{\alpha}w_{0}}\\
& &+([2]^{4}-2[2]^{2})C_{x_{\beta}w_{0}}+([2]^{4}-2[2]^{2})C_{w_{0}}\\
& &+[2]^{2}C_{x_{\alpha}^{2}w_{0}}+([2]^{4}-2[2]^{2})C_{w_{0}}\\
&=&C_{x_{\alpha}x_{\beta}w_{0}}+[2]^{2}C_{x_{\alpha}^{2}w_{0}}+(3[2]^{2})C_{x_{\alpha}w_{0}}\\
& &+([2]^{4}-2[2]^{2})C_{x_{\beta}w_{0}}+(2[2]^{4}-4[2]^{2})C_{w_{0}},\\
\end{eqnarray*}

(24) Computing $C_{1212120}C_{0121201210121212}$
\begin{eqnarray*}
C_{0}C_{0121201210121212}&=&[2]C_{0121201210121212}.\\
C_{1212120}C_{0121201210121212}&=&C_{121212}C_{0}C_{0121201210121212},\\
&=&[2]C_{121212}C_{0121201210121212},\\
\end{eqnarray*}

(25) Computing $C_{12121201}C_{121212}$
\begin{eqnarray*}
C_{12121201}C_{121212}&=&C_{121212}C_{10121212}=[2]^{2}C_{x_{\alpha}w_{0}}+([2]^{4}-2[2]^{2})C_{w_{0}}.
\end{eqnarray*}

(26) Computing $C_{12121201}C_{0121212}$
\begin{eqnarray*}
C_{12121201}C_{0121212}&=&C_{1212120}C_{10121212}=[2]^{3}C_{x_{\alpha}w_{0}}+([2]^{5}-2[2]^{3})C_{w_{0}}.
\end{eqnarray*}

(27) Computing $C_{12121201}C_{10121212}$
\begin{eqnarray*}
C_{1212120}C_{1}&=&C_{12121201}+C_{121212},\\
C_{12121201}&=&C_{1212120}C_{1}-C_{121212}.\\
C_{12121201}C_{10121212}&=&(C_{1212120}C_{1}-C_{121212})C_{10121212},\\
&=&C_{1212120}C_{1}C_{10121212}-C_{121212}C_{10121212},\\
&=&[2]C_{1212120}C_{10121212}-C_{121212}C_{10121212},\\
&=&([2]^{4}-[2]^{2})C_{x_{\alpha}w_{0}}+([2]^{6}-3[2]^{4}+2[2]^{2})C_{w_{0}},\\
\end{eqnarray*}

(28) Computing $C_{12121201}C_{210121212}$
\begin{eqnarray*}
C_{12121201}&=&C_{1212120}C_{1}-C_{121212},\\
C_{1}C_{210121212}&=&C_{1210121212}+C_{10121212}.\\
C_{12121201}C_{210121212}&=&(C_{1212120}C_{1}-C_{121212})C_{210121212},\\
&=&C_{1212120}C_{1}C_{210121212}-C_{121212}C_{210121212},\\
&=&C_{1212120}C_{1210121212}+C_{1212120}C_{10121212}-C_{121212}C_{210121212},\\
&=&(2[2]^{3}-[2])C_{x_{\alpha}w_{0}}+[2]C_{x_{\beta}w_{0}}+([2]^{5}-2[2]^{3}+[2])C_{w_{0}},\\
\end{eqnarray*}

(29) Computing $C_{12121201}C_{1210121212}$
\begin{eqnarray*}
C_{12121201}&=&C_{1212120}C_{1}-C_{121212}.\\
C_{12121201}C_{1210121212}&=&(C_{1212120}C_{1}-C_{121212})C_{1210121212},\\
&=&C_{1212120}C_{1}C_{1210121212}-C_{121212}C_{1210121212},\\
&=&[2]C_{1212120}C_{1210121212}-C_{121212}C_{1210121212},\\
&=&[2]^{4}C_{x_{\alpha}w_{0}}+[2]^{2}C_{x_{\beta}w_{0}}+([2]^{4}-[2]^{2})C_{w_{0}},\\
\end{eqnarray*}

(30) Computing $C_{12121201}C_{21210121212}$
\begin{eqnarray*}
C_{1}C_{21210121212}&=&C_{121210121212}+C_{1210121212},\\
C_{0}C_{121210121212}&=&C_{0121210121212},\\
C_{121212}C_{0121210121212}&=&[2]C_{x_{\alpha}^{2}w_{0}}+([2]^{5}-3[2]^{3}+2[2])C_{x_{\alpha}w_{0}}+[2]C_{x_{\beta}w_{0}}+[2]C_{w_{0}},\\
C_{12121201}C_{21210121212}&=&(C_{1212120}C_{1}-C_{121212})C_{21210121212},\\
&=&C_{1212120}C_{1}C_{21210121212}-C_{121212}C_{21210121212},\\
&=&C_{1212120}C_{121210121212}+C_{1212120}C_{1210121212}-C_{121212}C_{21210121212},\\
&=&C_{121212}C_{0121210121212}+C_{1212120}C_{1210121212}-C_{121212}C_{21210121212},\\
&=&[2]C_{x_{\alpha}^{2}w_{0}}+(2[2]^{3}-[2])C_{x_{\alpha}w_{0}}+2[2]C_{x_{\beta}w_{0}}+[2]^{3}C_{w_{0}},\\
\end{eqnarray*}

(31) Computing $C_{12121201}C_{01210121212}$
\begin{eqnarray*}
C_{12121201}C_{01210121212}&=&(C_{1212120}C_{1}-C_{121212})C_{01210121212},\\
&=&C_{1212120}C_{1}C_{01210121212}-C_{121212}C_{01210121212},\\
&=&[2]C_{1212120}C_{01210121212}-C_{121212}C_{01210121212},\\
&=&([2]^{5}-2[2]^{3}+[2])C_{x_{\alpha}w_{0}}+([2]^{3}-[2])C_{x_{\beta}w_{0}}+([2]^{3}-[2])C_{w_{0}},\\
\end{eqnarray*}

(32) Computing $C_{12121201}C_{201210121212}$
\begin{eqnarray*}
C_{1}C_{201210121212}&=&C_{1201210121212}+C_{01210121212}.\\
C_{12121201}C_{201210121212}&=&(C_{1212120}C_{1}-C_{121212})C_{201210121212},\\
&=&C_{1212120}C_{1}C_{201210121212}-C_{121212}C_{201210121212},\\
&=&C_{1212120}C_{1201210121212}+C_{1212120}C_{01210121212}-C_{121212}C_{201210121212},\\
&=&[2]^{2}C_{x_{\alpha}^{2}w_{0}}+(2[2]^{4}-[2]^{2})C_{x_{\alpha}w_{0}}+2[2]^{2}C_{x_{\beta}w_{0}}+[2]^{4}C_{w_{0}},\\
\end{eqnarray*}

(33) Computing $C_{12121201}C_{1201210121212}$
\begin{eqnarray*}
C_{12121201}C_{1201210121212}&=&(C_{1212120}C_{1}-C_{121212})C_{1201210121212},\\
&=&C_{1212120}C_{1}C_{1201210121212}-C_{121212}C_{1201210121212},\\
&=&[2]C_{1212120}C_{1201210121212}-C_{121212}C_{1201210121212},\\
&=&[2]^{3}C_{x_{\alpha}^{2}w_{0}}+([2]^{5}+[2]^{3}-[2])C_{x_{\alpha}w_{0}}+([2]^{3}+[2])C_{x_{\beta}w_{0}},\\
& &+([2]^{5}-[2]^{3}+[2])C_{w_{0}},\\
\end{eqnarray*}

(34) Computing $C_{12121201}C_{21201210121212}$
\begin{eqnarray*}
C_{1}C_{21201210121212}&=&C_{121201210121212}+C_{1201210121212}.\\
C_{12121201}C_{21201210121212}&=&(C_{1212120}C_{1}-C_{121212})C_{21201210121212},\\
&=&C_{1212120}C_{1}C_{21201210121212}-C_{121212}C_{21201210121212},\\
&=&C_{1212120}C_{121201210121212}+C_{1212120}C_{1201210121212},\\
& &-C_{121212}C_{21201210121212},\\
&=&C_{x_{\alpha}x_{\beta}w_{0}}+2[2]^{2}C_{x_{\alpha}^{2}w_{0}}+2([2]^{4}-[2]^{2})C_{w_{0}},\\
& &+([2]^{4}+3[2]^{2})C_{x_{\alpha}w_{0}}+2[2]^{2}C_{x_{\beta}w_{0}},\\
\end{eqnarray*}

(35) Computing $C_{12121201}C_{121201210121212}$
\begin{eqnarray*}
C_{12121201}C_{121201210121212}&=&(C_{1212120}C_{1}-C_{121212})C_{121201210121212},\\
&=&C_{1212120}C_{1}C_{121201210121212}-C_{121212}C_{121201210121212},\\
&=&[2]C_{1212120}C_{121201210121212}-C_{121212}C_{121201210121212},\\
&=&[2]C_{x_{\alpha}x_{\beta}121212}+[2]^{3}C_{x_{\alpha}^{2}121212}+([2]^{5}-[2]^{3}-[2])C_{w_{0}},\\
& &+(2[2]^{3}+[2])C_{x_{\alpha}w_{0}}+([2]^{3}-[2])C_{x_{\beta}w_{0}},\\
\end{eqnarray*}

(36) Computing $C_{12121201}C_{0121201210121212}$
\begin{eqnarray*}
C_{1}C_{0121201210121212}&=&C_{10121201210121212}+C_{121201210121212},\\
&=&C_{01210121210121212}+C_{121201210121212},\\
C_{1212120}C_{121201210121212}&=&C_{121212}C_{0121201210121212}+C_{121212}C_{10121212}.\\
C_{12121201}C_{0121201210121212}&=&(C_{1212120}C_{1}-C_{121212})C_{0121201210121212},\\
&=&C_{1212120}C_{1}C_{0121201210121212}-C_{121212}C_{0121201210121212},\\
&=&C_{1212120}C_{01210121210121212}+C_{1212120}C_{121201210121212},\\
& &+C_{121212}C_{10121212}-C_{121212}C_{0121201210121212},\\
&=&[2]C_{121212}C_{01210121210121212}+C_{121212}C_{10121212},\\
&=&[2]C_{121212}C_{10121201210121212}+C_{121212}C_{10121212},\\
&=&[2]^{2}C_{x_{\alpha}x_{\beta}w_{0}}+([2]^{4}-[2]^{2})C_{w_{0}}+([2]^{4}+[2]^{2})C_{x_{\alpha}w_{0}},\\
& &+([2]^{4}-[2]^{2})C_{x_{\beta}w_{0}}+[2]^{4}C_{x_{\alpha}^{2}w_{0}},\\
& &+[2]^{2}C_{x_{\alpha}w_{0}}+([2]^{4}-[2]^{2})C_{w_{0}},\\
&=&[2]^{2}C_{x_{\alpha}x_{\beta}w_{0}}+(2[2]^{4}-3[2]^{2})C_{w_{0}}+([2]^{4}+2[2]^{2})C_{x_{\alpha}w_{0}},\\
& &+([2]^{4}-[2]^{2})C_{x_{\beta}w_{0}}+[2]^{4}C_{x_{\alpha}^{2}w_{0}},\\
\end{eqnarray*}
Computing $C_{121212}C_{10121201210121212}=C_{121212}C_{01210121210121212}$
\begin{eqnarray*}
C_{2}C_{10121201210121212}&=&C_{210121201210121212},\\
C_{1}C_{210121201210121212}&=&C_{1210121201210121212}+C_{10121201210121212},\\
C_{2}C_{1210121201210121212}&=&C_{21210121201210121212}+C_{210121201210121212}+C_{w_{0}},\\
& &+C_{x_{\alpha}w_{0}}+C_{x_{\beta}w_{0}}+C_{x_{\alpha}^{2}w_{0}},\\
C_{1}C_{21210121201210121212}&=&C_{121210121201210121212}+C_{1210121201210121212},\\
C_{2}C_{121210121201210121212}&=&C_{x_{\alpha}x_{\beta}w_{0}}+C_{21210121201210121212}+2C_{x_{\alpha}w_{0}}+C_{x_{\alpha}^{2}w_{0}},\\
\end{eqnarray*}
\begin{eqnarray*}
C_{1}C_{2}C_{10121201210121212}&=&C_{1210121201210121212}+C_{10121201210121212},\\
C_{2}C_{1}C_{2}C_{10121201210121212}&=&C_{21210121201210121212}+2C_{210121201210121212}+C_{w_{0}},\\
& &+C_{x_{\alpha}w_{0}}+C_{x_{\beta}w_{0}}+C_{x_{\alpha}^{2}w_{0}},\\
C_{1}C_{2}C_{1}C_{2}C_{10121201210121212}&=&C_{121210121201210121212}+3C_{1210121201210121212},\\
& &+2C_{10121201210121212}+[2]C_{w_{0}}+[2]C_{x_{\alpha}w_{0}}\\
& &+[2]C_{x_{\beta}w_{0}}+[2]C_{x_{\alpha}^{2}w_{0}},\\
C_{2}C_{1}C_{2}C_{1}C_{2}C_{0121201210121212}&=&C_{x_{\alpha}x_{\beta}w_{0}}+4C_{21210121201210121212}+5C_{210121201210121212}\\
& &+([2]^{2}+5)C_{x_{\alpha}w_{0}}+([2]^{2}+4)C_{x_{\alpha}^{2}w_{0}}\\
& &+([2]^{2}+3)C_{w_{0}}+([2]^{2}+3)C_{x_{\beta}w_{0}},\\
C_{1}C_{2}C_{1}C_{2}C_{1}C_{2}C_{0121201210121212}&=&[2]C_{x_{\alpha}x_{\beta}w_{0}}+4C_{121210121201210121212}+9C_{1210121201210121212}\\
& &+5C_{10121201210121212}+([2]^{3}+3[2])C_{w_{0}}+([2]^{3}+5[2])C_{x_{\alpha}w_{0}}\\
& &+([2]^{3}+3[2])C_{x_{\beta}w_{0}}+([2]^{3}+4[2])C_{x_{\alpha}^{2}w_{0}},\\
C_{121212}C_{10121201210121212}&=&(C_{2}C_{1}C_{1}C_{2}C_{1}C_{2}-4C_{2}C_{1}C_{2}C_{1}+3C_{2}C_{1})C_{10121210121212}\\
&=&[2]C_{x_{\alpha}x_{\beta}w_{0}}+4C_{121210121201210121212}+9C_{1210121201210121212}\\
& &+5C_{10121201210121212}+([2]^{3}+3[2])C_{w_{0}}+([2]^{3}+5[2])C_{x_{\alpha}w_{0}}\\
& &+([2]^{3}+3[2])C_{x_{\beta}w_{0}}+([2]^{3}+4[2])C_{x_{\alpha}^{2}w_{0}}\\
& &-4\{C_{121210121201210121212}+3C_{1210121201210121212}+[2]C_{x_{\beta}w_{0}}\\
& &+2C_{10121201210121212}+[2]C_{w_{0}}+[2]C_{x_{\alpha}w_{0}}+[2]C_{x_{\alpha}^{2}w_{0}}\}\\
& &+3\{C_{1210121201210121212}+C_{10121201210121212}\}\\
&=&[2]C_{x_{\alpha}x_{\beta}w_{0}}+([2]^{3}-[2])C_{w_{0}}+([2]^{3}+[2])C_{x_{\alpha}w_{0}}\\
& &+([2]^{3}-[2])C_{x_{\beta}w_{0}}+[2]^{3}C_{x_{\alpha}^{2}w_{0}}\\
&=&C_{1212120}C_{1210121210121212}.\\
\end{eqnarray*}

(37) Computing $C_{121212012}C_{121212}$
\begin{eqnarray*}
C_{121212012}C_{121212}&=&C_{121212}C_{210121212}=([2]^{3}-[2])C_{x_{\alpha}w_{0}}+([2]^{3}-[2])C_{w_{0}}.
\end{eqnarray*}

(38) Computing $C_{121212012}C_{0121212}$
\begin{eqnarray*}
C_{121212012}C_{0121212}&=&C_{1212120}C_{210121212}=([2]^{4}-[2]^{2})C_{x_{\alpha}w_{0}}+([2]^{4}-[2]^{2})C_{w_{0}}.
\end{eqnarray*}

(39) Computing $C_{121212012}C_{10121212}$
\begin{eqnarray*}
C_{121212012}C_{10121212}&=&C_{12121201}C_{210121212}==(2[2]^{3}-[2])C_{x_{\alpha}w_{0}}+[2]C_{x_{\beta}w_{0}}\\
& &+([2]^{5}-2[2]^{3}+[2])C_{w_{0}}.
\end{eqnarray*}

(40) Computing $C_{121212012}C_{210121212}$
\begin{eqnarray*}
C_{12121201}C_{2}&=&C_{121212012}+C_{1212120},\\
C_{121212012}&=&C_{12121201}C_{2}-C_{1212120}.\\
C_{121212012}C_{210121212}&=&(C_{12121201}C_{2}-C_{1212120})C_{210121212}\\
&=&C_{12121201}C_{2}C_{210121212}-C_{1212120}C_{210121212}\\
&=&[2]C_{12121201}C_{210121212}-C_{1212120}C_{210121212}\\
&=&[2]^{2}C_{x_{\beta}w_{0}}+[2]^{4}C_{x_{\alpha}w_{0}}+([2]^{6}-3[2]^{4}+2[2]^{2})C_{w_{0}}\\
\end{eqnarray*}

(41) Computing $C_{121212012}C_{1210121212}$
\begin{eqnarray*}
C_{2}C_{1210121212}&=&C_{21210121212}+C_{210121212},\\
C_{121212012}C_{1210121212}&=&(C_{12121201}C_{2}-C_{1212120})C_{1210121212}\\
&=&C_{12121201}C_{2}C_{1210121212}-C_{1212120}C_{1210121212}\\
&=&C_{12121201}C_{21210121212}+C_{212121201}C_{0121212}-C_{1212120}C_{1210121212}\\
&=&[2]C_{x_{\alpha}^{2}w_{0}}+2[2]C_{x_{\beta}w_{0}}+2[2]^{3}C_{x_{\alpha}w_{0}}+([2]^{5}-2[2]^{3}+[2])C_{w_{0}}\\
\end{eqnarray*}

(42) Computing $C_{121212012}C_{21210121212}$
\begin{eqnarray*}
C_{121212012}C_{21210121212}&=&(C_{12121201}C_{2}-C_{1212120})C_{21210121212}\\
&=&C_{12121201}C_{2}C_{21210121212}-C_{1212120}C_{21210121212}\\
&=&[2]C_{12121201}C_{21210121212}-C_{1212120}C_{21210121212}\\
&=&[2]^{2}C_{x_{\alpha}^{2}w_{0}}+[2]^{4}C_{x_{\alpha}w_{0}}+[2]^{2}C_{x_{\beta}w_{0}}+([2]^{4}-[2]^{2})C_{w_{0}}\\
\end{eqnarray*}

(43) Computing $C_{121212012}C_{01210121212}$
\begin{eqnarray*}
C_{2}C_{01210121212}&=&C_{201210121212}.\\
C_{121212012}C_{1210121212}&=&(C_{12121201}C_{2}-C_{1212120})C_{01210121212}\\
&=&C_{12121201}C_{2}C_{01210121212}-C_{1212120}C_{01210121212}\\
&=&C_{12121201}C_{201210121212}-C_{1212120}C_{1210121212}\\
&=&[2]^{2}C_{x_{\alpha}^{2}w_{0}}+[2]^{4}C_{x_{\alpha}w_{0}}+[2]^{2}C_{x_{\beta}w_{0}}+([2]^{4}-[2]^{2})C_{w_{0}}\\
\end{eqnarray*}

(44) Computing $C_{121212012}C_{201210121212}$
\begin{eqnarray*}
C_{121212012}C_{201210121212}&=&(C_{12121201}C_{2}-C_{1212120})C_{201210121212}\\
&=&C_{12121201}C_{2}C_{201210121212}-C_{1212120}C_{201210121212}\\
&=&[2]C_{12121201}C_{201210121212}-C_{1212120}C_{201210121212}\\
&=&[2]^{3}C_{x_{\alpha}^{2}w_{0}}+[2]^{5}C_{x_{\alpha}w_{0}}+[2]^{3}C_{x_{\beta}w_{0}}+([2]^{5}-[2]^{3})C_{w_{0}}\\
\end{eqnarray*}

(45) Computing $C_{121212012}C_{1201210121212}$
\begin{eqnarray*}
C_{12121201}C_{x_{\alpha}w_{0}}&=&(C_{1212120}C_{1}-C_{121212})C_{x_{\alpha}w_{0}}\\
&=&[2]C_{1212120}C_{x_{\alpha}w_{0}}-C_{121212}C_{x_{\alpha}w_{0}}\\
&=&[2]^{2}C_{x_{\alpha}^{2}w_{0}}+([2]^{6}-3[2]^{4}+2[2]^{2})C_{x_{\alpha}w_{0}}+[2]^{2}C_{x_{\beta}w_{0}}\\
& &+[2]^{2}C_{w_{0}}-([2]^{6}-4[2]^{4}+3[2]^{2})C_{x_{\alpha}w_{0}}\\
&=&[2]^{2}C_{x_{\alpha}^{2}w_{0}}+([2]^{4}-[2]^{2})C_{x_{\alpha}w_{0}}+[2]^{2}C_{x_{\beta}w_{0}}+[2]^{2}C_{w_{0}}.\\
\end{eqnarray*}
\begin{eqnarray*}
C_{2}C_{1201210121212}&=&C_{21201210121212}+C_{x_{\alpha}w_{0}}+C_{201210121212}\\
C_{121212012}C_{1201210121212}&=&(C_{12121201}C_{2}-C_{1212120})C_{1201210121212}\\
&=&C_{12121201}C_{2}C_{1201210121212}-C_{1212120}C_{1201210121212}\\
&=&C_{12121201}C_{21201210121212}+C_{12121201}C_{x_{\alpha}w_{0}}\\
& &+C_{12121201}C_{201210121212}-C_{1212120}C_{1201210121212}\\
&=&C_{x_{\alpha}x_{\beta}w_{0}}+3[2]^{2}C_{x_{\alpha}^{2}w_{0}}+2([2]^{4}+[2]^{2})C_{x_{\alpha}w_{0}}\\
& &+3[2]^{2}C_{x_{\beta}w_{0}}+(2[2]^{4}-[2]^{2})C_{w_{0}}.\\
\end{eqnarray*}

(46) Computing $C_{121212012}C_{21201210121212}$
\begin{eqnarray*}
C_{121212012}C_{201210121212}&=&(C_{12121201}C_{2}-C_{1212120})C_{21201210121212}\\
&=&C_{12121201}C_{2}C_{21201210121212}-C_{1212120}C_{21201210121212}\\
&=&[2]C_{x_{\alpha}x_{\beta}w_{0}}+([2]^{3}+[2])C_{x_{\alpha}^{2}w_{0}}+([2]^{5}+2[2])C_{x_{\alpha}w_{0}}\\
& &+([2]^{3}+[2])C_{x_{\beta}w_{0}}+([2]^{5}-[2]^{3}+[2])C_{w_{0}}.\\
\end{eqnarray*}

(47) Computing $C_{121212012}C_{121201210121212}$
\begin{eqnarray*}
C_{2}C_{121201210121212}&=&C_{x_{\beta}w_{0}}+C_{21201210121212}+C_{w_{0}},\\
C_{121212012}C_{121201210121212}&=&(C_{12121201}C_{2}-C_{1212120})C_{121201210121212}\\
&=&C_{12121201}C_{2}C_{121201210121212}-C_{1212120}C_{121201210121212}\\
&=&C_{12121201}C_{x_{\beta}w_{0}}+C_{12121201}C_{21201210121212}\\
& &+C_{12121201}C_{w_{0}}-C_{1212120}C_{121201210121212}\\
&=&[2]^{2}C_{x_{\alpha}x_{\beta}w_{0}}+2[2]^{2}C_{x_{\beta}w_{0}}+([2]^{4}+2[2]^{2})C_{x_{\alpha}w_{0}}\\
& &+2[2]^{2}C_{x_{\alpha}^{2}w_{0}}+[2]^{4}C_{w_{0}}.\\
\end{eqnarray*}

\begin{eqnarray*}
C_{12121201}C_{x_{\beta}w_{0}}&=&(C_{121212}C_{0}C_{1}-C_{121212})C_{x_{\beta}w_{0}}\\
&=&[2]C_{121212}C_{0}C_{x_{\beta}w_{0}}-C_{121212}C_{2121201210121212}\\
&=&[2]C_{121212}C_{02121021210121212}-([2]^{6}-4[2]^{4}+3[2]^{2})C_{2121201210121212}\\
&=&[2]^{2}C_{x_{\alpha}x_{\beta}w_{0}}+([2]^{6}-3[2]^{4}+[2]^{2})C_{x_{\beta}w_{0}}+[2]^{2}C_{x_{\alpha}w_{0}}\\
& &+[2]^{2}C_{x_{\alpha}^{2}w_{0}}-([2]^{6}-4[2]^{4}+3[2]^{2})C_{x_{\beta}w_{0}}\\
&=&[2]^{2}C_{x_{\alpha}x_{\beta}w_{0}}+([2]^{4}-2[2]^{2})C_{x_{\beta}w_{0}}+[2]^{2}C_{x_{\alpha}w_{0}}+[2]^{2}C_{x_{\alpha}^{2}w_{0}}\\
\end{eqnarray*}

(48) Computing $C_{121212012}C_{0121201210121212}$
\begin{eqnarray*}
C_{2}C_{0121201210121212}&=&C_{02121201210121212}+C_{021201210121212}+C_{0121212}\\
&=&C_{02121201210121212}+C_{210121210121212}+C_{0121212},\\
C_{121212012}C_{0121201210121212}&=&(C_{12121201}C_{2}-C_{1212120})C_{0121201210121212}\\
&=&C_{12121201}C_{2}C_{0121201210121212}-C_{1212120}C_{0121201210121212}\\
&=&C_{12121201}C_{02121201210121212}+C_{12121201}C_{210121210121212}\\
& &+C_{12121201}C_{0121212}-C_{1212120}C_{0121201210121212}\\
&=&[2]^{3}C_{x_{\alpha}x_{\beta}w_{0}}+(2[2]^{3}-[2])C_{x_{\beta}w_{0}}+([2]^{5}+[2])C_{x_{\alpha}w_{0}}\\
& &+2[2]^{3}C_{x_{\alpha}^{2}w_{0}}+(2[2]^{3}-[2])C_{w_{0}}\\
\end{eqnarray*}

\begin{eqnarray*}
C_{02121201210121212}&=&C_{0}C_{2121201210121212},\\
C_{12121201}C_{02121201210121212}&=&C_{12121201}C_{0}C_{2121201210121212}\\
&=&[2]C_{12121201}C_{x_{\beta}w_{0}}\\
&=&[2]^{3}C_{x_{\alpha}x_{\beta}w_{0}}+([2]^{5}-2[2]^{3})C_{x_{\beta}w_{0}}+[2]^{3}C_{x_{\alpha}w_{0}}+[2]^{3}C_{x_{\alpha}^{2}w_{0}}.\\
\end{eqnarray*}

\begin{eqnarray*}
C_{12121201}C_{210121210121212}&=&(C_{1212120}C_{1}-C_{121212})C_{210121210121212}\\
&=&C_{1212120}C_{1210121210121212}+C_{1212120}C_{10121210121212}\\
& &-C_{121212}C_{210121210121212},\\
&=&C_{121212}C_{01210121210121212}+C_{121212}C_{210121210121212}\\
& &+[2]C_{121212}C_{10121210121212}-C_{121212}C_{210121210121212}\\
&=&C_{121212}C_{01210121210121212}+[2]C_{121212}C_{10121210121212}\\
&=&[2]C_{x_{\alpha}x_{\beta}w_{0}}+(2[2]^{3}-[2])C_{w_{0}}+([2]^{5}+[2])C_{x_{\alpha}w_{0}}\\
& &+(2[2]^{3}-[2])C_{x_{\beta}w_{0}}+2[2]^{3}C_{x_{\alpha}^{2}w_{0}}.\\
\end{eqnarray*}

(49) Computing $C_{1212120121}C_{121212}$
\begin{eqnarray*}
C_{1212120121}C_{121212}&=&C_{121212}C_{1210121212}=([2]^{4}-2[2]^{2})C_{x_{\alpha}w_{0}}+[2]^{2}C_{w_{0}}.
\end{eqnarray*}

(50) Computing $C_{1212120121}C_{0121212}$
\begin{eqnarray*}
C_{1212120121}C_{0121212}&=&C_{1212120}C_{1210121212}=[2]C_{x_{\beta}w_{0}}+2([2]^{3}-[2])C_{x_{\alpha}w_{0}}+[2]^{3}C_{w_{0}}.\\
\end{eqnarray*}

(51) Computing $C_{1212120121}C_{10121212}$
\begin{eqnarray*}
C_{1212120121}C_{10121212}&=&C_{12121201}C_{1210121212}=[2]^{4}C_{x_{\alpha}w_{0}}+[2]^{2}C_{x_{\beta}w_{0}}+([2]^{4}-[2]^{2})C_{w_{0}}.\\
\end{eqnarray*}

(52) Computing $C_{1212120121}C_{210121212}$
\begin{eqnarray*}
C_{1212120121}C_{210121212}&=&C_{121212012}C_{1210121212}=[2]C_{x_{\alpha}^{2}w_{0}}+2[2]C_{x_{\beta}w_{0}}\\
& &+2[2]^{3}C_{x_{\alpha}w_{0}}+([2]^{5}-2[2]^{3}+[2])C_{w_{0}}.\\
\end{eqnarray*}

(53) Computing $C_{1212120121}C_{1210121212}$
\begin{eqnarray*}
&&C_{121212012}C_{1}=C_{1212120121}+C_{12121201},\\
&&C_{1212120121}=C_{121212012}C_{1}-C_{12121201},\\
C_{1212120121}C_{1210121212}&=&(C_{121212012}C_{1}-C_{12121201})C_{1210121212}\\
&=&C_{121212012}C_{1}C_{1210121212}-C_{12121201}C_{1210121212}\\
&=&[2]C_{121212012}C_{1210121212}-C_{12121201}C_{1210121212}\\
&=&[2]^{2}C_{x_{\alpha}^{2}w_{0}}+[2]^{4}C_{x_{\alpha}w_{0}}+[2]^{4}C_{x_{\beta}w_{0}}+([2]^{6}-3[2]^{4}+2[2]^{2})C_{w_{0}}.\\
\end{eqnarray*}

(54) Computing $C_{1212120121}C_{21210121212}$
\begin{eqnarray*}
C_{1}C_{21210121212}&=&C_{x_{\alpha}121212}+C_{1210121212},\\
C_{121212012}C_{x_{\alpha}121212}&=&(C_{121212}C_{0}C_{1}C_{2}-[2]C_{121212}-C_{1212120})C_{x_{\alpha}121212}\\
&=&([2]^{2}-1)C_{121212}C_{0121210121212}-[2]C_{121212}C_{x_{\alpha}121212}\\
&=&([2]^{3}-[2])C_{x_{\alpha}^{2}w_{0}}+2([2]^{3}-[2])C_{x_{\alpha}w_{0}}+([2]^{3}-[2])C_{x_{\beta}w_{0}}\\
& &+([2]^{3}-[2])C_{w_{0}}\\
\end{eqnarray*}

\begin{eqnarray*}
C_{1212120121}C_{21210121212}&=&(C_{121212012}C_{1}-C_{12121201})C_{21210121212}\\
&=&C_{121212012}C_{1}C_{21210121212}-C_{12121201}C_{21210121212}\\
&=&C_{121212012}C_{x_{\alpha}121212}+C_{121212012}C_{1210121212}-C_{12121201}C_{21210121212}\\
&=&([2]^{3}-[2])C_{x_{\alpha}^{2}w_{0}}+(2[2]^{3}-[2])C_{x_{\alpha}w_{0}}+2([2]^{3}-[2])C_{x_{\beta}w_{0}}\\
& &+([2]^{5}-2[2]^{3})C_{w_{0}}\\
\end{eqnarray*}

(55) Computing $C_{1212120121}C_{01210121212}$
\begin{eqnarray*}
C_{1212120121}C_{01210121212}&=&(C_{121212012}C_{1}-C_{12121201})C_{01210121212}\\
&=&C_{121212012}C_{1}C_{01210121212}-C_{12121201}C_{01210121212}\\
&=&[2]C_{121212012}C_{01210121212}-C_{12121201}C_{01210121212}\\
&=&[2]^{3}C_{x_{\alpha}^{2}w_{0}}+(2[2]^{3}-[2])C_{x_{\alpha}w_{0}}+[2]C_{x_{\beta}w_{0}}\\
& &+([2]^{5}-2[2]^{3}+[2])C_{w_{0}}\\
\end{eqnarray*}

(56) Computing $C_{1212120121}C_{201210121212}$
\begin{eqnarray*}
C_{1}C_{201210121212}&=&C_{1201210121212}+C_{01210121212},\\
C_{1212120121}C_{201210121212}&=&(C_{121212012}C_{1}-C_{12121201})C_{201210121212}\\
&=&C_{121212012}C_{1}C_{201210121212}-C_{12121201}C_{201210121212}\\
&=&C_{121212012}C_{1201210121212}+C_{121212012}C_{01210121212}\\
& &-C_{12121201}C_{201210121212}\\
&=&C_{x_{\alpha}x_{\beta}w_{0}}+3[2]^{2}C_{x_{\alpha}^{2}w_{0}}+2([2]^{4}-[2]^{2})C_{w_{0}}\\
& &+([2]^{4}+3[2]^{2})C_{x_{\alpha}w_{0}}+2[2]^{2}C_{x_{\beta}w_{0}}\\
\end{eqnarray*}

(57) Computing $C_{1212120121}C_{1201210121212}$
\begin{eqnarray*}
C_{1212120121}C_{1201210121212}&=&(C_{121212012}C_{1}-C_{12121201})C_{1201210121212},\\
&=&C_{121212012}C_{1}C_{1201210121212}-C_{12121201}C_{1201210121212}\\
&=&[2]C_{121212012}C_{1201210121212}-C_{12121201}C_{1201210121212}\\
&=&[2]C_{x_{\alpha}x_{\beta}w_{0}}+2[2]^{3}C_{x_{\alpha}^{2}w_{0}}+([2]^{5}-[2])C_{w_{0}}\\
& &+([2]^{5}+[2]^{3}+[2])C_{x_{\alpha}w_{0}}+(2[2]^{3}-[2])C_{x_{\beta}w_{0}}.\\
\end{eqnarray*}

(58) Computing $C_{1212120121}C_{21201210121212}$
\begin{eqnarray*}
C_{1}C_{21201210121212}&=&C_{121201210121212}+C_{1201210121212}.\\
C_{1212120121}C_{21201210121212}&=&(C_{121212012}C_{1}-C_{12121201})C_{21201210121212}\\
&=&C_{121212012}C_{1}C_{21201210121212}-C_{12121201}C_{21201210121212}\\
&=&C_{121212012}C_{121201210121212}+C_{121212012}C_{1201210121212}\\
& &-C_{12121201}C_{21201210121212}\\
&=&[2]^{2}C_{x_{\alpha}x_{\beta}w_{0}}+3[2]^{2}C_{x_{\alpha}^{2}w_{0}}+([2]^{4}+[2]^{2})C_{w_{0}}\\
& &+(2[2]^{4}+[2]^{2})C_{x_{\alpha}w_{0}}+3[2]^{2}C_{x_{\beta}w_{0}}.\\
\end{eqnarray*}

(59) Computing $C_{1212120121}C_{121201210121212}$
\begin{eqnarray*}
C_{1212120121}C_{121201210121212}&=&(C_{121212012}C_{1}-C_{12121201})C_{121201210121212}\\
&=&C_{121212012}C_{1}C_{121201210121212}-C_{12121201}C_{121201210121212}\\
&=&[2]C_{121212012}C_{121201210121212}-C_{12121201}C_{121201210121212}\\
&=&([2]^{3}-[2])C_{x_{\alpha}x_{\beta}w_{0}}+[2]^{3}C_{x_{\alpha}^{2}w_{0}}+([2]^{3}+[2])C_{w_{0}}\\
& &+([2]^{5}-[2])C_{x_{\alpha}w_{0}}+([2]^{3}+[2])C_{x_{\beta}w_{0}}.\\
\end{eqnarray*}

(60) Computing $C_{1212120121}C_{0121201210121212}$
\begin{eqnarray*}
C_{1}C_{0121201210121212}&=&C_{10121201210121212}+C_{121201210121212},\\
C_{1212120121}C_{0121201210121212}&=&(C_{121212012}C_{1}-C_{12121201})C_{0121201210121212}\\
&=&C_{121212012}(C_{10121201210121212}+C_{121201210121212})\\
& &-C_{12121201}C_{0121201210121212}\\
&=&[2]^{2}C_{x_{\alpha}^{3}w_{0}}+2[2]^{4}C_{x_{\alpha}w_{0}}+4[2]^{2}C_{x_{\beta}w_{0}}\\
& &+2[2]^{2}C_{x_{\alpha}x_{\beta}w_{0}}+4[2]^{2}C_{x_{\alpha}^{2}w_{0}}+3[2]^{2}C_{w_{0}}.\\
\end{eqnarray*}

\begin{eqnarray*}
C_{121212012}C_{10121201210121212}&=&(C_{12121201}C_{2}-C_{1212120})C_{10121201210121212}\\
&=&C_{12121201}C_{210121201210121212}-C_{1212120}C_{10121201210121212}\\
&=&C_{1212120}C_{1210121201210121212}+C_{1212120}C_{10121201210121212}\\
& &-C_{121212}C_{2101210121210121212}-C_{1212120}C_{10121201210121212}\\
&=&C_{1212120}C_{1210121201210121212}-C_{121212}C_{2101210121210121212}\\
&=&[2]^{2}C_{x_{\alpha}^{3}w_{0}}+(2[2]^{4}+[2]^{2})C_{x_{\alpha}w_{0}}+([2]^{4}+[2]^{2})C_{x_{\beta}w_{0}}\\
& &+2[2]^{2}C_{x_{\alpha}x_{\beta}w_{0}}+([2]^{4}+2[2]^{2})C_{x_{\alpha}^{2}w_{0}}+[2]^{4}C_{w_{0}}.\\
\end{eqnarray*}

\begin{eqnarray*}
C_{0}C_{1210121201210121212}&=&C_{01210121201210121212}+C_{210121201210121212}\\
& &+C_{102121021210121212}+2C_{10121210121212}+C_{10121212},\\
C_{1212120}C_{1210121201210121212}&=&C_{121212}C_{0}C_{1210121201210121212}\\
&=&C_{121212}(C_{01210121201210121212}+C_{210121201210121212}\\
& &+C_{102121021210121212}+2C_{10121210121212}+C_{10121212})\\
&=&[2]^{2}C_{x_{\alpha}^{3}w_{0}}+[2]^{2}C_{x_{\alpha}w_{0}}+[2]^{2}C_{x_{\beta}w_{0}}\\
& &+[2]^{2}C_{x_{\alpha}x_{\beta}w_{0}}+([2]^{4}-[2]^{2})C_{x_{\alpha}^{2}w_{0}}\\
& &+[2]^{2}C_{x_{\alpha}x_{\beta}w_{0}}+[2]^{4}C_{x_{\alpha}^{2}w_{0}}+([2]^{4}+[2]^{2})C_{x_{\alpha}w_{0}}\\
& &+([2]^{4}-[2]^{2})C_{x_{\beta}w_{0}}+([2]^{4}-[2]^{2})C_{w_{0}}\\
& &+[2]^{2}C_{x_{\alpha}x_{\beta}w_{0}}+[2]^{2}C_{x_{\alpha}^{2}w_{0}}+[2]^{2}C_{x_{\alpha}w_{0}}\\
& &+([2]^{4}-2[2]^{2})C_{x_{\beta}w_{0}}+2\{[2]^{2}C_{x_{\alpha}^{2}w_{0}}+[2]^{2}C_{x_{\beta}w_{0}}\\
& &+([2]^{4}-[2]^{2})C_{x_{\alpha}w_{0}}+[2]^{2}C_{w_{0}}\}\\
& &+[2]^{2}C_{x_{\alpha}w_{0}}+([2]^{4}-2[2]^{2})C_{w_{0}}\\
&=&[2]^{2}C_{x_{\alpha}^{3}w_{0}}+(3[2]^{4}+2[2]^{2})C_{x_{\alpha}w_{0}}+2[2]^{4}C_{x_{\beta}w_{0}}+3[2]^{2}C_{x_{\alpha}x_{\beta}w_{0}}\\
& &+(2[2]^{4}+2[2]^{2})C_{x_{\alpha}^{2}w_{0}}+(2[2]^{4}-[2]^{2})C_{w_{0}}.\\
\end{eqnarray*}

\begin{eqnarray*}
C_{121212}C_{102121021210121212}&=&C_{121212}(C_{1}C_{02121021210121212}-C_{x_{\beta}w_{0}})\\
&=&[2]C_{121212}C_{02121021210121212}-C_{121212}C_{x_{\beta}w_{0}}\\
&=&[2]^{2}C_{x_{\alpha}x_{\beta}w_{0}}+[2]^{2}C_{x_{\alpha}^{2}w_{0}}+[2]^{2}C_{x_{\alpha}w_{0}}+([2]^{4}-2[2]^{2})C_{x_{\beta}w_{0}}\\
\end{eqnarray*}

\begin{eqnarray*}
C_{121212}C_{210121201210121212}&=&C_{121212}C_{2}C_{10121201210121212}\\
&=&[2]C_{121212}C_{10121201210121212}\\
&=&[2]^{2}C_{x_{\alpha}x_{\beta}w_{0}}+[2]^{4}C_{x_{\alpha}^{2}w_{0}}+([2]^{4}+[2]^{2})C_{x_{\alpha}w_{0}}\\
& &+([2]^{4}-[2]^{2})C_{x_{\beta}w_{0}}+([2]^{4}-[2]^{2})C_{w_{0}}\\
\end{eqnarray*}

\begin{eqnarray*}
C_{121212}C_{01210121201210121212}&=&C_{121212}C_{10121210121210121212},\\
C_{2}C_{10121210121210121212}&=&C_{210121210121210121212}+C_{0121210121210121212},\\
C_{1}C_{210121210121210121212}&=&C_{1210121210121210121212}+C_{10121210121210121212},\\
C_{1}C_{0121210121210121212}&=&C_{10121210121210121212}+C_{x_{\alpha}^{2}w_{0}},\\
C_{1}C_{2}C_{0121210121210121212}&=&C_{1210121210121210121212}+2C_{10121210121210121212}\\
& &+C_{x_{\alpha}^{2}w_{0}},\\
C_{2}C_{1210121210121210121212}&=&C_{21210121210121210121212}+C_{210121210121210121212},\\
C_{2}C_{1}C_{2}C_{0121210121210121212}&=&C_{21210121210121210121212}+3C_{210121210121210121212}\\
& &+2C_{0121210121210121212}+[2]C_{x_{\alpha}^{2}w_{0}},\\
C_{1}C_{21210121210121210121212}&=&C_{x_{\alpha}^{3}w_{0}}+C_{x_{\alpha}w_{0}}+C_{x_{\alpha}^{2}w_{0}}\\
& &+C_{x_{\alpha}x_{\beta}w_{0}}+C_{x_{\beta}w_{0}}+C_{1210121210121210121212},\\
C_{1}C_{2}C_{1}C_{2}C_{0121210121210121212}&=&C_{x_{\alpha}^{3}w_{0}}+C_{x_{\alpha}w_{0}}+C_{x_{\beta}w_{0}}\\
& &+C_{x_{\alpha}x_{\beta}w_{0}}+([2]^{2}+3)C_{x_{\alpha}^{2}w_{0}}\\
& &+4C_{1210121210121210121212}+5C_{10121210121210121212},\\
C_{2}C_{1}C_{2}C_{1}C_{2}C_{0121210121210121212}&=&4C_{21210121210121210121212}+9C_{210121210121210121212}\\
& &+5C_{0121210121210121212}+[2]C_{x_{\alpha}^{3}w_{0}}+[2]C_{x_{\alpha}w_{0}}\\
& &+([2]^{3}+3[2])C_{x_{\alpha}^{2}w_{0}}+[2]C_{x_{\alpha}x_{\beta}w_{0}}+[2]C_{x_{\beta}w_{0}},\\
C_{1}C_{2}C_{1}C_{2}C_{1}C_{2}C_{0121210121210121212}&=&([2]^{2}+4)C_{x_{\alpha}^{3}w_{0}}+([2]^{2}+4)C_{x_{\alpha}w_{0}}\\
& &+([2]^{2}+4)C_{x_{\beta}w_{0}}+([2]^{2}+4)C_{x_{\alpha}x_{\beta}w_{0}}\\
& &+([2]^{4}+3[2]^{2}+9)C_{x_{\alpha}^{2}w_{0}}\\
& &+13C_{1210121210121210121212}+14C_{10121210121210121212},\\
C_{121212}C_{01210121201210121212}&=&C_{121212}C_{10121210121210121212}\\
&=&([2]^{2}+4)C_{x_{\alpha}^{3}w_{0}}+([2]^{2}+4)C_{x_{\alpha}w_{0}}\\
& &+([2]^{2}+4)C_{x_{\beta}w_{0}}+([2]^{2}+4)C_{x_{\alpha}x_{\beta}w_{0}}\\
& &+([2]^{4}+3[2]^{2}+9)C_{x_{\alpha}^{2}w_{0}}\\
& &+13C_{1210121210121210121212}+14C_{10121210121210121212}\\
& &-4\{C_{x_{\alpha}^{3}w_{0}}+C_{x_{\alpha}w_{0}}+C_{x_{\beta}w_{0}}\\
& &+C_{x_{\alpha}x_{\beta}w_{0}}+([2]^{2}+3)C_{x_{\alpha}^{2}w_{0}}\\
& &+4C_{1210121210121210121212}+5C_{10121210121210121212}\}\\
& &+3\{C_{1210121210121210121212}+2C_{10121210121210121212}\\
& &+C_{x_{\alpha}^{2}w_{0}}\}\\
&=&[2]^{2}C_{x_{\alpha}^{3}w_{0}}+[2]^{2}C_{x_{\alpha}w_{0}}+[2]^{2}C_{x_{\beta}w_{0}}+[2]^{2}C_{x_{\alpha}x_{\beta}w_{0}}\\
& &+([2]^{4}-[2]^{2})C_{x_{\alpha}^{2}w_{0}}\\
\end{eqnarray*}

(61) Computing $C_{12121201212}C_{121212}$
\begin{eqnarray*}
C_{12121201212}C_{121212}&=&C_{121212}C_{21210121212}=([2]^{5}-3[2]^{3}+[2])C_{x_{\alpha}w_{0}}+[2]C_{w_{0}}.
\end{eqnarray*}

(62) Computing $C_{12121201212}C_{0121212}$
\begin{eqnarray*}
C_{12121201212}C_{0121212}&=&C_{1212120}C_{21210121212}=[2]^{2}C_{x_{\beta}w_{0}}+([2]^{4}-[2]^{2})C_{x_{\alpha}w_{0}}+[2]^{2}C_{w_{0}}.
\end{eqnarray*}

(63) Computing $C_{12121201212}C_{10121212}$
\begin{eqnarray*}
C_{12121201212}C_{10121212}&=&C_{12121201}C_{21210121212}==[2]C_{x_{\alpha}^{2}w_{0}}+(2[2]^{3}-[2])C_{x_{\alpha}w_{0}}\\
& &+2[2]C_{x_{\beta}w_{0}}+[2]^{3}C_{w_{0}}.
\end{eqnarray*}

(64) Computing $C_{12121201212}C_{210121212}$
\begin{eqnarray*}
C_{12121201212}C_{210121212}&=&C_{121212012}C_{21210121212}==[2]^{2}C_{x_{\alpha}^{2}w_{0}}+[2]^{4}C_{x_{\alpha}w_{0}}\\
& &+[2]^{2}C_{x_{\beta}w_{0}}+([2]^{4}-[2]^{2})C_{w_{0}}.
\end{eqnarray*}

(65) Computing $C_{12121201212}C_{1210121212}$
\begin{eqnarray*}
C_{12121201212}C_{1210121212}&=&C_{12121201212}C_{1210121212}==([2]^{3}-[2])C_{x_{\alpha}^{2}w_{0}}+(2[2]^{3}-[2])C_{x_{\alpha}w_{0}}\\
& &+2([2]^{3}-[2])C_{x_{\beta}w_{0}}+([2]^{5}-2[2]^{3})C_{w_{0}}.
\end{eqnarray*}

(66) Computing $C_{12121201212}C_{21210121212}$
\begin{eqnarray*}
C_{1212120121}C_{2}&=&C_{12121201212}+C_{121212012},\\
C_{12121201212}&=&C_{1212120121}C_{2}-C_{121212012}.\\
C_{12121201212}C_{21210121212}&=&(C_{1212120121}C_{2}-C_{121212012})C_{21210121212}\\
&=&C_{1212120121}C_{2}C_{21210121212}-C_{121212012}C_{21210121212}\\
&=&[2]C_{1212120121}C_{21210121212}-C_{121212012}C_{21210121212}\\
&=&([2]^{4}-2[2]^{2})C_{x_{\alpha}^{2}w_{0}}+([2]^{4}-[2]^{2})C_{x_{\alpha}w_{0}}\\
& &+(2[2]^{4}-3[2]^{2})C_{x_{\beta}w_{0}}+([2]^{6}-3[2]^{4}+[2]^{2})C_{w_{0}}\\
\end{eqnarray*}

(67) Computing $C_{12121201212}C_{01210121212}$
\begin{eqnarray*}
C_{12121201212}C_{01210121212}&=&(C_{1212120121}C_{2}-C_{121212012})C_{01210121212}\\
&=&C_{1212120121}C_{2}C_{01210121212}-C_{121212012}C_{01210121212}\\
&=&C_{1212120121}C_{201210121212}-C_{121212012}C_{01210121212}\\
&=&C_{x_{\alpha}x_{\beta}w_{0}}+2[2]^{2}C_{x_{\alpha}^{2}w_{0}}+3[2]^{2}C_{x_{\alpha}w_{0}}+[2]^{2}C_{x_{\beta}w_{0}}\\
& &+([2]^{4}-[2]^{2})C_{w_{0}}.\\
\end{eqnarray*}

(68) Computing $C_{12121201212}C_{201210121212}$
\begin{eqnarray*}
C_{12121201212}C_{201210121212}&=&C_{12121201212}C_{2}C_{01210121212}\\
&=&[2]C_{12121201212}C_{01210121212}\\
&=&[2]C_{x_{\alpha}x_{\beta}w_{0}}+2[2]^{3}C_{x_{\alpha}^{2}w_{0}}+3[2]^{3}C_{x_{\alpha}w_{0}}+[2]^{3}C_{x_{\beta}w_{0}}\\\
& &+([2]^{5}-[2]^{3})C_{w_{0}}.\\
\end{eqnarray*}

(69) Computing $C_{12121201212}C_{1201210121212}$
\begin{eqnarray*}
C_{12121201212}C_{1201210121212}&=&(C_{1212120121}C_{2}-C_{121212012})C_{1201210121212}\\
&=&C_{1212120121}C_{2}C_{1201210121212}-C_{121212012}C_{1201210121212}\\
&=&C_{1212120121}C_{21201210121212}+C_{1212120121}C_{201210121212}\\
& &+C_{1212120121}C_{121210121212}-C_{121212012}C_{1201210121212}\\
&=&[2]^{2}C_{x_{\alpha}x_{\beta}w_{0}}+([2]^{4}+[2]^{2})C_{x_{\alpha}^{2}w_{0}}+(2[2]^{4}-2[2]^{2})C_{w_{0}}\\
& &+(2[2]^{4}+[2]^{2})C_{x_{\alpha}w_{0}}+[2]^{4}C_{x_{\beta}w_{0}}.\\
\end{eqnarray*}

\begin{eqnarray*}
C_{1212120121}C_{121210121212}&=&(C_{121212012}C_{1}-C_{12121201})C_{121210121212}\\
&=&C_{121212012}C_{1}C_{121210121212}-C_{12121201}C_{121210121212}\\
&=&[2]C_{121212012}C_{121210121212}-C_{12121201}C_{121210121212}\\
&=&([2]^{2}-1)C_{12121201}C_{121210121212}-[2]C_{121212}C_{0121210121212}\\
&=&([2]^{4}-2[2]^{2})C_{x_{\alpha}^{2}w_{0}}+([2]^{4}-[2]^{2})C_{x_{\alpha}w_{0}}\\
& &+([2]^{4}-2[2]^{2})C_{x_{\beta}w_{0}}+([2]^{4}-[2]^{2})C_{w_{0}}.\\
\end{eqnarray*}

(70) Computing $C_{12121201212}C_{21201210121212}$
\begin{eqnarray*}
C_{12121201212}C_{21201210121212}&=&(C_{1212120121}C_{2}-C_{121212012})C_{21201210121212}\\
&=&C_{1212120121}C_{2}C_{21201210121212}-C_{121212012}C_{21201210121212}\\
&=&[2]C_{1212120121}C_{21201210121212}-C_{121212012}C_{21201210121212}\\
&=&([2]^{3}-[2])C_{x_{\alpha}x_{\beta}w_{0}}+(2[2]^{3}-[2])C_{x_{\alpha}^{2}w_{0}}+(2[2]^{3}-[2])C_{w_{0}}\\
& &+([2]^{5}+[2]^{3}-2[2])C_{x_{\alpha}w_{0}}+(2[2]^{3}-[2])C_{x_{\beta}w_{0}}.\\
\end{eqnarray*}

(71) Computing $C_{12121201212}C_{121201210121212}$
\begin{eqnarray*}
C_{12121201212}C_{121201210121212}&=&(C_{1212120121}C_{2}-C_{121212012})C_{121201210121212}\\
&=&C_{1212120121}C_{2}C_{121201210121212}-C_{121212012}C_{121201210121212}\\
&=&C_{1212120121}C_{x_{\beta}121212}+C_{1212120121}C_{21201210121212}\\
& &+C_{1212120121}C_{121212}-C_{121212012}C_{121201210121212}\\
&=&([2]^{4}-2[2]^{2})C_{x_{\alpha}x_{\beta}w_{0}}+2[2]^{2}C_{x_{\beta}w_{0}}+2[2]^{2}C_{w_{0}}\\
& &+(3[2]^{4}-5[2]^{2})C_{x_{\alpha}w_{0}}+([2]^{4}-[2]^{2})C_{x_{\alpha}^{2}w_{0}}\\
\end{eqnarray*}

\begin{eqnarray*}
C_{1212120121}C_{2121021210121212}&=&(C_{121212012}C_{1}-C_{12121201})C_{2121021210121212}\\
&=&C_{121212012}C_{1}C_{2121021210121212}-C_{12121201}C_{2121021210121212}\\
&=&[2]C_{121212012}C_{2121021210121212}-C_{12121201}C_{2121021210121212}\\
&=&([2]^{2}-1)C_{12121201}C_{2121021210121212}-[2]C_{121212}C_{02121021210121212}\\
&=&([2]^{4}-2[2]^{2})C_{x_{\alpha}x_{\beta}w_{0}}+[2]^{2}C_{x_{\beta}w_{0}}+([2]^{4}-2[2]^{2})C_{x_{\alpha}w_{0}}\\
& &+([2]^{4}-2[2]^{2})C_{x_{\alpha}^{2}w_{0}}\\
\end{eqnarray*}

(72) Computing $C_{12121201212}C_{0121201210121212}$
\begin{eqnarray*}
C_{2}C_{0121201210121212}&=&C_{20121201210121212}+C_{210121210121212}+C_{0121212},\\
C_{12121201212}C_{0121201210121212}&=&(C_{1212120121}C_{2}-C_{121212012})C_{0121201210121212}\\
&=&C_{1212120121}C_{2}C_{0121201210121212}-C_{121212012}C_{0121201210121212}\\
&=&C_{1212120121}C_{20121201210121212}+C_{1212120121}C_{210121210121212}\\
& &+C_{1212120121}C_{0121212}-C_{121212012}C_{0121201210121212}\\
&=&2([2]^{3}-[2])C_{x_{\alpha}x_{\beta}w_{0}}+(2[2]^{3}-[2])C_{x_{\alpha}^{2}w_{0}}+2([2]^{3}-[2])C_{x_{\alpha}w_{0}}\\
& &+([2]^{3}+[2])C_{x_{\beta}w_{0}}+[2]C_{w_{0}}+[2]C_{x_{\alpha}^{3}w_{0}}\\
& &+[2]C_{x_{\alpha}^{3}w_{0}}+([2]^{5}+4[2])C_{x_{\alpha}w_{0}}+(2[2]^{3}+[2])C_{x_{\beta}w_{0}}\\
& &+3[2]C_{x_{\alpha}x_{\beta}w_{0}}+(2[2]^{3}+3[2])C_{x_{\alpha}^{2}w_{0}}+2[2]^{3}C_{w_{0}}\\
& &+[2]C_{x_{\beta}w_{0}}+2([2]^{3}-[2])C_{x_{\alpha}w_{0}}+[2]^{3}C_{w_{0}}\\
& &-\{[2]^{3}C_{x_{\alpha}x_{\beta}w_{0}}+(2[2]^{3}-[2])C_{x_{\beta}w_{0}}+([2]^{5}+[2])C_{x_{\alpha}w_{0}}\\
& &+2[2]^{3}C_{x_{\alpha}^{2}w_{0}}+(2[2]^{3}-[2])C_{w_{0}}\}\\
&=&([2]^{3}+[2])C_{x_{\alpha}x_{\beta}w_{0}}+(2[2]^{3}+[2])C_{x_{\alpha}^{2}w_{0}}+(4[2]^{3}-[2])C_{x_{\alpha}w_{0}}\\
& &+([2]^{3}+4[2])C_{x_{\beta}w_{0}}+([2]^{3}+2[2])C_{w_{0}}+2[2]C_{x_{\alpha}^{3}w_{0}}\\
\end{eqnarray*}

\begin{eqnarray*}
C_{1212120121}C_{20121201210121212}&=&(C_{121212012}C_{1}-C_{12121201})C_{20121201210121212}\\
&=&C_{121212012}C_{1}C_{20121201210121212}-C_{12121201}C_{20121201210121212}\\
&=&C_{121212012}C_{120121201210121212}+C_{121212012}C_{2121201210121212}\\
& &-C_{12121201}C_{20121201210121212}\\
&=&(C_{12121201}C_{2}-C_{1212120})C_{120121201210121212}\\
& &+C_{121212012}C_{2121201210121212}-C_{12121201}C_{20121201210121212}\\
&=&C_{12121201}C_{2}C_{120121201210121212}-C_{1212120}C_{120121201210121212}\\
& &+C_{121212012}C_{2121201210121212}-C_{12121201}C_{20121201210121212}\\
&=&C_{12121201}C_{2120121201210121212}+C_{12121201}C_{20121201210121212}\\
& &-C_{1212120}C_{120121201210121212}+C_{121212012}C_{2121201210121212}\\
& &-C_{12121201}C_{20121201210121212}\\
&=&C_{12121201}C_{2120121201210121212}-C_{1212120}C_{120121201210121212}\\
& &+C_{121212012}C_{2121201210121212}\\
&=&(C_{1212120}C_{1}-C_{121212})C_{2120121201210121212}\\
& &-C_{1212120}C_{120121201210121212}+C_{121212012}C_{2121201210121212}\\
&=&C_{1212120}C_{12120121201210121212}+C_{1212120}C_{120121201210121212}\\
& &-C_{121212}C_{2120121201210121212}-C_{1212120}C_{120121201210121212}\\
& &+C_{121212012}C_{2121201210121212}\\
&=&C_{1212120}C_{12120121201210121212}-C_{121212}C_{2120121201210121212}\\
& &+C_{121212012}C_{2121201210121212}\\
&=&C_{121212}C_{0}C_{12120121201210121212}-C_{121212}C_{2120121201210121212}\\
& &+C_{121212012}C_{2121201210121212}\\
&=&C_{121212}C_{012120121201210121212}+C_{121212}C_{2120121201210121212}\\
& &-C_{121212}C_{2120121201210121212}+C_{121212012}C_{2121201210121212}\\
&=&C_{121212}C_{012120121201210121212}+C_{121212012}C_{2121201210121212}\\
&=&[2]C_{x_{\beta}^{2}w_{0}}+[2]C_{w_{0}}+2[2]C_{x_{\beta}w_{0}}+[2]^{3}C_{x_{\alpha}^{2}w_{0}}\\
& &+[2]C_{x_{\alpha}^{3}w_{0}}+([2]^{3}-[2])C_{x_{\alpha}w_{0}}+([2]^{3}-[2])C_{x_{\alpha}x_{\beta}w_{0}}\\
& &+([2]^{3}-[2])C_{x_{\alpha}x_{\beta}w_{0}}+([2]^{3}-[2])C_{x_{\alpha}^{2}w_{0}}+([2]^{3}-[2])C_{x_{\alpha}w_{0}}\\
& &+([2]^{3}-[2])C_{x_{\beta}w_{0}}\\
&=&[2]C_{x_{\beta}^{2}w_{0}}+2([2]^{3}-[2])C_{x_{\alpha}x_{\beta}w_{0}}+(2[2]^{3}-[2])C_{x_{\alpha}^{2}w_{0}}\\
& &+2([2]^{3}-[2])C_{x_{\alpha}w_{0}}+([2]^{3}+[2])C_{x_{\beta}w_{0}}+[2]C_{w_{0}}+[2]C_{x_{\alpha}^{3}w_{0}}\\
\end{eqnarray*}

\begin{eqnarray*}
C_{121212}C_{012120121201210121212}&=&C_{121212}C_{012102121021210121212},\\
C_{2}C_{012120121201210121212}&=&C_{0212102121021210121212},\\
C_{1}C_{0212102121021210121212}&=&C_{10212102121021210121212}+C_{012102121021210121212},\\
C_{1}C_{2}C_{012120121201210121212}&=&C_{10212102121021210121212}+C_{012102121021210121212},\\
C_{2}C_{10212102121021210121212}&=&C_{210212102121021210121212}+C_{x_{\alpha}w_{0}}+C_{x_{\alpha}^{2}w_{0}}\\
& &+C_{x_{\alpha}x_{\beta}w_{0}}+C_{0212102121021210121212},\\
C_{2}C_{1}C_{2}C_{012102121021210121212}&=&C_{210212102121021210121212}+C_{x_{\alpha}w_{0}}+C_{x_{\alpha}^{2}w_{0}}\\
& &+C_{x_{\alpha}x_{\beta}w_{0}}+2C_{0212102121021210121212},\\
C_{1}C_{210212102121021210121212}&=&C_{1210212102121021210121212}+C_{10212102121021210121212},\\
C_{1}C_{2}C_{1}C_{2}C_{012120121201210121212}&=&C_{1210212102121021210121212}+C_{10212102121021210121212}\\
& &+[2]C_{x_{\alpha}w_{0}}+[2]C_{x_{\alpha}^{2}w_{0}}+[2]C_{x_{\alpha}x_{\beta}w_{0}}\\
C_{2}C_{1210212102121021210121212}&=&C_{x_{\beta}^{2}w_{0}}+C_{w_{0}}+2C_{x_{\beta}w_{0}}\\
& &+C_{x_{\alpha}^{2}w_{0}}++C_{x_{\alpha}^{3}w_{0}}+C_{210212102121021210121212},\\
C_{2}C_{1}C_{2}C_{1}C_{2}C_{012120121201210121212}&=&C_{x_{\beta}^{2}w_{0}}+C_{w_{0}}+2C_{x_{\beta}w_{0}}+5C_{0212102121021210121212}\\
& &+([2]^{2}+4)C_{x_{\alpha}^{2}w_{0}}+4C_{210212102121021210121212}\\
& &+([2]^{2}+3)C_{x_{\alpha}^{2}w_{0}}+([2]^{2}+3)C_{x_{\alpha}x_{\beta}w_{0}}+C_{x_{\alpha}^{3}w_{0}}\\
C_{1}C_{2}C_{1}C_{2}C_{1}C_{2}C_{012120121201210121212}&=&[2]C_{x_{\beta}^{2}w_{0}}+[2]C_{w_{0}}+2[2]C_{x_{\beta}w_{0}}+([2]^{3}+4[2])C_{x_{\alpha}^{2}w_{0}}\\
& &+[2]C_{x_{\alpha}^{3}w_{0}}+([2]^{3}+3[2])C_{x_{\alpha}^{2}w_{0}}+([2]^{3}+3[2])C_{x_{\alpha}x_{\beta}w_{0}}\\
& &+4C_{1210212102121021210121212}+9C_{10212102121021210121212}\\
& &+5C_{012102121021210121212}\\
C_{121212}C_{012120121201210121212}&=&[2]C_{x_{\beta}^{2}w_{0}}+[2]C_{w_{0}}+2[2]C_{x_{\beta}w_{0}}+([2]^{3}+4[2])C_{x_{\alpha}^{2}w_{0}}\\
& &+[2]C_{x_{\alpha}^{3}w_{0}}+([2]^{3}+3[2])C_{x_{\alpha}^{2}w_{0}}+([2]^{3}+3[2])C_{x_{\alpha}x_{\beta}w_{0}}\\
& &+4C_{1210212102121021210121212}+9C_{10212102121021210121212}\\
& &-4\{C_{1210212102121021210121212}+C_{10212102121021210121212}\\
& &+[2]C_{x_{\alpha}w_{0}}+[2]C_{x_{\alpha}^{2}w_{0}}+[2]C_{x_{\alpha}x_{\beta}w_{0}}\}\\
& &+3\{C_{10212102121021210121212}+C_{012102121021210121212}\}\\
&=&[2]C_{x_{\beta}^{2}w_{0}}+[2]C_{w_{0}}+2[2]C_{x_{\beta}w_{0}}+[2]^{3}C_{x_{\alpha}^{2}w_{0}}\\
& &+[2]C_{x_{\alpha}^{3}w_{0}}+([2]^{3}-[2])C_{x_{\alpha}w_{0}}+([2]^{3}-[2])C_{x_{\alpha}x_{\beta}w_{0}}.\\
\end{eqnarray*}

\begin{eqnarray*}
C_{121212012}C_{2121021210121212}&=&(C_{12121201}C_{2}-C_{1212120})C_{2121021210121212}\\
&=&C_{12121201}C_{2}C_{2121021210121212}-C_{1212120}C_{2121021210121212}\\
&=&[2]C_{12121201}C_{2121021210121212}-C_{1212120}C_{2121021210121212}\\
&=&([2]^{3}-[2])C_{x_{\alpha}x_{\beta}w_{0}}+([2]^{3}-[2])C_{x_{\alpha}^{2}w_{0}}+([2]^{3}-[2])C_{x_{\alpha}w_{0}}\\
& &+([2]^{3}-[2])C_{x_{\beta}w_{0}}.\\
\end{eqnarray*}
\begin{eqnarray*}
C_{12121201}C_{0121210121212}&=&(C_{1212120}C_{1}-C_{121212})C_{0121210121212}\\
&=&C_{1212120}C_{10121210121212}+C_{1212120}C_{121210121212}\\
& &-C_{121212}C_{0121210121212}\\
&=&C_{1212120}C_{10121210121212}.\\
\end{eqnarray*}

\begin{eqnarray*}
C_{1212120121}C_{210121210121212}&=&(C_{121212012}C_{1}-C_{12121201})C_{210121210121212}\\
&=&C_{121212012}C_{1}C_{210121210121212}-C_{12121201}C_{210121210121212}\\
&=&C_{121212012}C_{1210121210121212}+C_{121212012}C_{10121210121212}\\
& &-C_{12121201}C_{210121210121212}\\
&=&C_{121212012}C_{1210121210121212}-C_{12121201}C_{210121210121212}\\
& &+C_{12121201}C_{210121210121212}+C_{12121201}C_{0121210121212}\\
& &-C_{1212120}C_{10121210121212}\\
&=&C_{121212012}C_{1210121210121212}
\end{eqnarray*}

\begin{eqnarray*}
C_{121212012}C_{1210121210121212}&=&(C_{12121201}C_{2}-C_{1212120})C_{1210121210121212}\\
&=&C_{12121201}C_{2}C_{1210121210121212}-C_{1212120}C_{1210121210121212}\\
&=&C_{12121201}C_{21210121210121212}+C_{12121201}C_{210121210121212}\\
& &-C_{1212120}C_{1210121210121212}\\
&=&C_{12121201}C_{21210121210121212}-C_{1212120}C_{1210121210121212}\\
& &+C_{1212120}C_{1210121210121212}+C_{1212120}C_{10121210121212}\\
& &-C_{121212}C_{210121210121212}\\
&=&C_{12121201}C_{21210121210121212}+C_{1212120}C_{10121210121212}\\
& &-C_{121212}C_{210121210121212}\\
&=&C_{1212120}(C_{x_{\alpha}^{2}w_{0}}+C_{x_{\alpha}w_{0}}+C_{x_{\beta}w_{0}}\\
& &+C_{w_{0}})+C_{121212}C_{01210121210121212}\\
& &+C_{121212}C_{210121210121212}-C_{121212}C_{21210121210121212}\\
& &+C_{1212120}C_{10121210121212}-C_{121212}C_{210121210121212}\\
&=&C_{1212120}(C_{x_{\alpha}^{2}w_{0}}+C_{x_{\alpha}w_{0}}+C_{x_{\beta}w_{0}}\\
& &+C_{w_{0}})+C_{121212}C_{01210121210121212}\\
& &-C_{121212}C_{21210121210121212}+C_{1212120}C_{10121210121212}\\
&=&C_{1212120}(C_{x_{\alpha}^{2}w_{0}}+C_{x_{\alpha}w_{0}}+C_{x_{\beta}w_{0}}\\
& &+C_{w_{0}})+C_{121212}C_{01210121210121212}-C_{121212}C_{21210121210121212}\\
& &+[2]C_{121212}C_{10121210121212}\\
&=&[2]C_{x_{\alpha}^{3}w_{0}}+[2]C_{x_{\alpha}w_{0}}+[2]C_{x_{\beta}w_{0}}+[2]C_{x_{\alpha}x_{\beta}w_{0}}\\
& &+([2]^{5}-3[2]^{3}+2[2])C_{x_{\alpha}^{2}w_{0}}+[2]C_{x_{\alpha}121212}\\
& &+([2]^{5}-3[2]^{3}+[2])C_{121212}+[2]C_{x_{\alpha}^{2}w_{0}}\\
& &+([2]^{5}-3[2]^{3}+2[2])C_{x_{\alpha}w_{0}}+[2]C_{x_{\beta}w_{0}}+[2]C_{w_{0}}\\
& &+[2]C_{x_{\alpha}x_{\beta}w_{0}}+([2]^{5}-3[2]^{3}+[2])C_{x_{\beta}w_{0}}+[2]C_{x_{\alpha}w_{0}}\\
& &+[2]C_{x_{\alpha}^{2}w_{0}}+[2]C_{x_{\alpha}x_{\beta}w_{0}}+([2]^{3}-[2])C_{w_{0}}+([2]^{3}+[2])C_{x_{\alpha}w_{0}}\\
& &+([2]^{3}-[2])C_{x_{\beta}w_{0}}+[2]^{3}C_{x_{\alpha}^{2}w_{0}}\\
& &+[2]^{3}C_{x_{\alpha}^{2}w_{0}}+([2]^{5}-[2]^{3})C_{x_{\alpha}w_{0}}+[2]^{3}C_{x_{\beta}w_{0}}+[2]^{3}C_{w_{0}}\\
& &-\{([2]^{5}-3[2]^{3}+[2])C_{x_{\alpha}^{2}w_{0}}+([2]^{5}-3[2]^{3}+2[2])C_{x_{\alpha}w_{0}}\\
& &+([2]^{5}-3[2]^{3}+[2])C_{x_{\beta}w_{0}}+([2]^{5}-3[2]^{3}+[2])C_{w_{0}}\}\\
&=&[2]C_{x_{\alpha}^{3}w_{0}}+([2]^{5}+4[2])C_{x_{\alpha}w_{0}}+(2[2]^{3}+[2])C_{x_{\beta}w_{0}}\\
& &+3[2]C_{x_{\alpha}x_{\beta}w_{0}}+(2[2]^{3}+3[2])C_{x_{\alpha}^{2}w_{0}}+2[2]^{3}C_{w_{0}}\\
\end{eqnarray*}

\begin{eqnarray*}
C_{1212120}C_{1210121210121212}&=&C_{121212}C_{01210121210121212}+C_{121212}C_{210121210121212}\\
&=&[2]C_{x_{\alpha}x_{\beta}w_{0}}+([2]^{3}-[2])C_{w_{0}}+([2]^{3}+[2])C_{x_{\alpha}w_{0}}\\
& &+([2]^{3}-[2])C_{x_{\beta}w_{0}}+[2]^{3}C_{x_{\alpha}^{2}w_{0}}+([2]^{3}-[2])C_{x_{\alpha}^{2}w_{0}}\\
& &+2([2]^{3}-[2])C_{x_{\alpha}w_{0}}+([2]^{3}-[2])C_{x_{\beta}w_{0}}+([2]^{3}-[2])C_{w_{0}}\\
&=&[2]C_{x_{\alpha}x_{\beta}w_{0}}+2([2]^{3}-[2])C_{w_{0}}+(3[2]^{3}-[2])C_{x_{\alpha}w_{0}}\\
& &+2([2]^{3}-[2])C_{x_{\beta}w_{0}}+(2[2]^{3}-[2])C_{x_{\alpha}^{2}w_{0}}\\
\end{eqnarray*}

\begin{eqnarray*}
C_{1212120}C_{x_{\alpha}^{2}w_{0}}&=&C_{121212}C_{0121210121210121212},\\
C_{121212}C_{10121210121210121212}&=&C_{121212}(C_{1}C_{0121210121210121212}-C_{121210121210121212})\\
&=&[2]C_{121212}C_{0121210121210121212}-C_{121212}C_{121210121210121212}),\\
C_{121212}C_{0121210121210121212}[2]&=&C_{121212}C_{10121210121210121212}+C_{121212}C_{121210121210121212}),\\
C_{121212}C_{0121210121210121212}[2]&=&[2]^{2}C_{x_{\alpha}^{3}w_{0}}+[2]^{2}C_{x_{\alpha}w_{0}}+[2]^{2}C_{x_{\beta}w_{0}}+[2]^{2}C_{x_{\alpha}x_{\beta}w_{0}}\\
& &+([2]^{4}-[2]^{2})C_{x_{\alpha}^{2}w_{0}}+([2]^{6}-4[2]^{4}+3[2]^{2})C_{x_{\alpha}^{2}w_{0}},\\
C_{121212}C_{0121210121210121212}&=&[2]C_{x_{\alpha}^{3}w_{0}}+[2]C_{x_{\alpha}w_{0}}+[2]C_{x_{\beta}w_{0}}+[2]C_{x_{\alpha}x_{\beta}w_{0}}\\
& &+([2]^{5}-3[2]^{3}+2[2])C_{x_{\alpha}^{2}w_{0}}.\\
\end{eqnarray*}

\begin{eqnarray*}
C_{121212}C_{21210121210121212}&=&C_{121212}(C_{2}C_{1210121210121212}-C_{210121210121212})\\
&=&[2]C_{121212}C_{1210121210121212}-C_{121212}C_{210121210121212}\\
&=&[2]C_{121212}(C_{1}C_{210121210121212}-C_{10121210121212})\\
& &-C_{121212}C_{210121210121212}\\
&=&([2]^{2}-1)C_{121212}C_{210121210121212}-[2]C_{121212}C_{10121210121212}\\
&=&([2]^{5}-2[2]^{3}+[2])C_{x_{\alpha}^{2}w_{0}}+2([2]^{5}-2[2]^{3}+[2])C_{x_{\alpha}w_{0}}\\
& &+([2]^{5}-2[2]^{3}+[2])C_{x_{\beta}w_{0}}+([2]^{5}-2[2]^{3}+[2])C_{w_{0}}\\
& &-\{[2]^{3}C_{x_{\alpha}^{2}w_{0}}+([2]^{5}-[2]^{3})C_{x_{\alpha}w_{0}}+[2]^{3}C_{x_{\beta}w_{0}}+[2]^{3}C_{w_{0}}\}\\
&=&([2]^{5}-3[2]^{3}+[2])C_{x_{\alpha}^{2}w_{0}}+([2]^{5}-3[2]^{3}+2[2])C_{x_{\alpha}w_{0}}\\
& &+([2]^{5}-3[2]^{3}+[2])C_{x_{\beta}w_{0}}+([2]^{5}-3[2]^{3}+[2])C_{w_{0}}\\
\end{eqnarray*}

(73) Computing $C_{12121201210}C_{121212}$
\begin{eqnarray*}
C_{12121201210}C_{121212}&=&C_{121212}C_{01210121212}=[2]C_{x_{\beta}w_{0}}+([2]^{3}-[2])C_{x_{\alpha}w_{0}}+[2]C_{w_{0}}.
\end{eqnarray*}

(74) Computing $C_{12121201210}C_{0121212}$
\begin{eqnarray*}
C_{12121201210}C_{0121212}&=&C_{1212120}C_{01210121212}=[2]^{2}C_{x_{\beta}w_{0}}+([2]^{4}-[2]^{2})C_{x_{\alpha}w_{0}}+[2]^{2}C_{w_{0}}.
\end{eqnarray*}

(75) Computing $C_{12121201210}C_{10121212}$
\begin{eqnarray*}
C_{12121201210}C_{10121212}&=&C_{12121201}C_{01210121212}=([2]^{5}-2[2]^{3}+[2])C_{x_{\alpha}w_{0}}\\
& &+([2]^{3}-[2])C_{x_{\beta}w_{0}}+([2]^{3}-[2])C_{w_{0}}.
\end{eqnarray*}

(76) Computing $C_{12121201210}C_{210121212}$
\begin{eqnarray*}
C_{12121201210}C_{210121212}&=&C_{121212012}C_{01210121212}=[2]^{2}C_{x_{\alpha}^{2}w_{0}}+[2]^{4}C_{x_{\alpha}w_{0}}\\
& &+[2]^{2}C_{x_{\beta}w_{0}}+([2]^{4}-[2]^{2})C_{w_{0}}.
\end{eqnarray*}

(77) Computing $C_{12121201210}C_{1210121212}$
\begin{eqnarray*}
C_{12121201210}C_{1210121212}&=&C_{1212120121}C_{01210121212}=[2]^{3}C_{x_{\alpha}^{2}w_{0}}+(2[2]^{3}-[2])C_{x_{\alpha}w_{0}}\\
& &+[2]C_{x_{\beta}w_{0}}+([2]^{5}-2[2]^{3}+[2])C_{w_{0}}.
\end{eqnarray*}

(78) Computing $C_{12121201210}C_{21210121212}$
\begin{eqnarray*}
C_{12121201210}C_{21210121212}&=&C_{12121201212}C_{01210121212}=C_{x_{\alpha}x_{\beta}w_{0}}+2[2]^{2}C_{x_{\alpha}^{2}w_{0}}+3[2]^{2}C_{x_{\alpha}w_{0}}\\\
& &+[2]^{2}C_{x_{\beta}w_{0}}+([2]^{4}-[2]^{2})C_{w_{0}}.\\
\end{eqnarray*}

(79) Computing $C_{12121201210}C_{01210121212}$
\begin{eqnarray*}
C_{1212120121}C_{0}&=&C_{12121201210}+C_{121212012},\\
C_{12121201210}C_{01210121212}&=&C_{1212120121}C_{0}C_{01210121212}\\
&=&C_{1212120121}C_{0}C_{01210121212}-C_{121212012}C_{01210121212}\\
&=&[2]C_{1212120121}C_{01210121212}-C_{121212012}C_{01210121212}\\
&=&[2]^{4}C_{x_{\alpha}^{2}w_{0}}+(2[2]^{4}-[2]^{2})C_{x_{\alpha}w_{0}}+[2]^{2}C_{x_{\beta}w_{0}}\\
& &+([2]^{6}-2[2]^{4}+[2]^{2})C_{w_{0}}-\{[2]^{2}C_{x_{\alpha}^{2}w_{0}}+[2]^{4}C_{x_{\alpha}w_{0}}\\
& &+[2]^{2}C_{x_{\beta}w_{0}}+([2]^{4}-[2]^{2})C_{w_{0}}\}\\
&=&([2]^{4}-[2]^{2})C_{x_{\alpha}^{2}w_{0}}+([2]^{4}-[2]^{2})C_{x_{\alpha}w_{0}}+([2]^{6}-3[2]^{4}+2[2]^{2})C_{w_{0}}.\\
\end{eqnarray*}

(80) Computing $C_{12121201210}C_{201210121212}$
\begin{eqnarray*}
C_{12121201210}C_{201210121212}&=&C_{1212120121}C_{0}C_{201210121212}\\
&=&C_{1212120121}C_{0}C_{201210121212}-C_{121212012}C_{201210121212}\\
&=&[2]C_{1212120121}C_{201210121212}-C_{121212012}C_{201210121212}\\
&=&[2]C_{x_{\alpha}x_{\beta}w_{0}}+3[2]^{3}C_{x_{\alpha}^{2}w_{0}}+2([2]^{5}-[2]^{3})C_{w_{0}}\\
& &+([2]^{5}+3[2]^{3})C_{x_{\alpha}w_{0}}+2[2]^{3}C_{x_{\beta}w_{0}}-\{[2]^{3}C_{x_{\alpha}^{2}w_{0}}\\
& &+[2]^{5}C_{x_{\alpha}w_{0}}+[2]^{3}C_{x_{\beta}w_{0}}+([2]^{5}-[2]^{3})C_{w_{0}}\}\\
&=&[2]C_{x_{\alpha}x_{\beta}w_{0}}+2[2]^{3}C_{x_{\alpha}^{2}w_{0}}+([2]^{5}-[2]^{3})C_{w_{0}}\\
& &+3[2]^{3}C_{x_{\alpha}w_{0}}+[2]^{3}C_{x_{\beta}w_{0}}.\\
\end{eqnarray*}

(81) Computing $C_{12121201210}C_{1201210121212}$
\begin{eqnarray*}
C_{0}C_{1201210121212}&=&C_{10121210121212}+C_{201210121212}+C_{10121212},\\
C_{12121201210}C_{1201210121212}&=&C_{1212120121}C_{0}C_{1201210121212}-C_{121212012}C_{1201210121212}\\
&=&C_{1212120121}C_{0}C_{1201210121212}-C_{121212012}C_{1201210121212}\\
&=&C_{1212120121}(C_{10121210121212}+C_{201210121212}+C_{10121212})\\
& &-C_{121212012}C_{1201210121212}\\
&=&[2]^{2}C_{x_{\alpha}x_{\beta}w_{0}}+[2]^{4}C_{w_{0}}+2[2]^{4}C_{x_{\alpha}w_{0}}\\
& &+[2]^{4}C_{x_{\beta}w_{0}}+([2]^{4}+[2]^{2})C_{x_{\alpha}^{2}w_{0}}\\
& &+C_{x_{\alpha}x_{\beta}w_{0}}+3[2]^{2}C_{x_{\alpha}^{2}w_{0}}+2([2]^{4}-[2]^{2})C_{w_{0}}\\
& &+([2]^{4}+3[2]^{2})C_{x_{\alpha}w_{0}}+2[2]^{2}C_{x_{\beta}w_{0}}\\
& &+[2]^{4}C_{x_{\alpha}w_{0}}+[2]^{2}C_{x_{\beta}w_{0}}+([2]^{4}-[2]^{2})C_{w_{0}}\\
& &-\{C_{x_{\alpha}x_{\beta}w_{0}}+3[2]^{2}C_{x_{\alpha}^{2}w_{0}}+2([2]^{4}+[2]^{2})C_{x_{\alpha}w_{0}}\\
& &+3[2]^{2}C_{x_{\beta}w_{0}}+(2[2]^{4}-[2]^{2})C_{w_{0}}\}\\
&=&[2]^{2}C_{x_{\alpha}x_{\beta}w_{0}}+2([2]^{4}-[2]^{2})C_{w_{0}}+(2[2]^{4}+[2]^{2})C_{x_{\alpha}w_{0}}\\
& &+[2]^{4}C_{x_{\beta}w_{0}}+([2]^{4}+[2]^{2})C_{x_{\alpha}^{2}w_{0}}.\\
\end{eqnarray*}

\begin{eqnarray*}
C_{1212120121}C_{10121210121212}&=&(C_{121212012}C_{1}-C_{12121201})C_{10121210121212}\\
&=&[2]C_{121212012}C_{10121210121212}-C_{12121201}C_{10121210121212}\\
&=&[2]\{C_{12121201}C_{210121210121212}+C_{12121201}C_{0121210121212}\\
& &-C_{1212120}C_{10121210121212}\}-C_{12121201}C_{10121210121212}\\
&=&[2]C_{12121201}C_{210121210121212}-C_{12121201}C_{10121210121212}\\
&=&[2]^{2}C_{x_{\alpha}x_{\beta}w_{0}}+(2[2]^{4}-[2]^{2})C_{w_{0}}+([2]^{6}+[2]^{2})C_{x_{\alpha}w_{0}}\\
& &+(2[2]^{4}-[2]^{2})C_{x_{\beta}w_{0}}+2[2]^{4}C_{x_{\alpha}^{2}w_{0}}\\
& &-\{([2]^{4}-[2]^{2})C_{x_{\alpha}^{2}w_{0}}+([2]^{6}-2[2]^{4}+[2]^{2})C_{x_{\alpha}w_{0}}\\
& &+([2]^{4}-[2]^{2})C_{x_{\beta}w_{0}}+([2]^{4}-[2]^{2})C_{w_{0}}\}\\
&=&[2]^{2}C_{x_{\alpha}x_{\beta}w_{0}}+[2]^{4}C_{w_{0}}+2[2]^{4}C_{x_{\alpha}w_{0}}+[2]^{4}C_{x_{\beta}w_{0}}\\
& &+([2]^{4}+[2]^{2})C_{x_{\alpha}^{2}w_{0}}.\\
\end{eqnarray*}

\begin{eqnarray*}
C_{12121201}C_{10121210121212}&=&(C_{121212}C_{0}C_{1}-C_{121212})C_{10121210121212}\\
&=&([2]^{2}-1)C_{121212}C_{10121210121212}\\
&=&([2]^{4}-[2]^{2})C_{x_{\alpha}^{2}w_{0}}+([2]^{6}-2[2]^{4}+[2]^{2})C_{x_{\alpha}w_{0}}\\
& &+([2]^{4}-[2]^{2})C_{x_{\beta}w_{0}}+([2]^{4}-[2]^{2})C_{w_{0}}.\\
\end{eqnarray*}

(82) Computing $C_{12121201210}C_{21201210121212}$
\begin{eqnarray*}
C_{12121201210}C_{21201210121212}&=&(C_{1212120121}C_{0}-C_{121212012})C_{21201210121212}\\
&=&C_{1212120121}C_{210121210121212}+C_{1212120121}C_{210121212}\\
& &+C_{1212120121}C_{0121212}-C_{121212012}C_{21201210121212}\\
&=&[2]C_{x_{\alpha}^{3}w_{0}}+([2]^{5}+4[2])C_{x_{\alpha}w_{0}}+(2[2]^{3}+[2])C_{x_{\beta}w_{0}}\\
& &+3[2]C_{x_{\alpha}x_{\beta}w_{0}}+(2[2]^{3}+3[2])C_{x_{\alpha}^{2}w_{0}}+2[2]^{3}C_{w_{0}}\\
& &+[2]C_{x_{\beta}w_{0}}+2([2]^{3}-[2])C_{x_{\alpha}w_{0}}+[2]^{3}C_{w_{0}}+[2]C_{x_{\alpha}^{2}w_{0}}\\
& &+2[2]C_{x_{\beta}w_{0}}+2[2]^{3}C_{x_{\alpha}w_{0}}+([2]^{5}-2[2]^{3}+[2])C_{w_{0}}\\
& &-\{[2]C_{x_{\alpha}x_{\beta}w_{0}}+([2]^{3}+[2])C_{x_{\alpha}^{2}w_{0}}+([2]^{5}+2[2])C_{x_{\alpha}w_{0}}\\
& &+([2]^{3}+[2])C_{x_{\beta}w_{0}}+([2]^{5}-[2]^{3}+[2])C_{w_{0}}\}\\
&=&[2]C_{x_{\alpha}^{3}w_{0}}+4[2]^{3}C_{x_{\alpha}w_{0}}+([2]^{3}+3[2])C_{x_{\beta}w_{0}}\\
& &+2[2]C_{x_{\alpha}x_{\beta}w_{0}}+([2]^{3}+3[2])C_{x_{\alpha}^{2}w_{0}}+2[2]^{3}C_{w_{0}}.\\
\end{eqnarray*}

(83) Computing $C_{12121201210}C_{121201210121212}$
\begin{eqnarray*}
C_{12121201210}C_{121201210121212}&=&C_{1212120121}C_{0}C_{121201210121212}-C_{121212012}C_{121201210121212}\\
&=&C_{1212120121}C_{0121201210121212}+C_{1212120121}C_{10121212}\\
& &-C_{121212012}C_{121201210121212}\\
&=&[2]^{2}C_{x_{\alpha}^{3}w_{0}}+2[2]^{4}C_{x_{\alpha}w_{0}}+4[2]^{2}C_{x_{\beta}w_{0}}\\
& &+2[2]^{2}C_{x_{\alpha}x_{\beta}w_{0}}+4[2]^{2}C_{x_{\alpha}^{2}w_{0}}+3[2]^{2}C_{w_{0}}\\
& &+[2]^{4}C_{x_{\alpha}w_{0}}+[2]^{2}C_{x_{\beta}w_{0}}+([2]^{4}-[2]^{2})C_{w_{0}}\\
& &-\{[2]^{2}C_{x_{\alpha}x_{\beta}w_{0}}+2[2]^{2}C_{x_{\beta}w_{0}}+([2]^{4}+2[2]^{2})C_{x_{\alpha}w_{0}}\\
& &+2[2]^{2}C_{x_{\alpha}^{2}w_{0}}+[2]^{4}C_{w_{0}}\}\\
&=&[2]^{2}C_{x_{\alpha}^{3}w_{0}}+2([2]^{4}-[2]^{2})C_{x_{\alpha}w_{0}}+3[2]^{2}C_{x_{\beta}w_{0}}+[2]^{2}C_{x_{\alpha}x_{\beta}w_{0}}\\
& &+2[2]^{2}C_{x_{\alpha}^{2}w_{0}}+2[2]^{2}C_{w_{0}}.\\
\end{eqnarray*}

(84) Computing $C_{12121201210}C_{0121201210121212}$
\begin{eqnarray*}
C_{12121201210}C_{0121201210121212}&=&C_{1212120121}C_{0}C_{0121201210121212}-C_{121212012}C_{0121201210121212}\\
&=&[2]C_{1212120121}C_{0121201210121212}-C_{121212012}C_{0121201210121212}\\
&=&[2]^{3}C_{x_{\alpha}^{3}w_{0}}+2[2]^{5}C_{x_{\alpha}w_{0}}+4[2]^{3}C_{x_{\beta}w_{0}}\\
& &+2[2]^{3}C_{x_{\alpha}x_{\beta}w_{0}}+4[2]^{3}C_{x_{\alpha}^{2}w_{0}}+3[2]^{3}C_{w_{0}}\\
& &-\{[2]^{3}C_{x_{\alpha}x_{\beta}w_{0}}+(2[2]^{3}-[2])C_{x_{\beta}w_{0}}+([2]^{5}+[2])C_{x_{\alpha}w_{0}}\\
& &+2[2]^{3}C_{x_{\alpha}^{2}w_{0}}+(2[2]^{3}-[2])C_{w_{0}}\}\\
&=&[2]^{3}C_{x_{\alpha}^{3}w_{0}}+[2]^{3}C_{x_{\alpha}x_{\beta}w_{0}}+(2[2]^{3}+[2])C_{x_{\beta}w_{0}}+([2]^{5}-[2])C_{x_{\alpha}w_{0}}\\
& &+2[2]^{3}C_{x_{\alpha}^{2}w_{0}}+([2]^{3}+[2])C_{w_{0}}.\\
\end{eqnarray*}

(85) Computing $C_{121212012102}C_{121212}$
\begin{eqnarray*}
C_{121212012102}C_{121212}&=&C_{121212}C_{201210121212}=[2]^{2}C_{x_{\beta}w_{0}}+([2]^{4}-[2]^{2})C_{x_{\alpha}121212}+[2]^{2}C_{w_{0}}.
\end{eqnarray*}

(86) Computing $C_{121212012102}C_{0121212}$
\begin{eqnarray*}
C_{121212012102}C_{0121212}&=&C_{1212120}C_{201210121212}=[2]^{3}C_{x_{\beta}w_{0}}+([2]^{5}-[2]^{3})C_{x_{\alpha}w_{0}}+[2]^{3}C_{w_{0}}.
\end{eqnarray*}

(87) Computing $C_{121212012102}C_{10121212}$
\begin{eqnarray*}
C_{121212012102}C_{10121212}&=&C_{12121201}C_{201210121212}=[2]^{2}C_{x_{\alpha}^{2}w_{0}}+(2[2]^{4}-[2]^{2})C_{x_{\alpha}w_{0}}\\
& &+2[2]^{2}C_{x_{\beta}w_{0}}+[2]^{4}C_{w_{0}}.
\end{eqnarray*}

(88) Computing $C_{121212012102}C_{210121212}$
\begin{eqnarray*}
C_{121212012102}C_{210121212}&=&C_{121212012}C_{201210121212}=[2]^{3}C_{x_{\alpha}^{2}w_{0}}+[2]^{5}C_{x_{\alpha}w_{0}}\\
& &+[2]^{3}C_{x_{\beta}w_{0}}+([2]^{5}-[2]^{3})C_{w_{0}}.
\end{eqnarray*}

(89) Computing $C_{121212012102}C_{1210121212}$
\begin{eqnarray*}
C_{121212012102}C_{1210121212}&=&C_{1212120121}C_{201210121212}=C_{x_{\alpha}x_{\beta}w_{0}}+3[2]^{2}C_{x_{\alpha}^{2}w_{0}}\\
& &+2([2]^{4}-[2]^{2})C_{w_{0}}+([2]^{4}+3[2]^{2})C_{x_{\alpha}w_{0}}+2[2]^{2}C_{x_{\beta}w_{0}}.
\end{eqnarray*}

(90) Computing $C_{121212012102}C_{21210121212}$
\begin{eqnarray*}
C_{121212012102}C_{21210121212}&=&C_{12121201212}C_{201210121212}=[2]C_{x_{\alpha}x_{\beta}w_{0}}+2[2]^{3}C_{x_{\alpha}^{2}w_{0}}+3[2]^{3}C_{x_{\alpha}w_{0}}\\\
& &+[2]^{3}C_{x_{\beta}w_{0}}+([2]^{5}-[2]^{3})C_{w_{0}}.
\end{eqnarray*}

(91) Computing $C_{121212012102}C_{01210121212}$
\begin{eqnarray*}
C_{121212012102}C_{01210121212}&=&C_{12121201210}C_{201210121212}=[2]C_{x_{\alpha}x_{\beta}w_{0}}+2[2]^{3}C_{x_{\alpha}^{2}w_{0}}\\
& &+([2]^{5}-[2]^{3})C_{w_{0}}+3[2]^{3}C_{x_{\alpha}w_{0}}+[2]^{3}C_{x_{\beta}w_{0}}.
\end{eqnarray*}

(92) Computing $C_{121212012102}C_{201210121212}$
\begin{eqnarray*}
C_{12121201210}C_{2}&=&C_{121212012102}.\\
C_{121212012102}C_{201210121212}&=&C_{12121201210}C_{2}C_{201210121212}\\
&=&[2]C_{12121201210}C_{201210121212}.\\
\end{eqnarray*}

(93) Computing $C_{121212012102}C_{1201210121212}$
\begin{eqnarray*}
C_{2}C_{1201210121212}&=&C_{21201210121212}+C_{121210121212}+C_{201210121212},\\
C_{121212012102}C_{1201210121212}&=&C_{12121201210}C_{2}C_{1201210121212}\\
&=&C_{12121201210}C_{2}C_{1201210121212}\\
&=&C_{12121201210}C_{21201210121212}+C_{12121201210}C_{121210121212}\\
& &+C_{12121201210}C_{201210121212}\\
&=&[2]C_{x_{\alpha}^{3}w_{0}}+4[2]^{3}C_{x_{\alpha}w_{0}}+([2]^{3}+3[2])C_{x_{\beta}w_{0}}\\
& &+2[2]C_{x_{\alpha}x_{\beta}w_{0}}+([2]^{3}+3[2])C_{x_{\alpha}^{2}w_{0}}+2[2]^{3}C_{w_{0}}\\
& &[2]C_{x_{\alpha}x_{\beta}w_{0}}+[2]^{3}C_{x_{\alpha}^{2}w_{0}}+([2]^{3}+[2])C_{x_{\alpha}w_{0}}\\
& &+([2]^{3}-[2])C_{x_{\beta}w_{0}}+([2]^{3}-[2])C_{w_{0}}\\
& &+[2]C_{x_{\alpha}x_{\beta}w_{0}}+2[2]^{3}C_{x_{\alpha}^{2}w_{0}}+([2]^{5}-[2]^{3})C_{w_{0}}\\
& &+3[2]^{3}C_{x_{\alpha}w_{0}}+[2]^{3}C_{x_{\beta}w_{0}}\\
&=&[2]C_{x_{\alpha}^{3}w_{0}}+(8[2]^{3}+[2])C_{x_{\alpha}w_{0}}+(3[2]^{3}+2[2])C_{x_{\beta}w_{0}}\\
& &+4[2]C_{x_{\alpha}x_{\beta}w_{0}}+(4[2]^{3}+3[2])C_{x_{\alpha}^{2}w_{0}}+([2]^{5}+2[2]^{3}-[2])C_{w_{0}}.\\
\end{eqnarray*}

\begin{eqnarray*}
C_{12121201210}C_{121210121212}&=&C_{12121201210}(C_{1}C_{21210121212}-C_{1210121212})\\
&=&C_{12121201210}C_{1}C_{21210121212}-C_{12121201210}C_{1210121212}\\
&=&[2]C_{12121201210}C_{21210121212}-C_{12121201210}C_{1210121212}\\
&=&[2]C_{x_{\alpha}x_{\beta}w_{0}}+2[2]^{3}C_{x_{\alpha}^{2}w_{0}}+3[2]^{3}C_{x_{\alpha}w_{0}}\\\
& &+[2]^{3}C_{x_{\beta}w_{0}}+([2]^{5}-[2]^{3})C_{w_{0}}\\
& &-\{[2]^{3}C_{x_{\alpha}^{2}w_{0}}+(2[2]^{3}-[2])C_{x_{\alpha}w_{0}}\\
& &+[2]C_{x_{\beta}w_{0}}+([2]^{5}-2[2]^{3}+[2])C_{w_{0}}\}\\
&=&[2]C_{x_{\alpha}x_{\beta}w_{0}}+[2]^{3}C_{x_{\alpha}^{2}w_{0}}+([2]^{3}+[2])C_{x_{\alpha}w_{0}}\\\
& &+([2]^{3}-[2])C_{x_{\beta}w_{0}}+([2]^{3}-[2])C_{w_{0}}.\\
\end{eqnarray*}
(94) Computing $C_{121212012102}C_{21201210121212}$
\begin{eqnarray*}
C_{121212012102}C_{21201210121212}&=&C_{12121201210}C_{2}C_{21201210121212}\\
&=&C_{12121201210}C_{2}C_{21201210121212}\\
&=&[2]C_{12121201210}C_{21201210121212}\\
&=&[2]^{2}C_{x_{\alpha}^{3}w_{0}}+4[2]^{4}C_{x_{\alpha}w_{0}}+([2]^{4}+3[2]^{2})C_{x_{\beta}w_{0}}\\
& &+2[2]^{2}C_{x_{\alpha}x_{\beta}w_{0}}+([2]^{4}+3[2]^{2})C_{x_{\alpha}^{2}w_{0}}+2[2]^{4}C_{w_{0}}.\\
\end{eqnarray*}

(95) Computing $C_{121212012102}C_{121201210121212}$
\begin{eqnarray*}
C_{2}C_{121201210121212}&=&C_{x_{\beta}w_{0}}+C_{21201210121212}+C_{w_{0}}.\\
C_{121212012102}C_{121201210121212}&=&C_{12121201210}C_{2}C_{121201210121212}\\
&=&C_{12121201210}C_{x_{\beta}121212}+C_{12121201210}C_{21201210121212}\\
& &+C_{12121201210}C_{w_{0}}\\
&=&(C_{1212120121}C_{0}-C_{121212012})C_{x_{\beta}121212}\\
& &+C_{12121201210}C_{21201210121212}+C_{12121201210}C_{w_{0}}\\
&=&C_{1212120121}C_{02121021210121212}-C_{121212012}C_{2121021210121212}\\
& &+C_{12121201210}C_{21201210121212}+C_{12121201210}C_{w_{0}}\\
&=&[2]C_{x_{\beta}^{2}w_{0}}+2([2]^{3}-[2])C_{x_{\alpha}x_{\beta}w_{0}}+(2[2]^{3}-[2])C_{x_{\alpha}^{2}w_{0}}\\
& &+2([2]^{3}-[2])C_{x_{\alpha}w_{0}}+([2]^{3}+[2])C_{x_{\beta}w_{0}}+[2]C_{w_{0}}+[2]C_{x_{\alpha}^{3}w_{0}}\\
& &-\{([2]^{3}-[2])C_{x_{\alpha}x_{\beta}w_{0}}+([2]^{3}-[2])C_{x_{\alpha}^{2}w_{0}}+([2]^{3}-[2])C_{x_{\alpha}w_{0}}\\
& &+([2]^{3}-[2])C_{x_{\beta}w_{0}}\}\\
& &+[2]C_{x_{\alpha}^{3}w_{0}}+4[2]^{3}C_{x_{\alpha}w_{0}}+([2]^{3}+3[2])C_{x_{\beta}w_{0}}\\
& &+2[2]C_{x_{\alpha}x_{\beta}w_{0}}+([2]^{3}+3[2])C_{x_{\alpha}^{2}w_{0}}+2[2]^{3}C_{w_{0}}\\
& &+[2]C_{x_{\beta}121212}+([2]^{3}-[2])C_{x_{\alpha}121212}+[2]C_{w_{0}}\\
&=&[2]C_{x_{\beta}^{2}w_{0}}+([2]^{3}+[2])C_{x_{\alpha}x_{\beta}w_{0}}+(2[2]^{3}+3[2])C_{x_{\alpha}^{2}w_{0}}+2[2]C_{x_{\alpha}^{3}w_{0}}\\
& &+(6[2]^{3}-2[2])C_{x_{\alpha}w_{0}}+([2]^{3}+6[2])C_{x_{\beta}w_{0}}+(2[2]^{3}+2[2])C_{w_{0}}\\
\end{eqnarray*}

(96) Computing $C_{121212012102}C_{0121201210121212}$
\begin{eqnarray*}
C_{12121201210}C_{2}&=&C_{121212012102},\\
C_{2}C_{0121201210121212}&=&C_{20121201210121212}+C_{021201210121212}+C_{0121212}\\
C_{121212012102}C_{0121201210121212}&=&C_{12121201210}C_{2}C_{0121201210121212}\\
&=&C_{12121201210}C_{20121201210121212}+C_{12121201210}C_{210121210121212}\\
& &+C_{12121201210}C_{0121212}\\
&=&[2]^{2}C_{x_{\beta}^{2}w_{0}}+([2]^{4}-[2]^{2})C_{x_{\alpha}x_{\beta}w_{0}}+[2]^{4}C_{x_{\alpha}^{2}w_{0}}+([2]^{4}-[2]^{2})C_{x_{\alpha}w_{0}}\\
& &+2[2]^{2}C_{x_{\beta}w_{0}}+[2]^{2}C_{w_{0}}+[2]^{2}C_{x_{\alpha}^{3}w_{0}}\\
& &+[2]^{2}C_{x_{\alpha}^{3}w_{0}}+(2[2]^{4}+[2]^{2})C_{x_{\alpha}w_{0}}+([2]^{4}+[2]^{2})C_{x_{\beta}w_{0}}\\
& &+2[2]^{2}C_{x_{\alpha}x_{\beta}w_{0}}+([2]^{4}+2[2]^{2})C_{x_{\alpha}^{2}w_{0}}+[2]^{4}C_{w_{0}}\\
& &+[2]^{2}C_{x_{\beta}w_{0}}+([2]^{4}-[2]^{2})C_{x_{\alpha}w_{0}}+[2]^{2}C_{w_{0}}\\
&=&[2]^{2}C_{x_{\beta}^{2}w_{0}}+([2]^{4}+[2]^{2})C_{x_{\alpha}x_{\beta}w_{0}}+(2[2]^{4}+2[2]^{2})C_{x_{\alpha}^{2}w_{0}}\\
& &+(4[2]^{4}-[2]^{2})C_{x_{\alpha}w_{0}}+([2]^{4}+4[2]^{2})C_{x_{\beta}w_{0}}+([2]^{4}+[2]^{2})C_{w_{0}}\\
& &+2[2]^{2}C_{x_{\alpha}^{3}w_{0}}.\\
\end{eqnarray*}

\begin{eqnarray*}
C_{12121201210}C_{20121201210121212}&=&C_{12121201210}C_{0}C_{2121201210121212}\\
&=&[2]C_{12121201210}C_{2121201210121212}\\
&=&[2]^{2}C_{x_{\beta}^{2}w_{0}}+([2]^{4}-[2]^{2})C_{x_{\alpha}x_{\beta}w_{0}}+[2]^{4}C_{x_{\alpha}^{2}w_{0}}+([2]^{4}-[2]^{2})C_{x_{\alpha}w_{0}}\\
& &+2[2]^{2}C_{x_{\beta}w_{0}}+[2]^{2}C_{w_{0}}+[2]^{2}C_{x_{\alpha}^{3}w_{0}}\\
\end{eqnarray*}

\begin{eqnarray*}
C_{12121201210}C_{2121201210121212}&=&(C_{1212120121}C_{0}-C_{121212012})C_{2121201210121212}\\
&=&C_{1212120121}C_{02121201210121212}-C_{121212012}C_{2121201210121212}\\
&=&[2]C_{x_{\beta}^{2}w_{0}}+2([2]^{3}-[2])C_{x_{\alpha}x_{\beta}w_{0}}+(2[2]^{3}-[2])C_{x_{\alpha}^{2}w_{0}}\\
& &+2([2]^{3}-[2])C_{x_{\alpha}w_{0}}+([2]^{3}+[2])C_{x_{\beta}w_{0}}+[2]C_{w_{0}}+[2]C_{x_{\alpha}^{3}w_{0}}\\
& &-\{([2]^{3}-[2])C_{x_{\alpha}x_{\beta}w_{0}}+([2]^{3}-[2])C_{x_{\alpha}^{2}w_{0}}+([2]^{3}-[2])C_{x_{\alpha}w_{0}}\\
& &+([2]^{3}-[2])C_{x_{\beta}w_{0}}\}\\
&=&[2]C_{x_{\beta}^{2}w_{0}}+([2]^{3}-[2])C_{x_{\alpha}x_{\beta}w_{0}}+[2]^{3}C_{x_{\alpha}^{2}w_{0}}+([2]^{3}-[2])C_{x_{\alpha}w_{0}}\\
& &+2[2]C_{x_{\beta}w_{0}}+[2]C_{w_{0}}+[2]C_{x_{\alpha}^{3}w_{0}}\\
\end{eqnarray*}

\begin{eqnarray*}
C_{12121201210}C_{210121210121212}&=&(C_{1212120121}C_{0}-C_{121212012})C_{210121210121212}\\
&=&[2]C_{1212120121}C_{210121210121212}-C_{121212012}C_{210121210121212}\\
&=&[2]^{2}C_{x_{\alpha}^{3}w_{0}}+([2]^{6}+4[2]^{2})C_{x_{\alpha}w_{0}}+(2[2]^{4}+[2]^{2})C_{x_{\beta}w_{0}}\\
& &+3[2]^{2}C_{x_{\alpha}x_{\beta}w_{0}}+(2[2]^{4}+3[2]^{2})C_{x_{\alpha}^{2}w_{0}}+2[2]^{4}C_{w_{0}}\\
& &-\{[2]^{2}C_{x_{\alpha}x_{\beta}w_{0}}+[2]^{4}C_{w_{0}}+([2]^{6}-2[2]^{4}+3[2]^{2})C_{x_{\alpha}w_{0}}\\
& &+[2]^{4}C_{x_{\beta}w_{0}}+([2]^{4}+[2]^{2})C_{x_{\alpha}^{2}w_{0}}\}\\
&=&[2]^{2}C_{x_{\alpha}^{3}w_{0}}+(2[2]^{4}+[2]^{2})C_{x_{\alpha}w_{0}}+([2]^{4}+[2]^{2})C_{x_{\beta}w_{0}}\\
& &+2[2]^{2}C_{x_{\alpha}x_{\beta}w_{0}}+([2]^{4}+2[2]^{2})C_{x_{\alpha}^{2}w_{0}}+[2]^{4}C_{w_{0}}.\\
\end{eqnarray*}

\begin{eqnarray*}
C_{121212012}C_{210121210121212}&=&(C_{12121201}C_{2}-C_{1212120})C_{210121210121212}\\
&=&[2](C_{12121201}C_{210121210121212}-C_{121212}C_{210121210121212})\\
&=&[2]^{2}C_{x_{\alpha}x_{\beta}w_{0}}+(2[2]^{4}-[2]^{2})C_{w_{0}}+([2]^{6}+[2]^{2})C_{x_{\alpha}w_{0}}\\
& &+(2[2]^{4}-[2]^{2})C_{x_{\beta}w_{0}}+2[2]^{4}C_{x_{\alpha}^{2}w_{0}}\\
& &-\{([2]^{4}-[2]^{2})C_{x_{\alpha}^{2}w_{0}}+2([2]^{4}-[2]^{2})C_{x_{\alpha}w_{0}}\\
& &+([2]^{4}-[2]^{2})C_{x_{\beta}w_{0}}+([2]^{4}-[2]^{2})C_{w_{0}}\}\\
&=&[2]^{2}C_{x_{\alpha}x_{\beta}w_{0}}+[2]^{4}C_{w_{0}}+([2]^{6}-2[2]^{4}+3[2]^{2})C_{x_{\alpha}w_{0}}\\
& &+[2]^{4}C_{x_{\beta}w_{0}}+([2]^{4}+[2]^{2})C_{x_{\alpha}^{2}w_{0}}.\\
\end{eqnarray*}

(97) Computing $C_{1212120121021}C_{121212}$
\begin{eqnarray*}
C_{1212120121021}C_{121212}&=&C_{121212}C_{1201210121212}=([2]^{3}-[2])C_{x_{\beta}w_{0}}+([2]^{3}-[2])C_{w_{0}}\\
& &+([2]^{5}-2[2]^{3}+[2])C_{x_{\alpha}w_{0}}.
\end{eqnarray*}

(98) Computing $C_{1212120121021}C_{0121212}$
\begin{eqnarray*}
C_{1212120121021}C_{0121212}&=&C_{1212120}C_{1201210121212}=[2]^{2}C_{x_{\alpha}^{2}121212}+(2[2]^{4}-[2]^{2})C_{x_{\alpha}121212}\\
& &+2[2]^{2}C_{x_{\beta}121212}+[2]^{4}C_{w_{0}}.
\end{eqnarray*}

(99) Computing $C_{1212120121021}C_{10121212}$
\begin{eqnarray*}
C_{1212120121021}C_{10121212}&=&C_{12121201}C_{1201210121212}=[2]^{3}C_{x_{\alpha}^{2}w_{0}}+([2]^{5}+[2]^{3}-[2])C_{x_{\alpha}w_{0}}\\
& &+([2]^{3}+[2])C_{x_{\beta}w_{0}}+([2]^{5}-[2]^{3}+[2])C_{w_{0}}.
\end{eqnarray*}

(100) Computing $C_{1212120121021}C_{210121212}$
\begin{eqnarray*}
C_{1212120121021}C_{210121212}&=&C_{121212012}C_{2101210121212}=C_{x_{\alpha}x_{\beta}w_{0}}+3[2]^{2}C_{x_{\alpha}^{2}w_{0}}\\
& &+2([2]^{4}+[2]^{2})C_{x_{\alpha}w_{0}}+3[2]^{2}C_{x_{\beta}w_{0}}+(2[2]^{4}-[2]^{2})C_{w_{0}}.
\end{eqnarray*}

(101) Computing $C_{1212120121021}C_{1210121212}$
\begin{eqnarray*}
C_{1212120121021}C_{1210121212}&=&C_{1212120121}C_{1201210121212}=[2]C_{x_{\alpha}x_{\beta}w_{0}}+2[2]^{3}C_{x_{\alpha}^{2}w_{0}}\\
& &+([2]^{5}-[2])C_{w_{0}}+([2]^{5}+[2]^{3}+[2])C_{x_{\alpha}w_{0}}+(2[2]^{3}-[2])C_{x_{\beta}w_{0}}.
\end{eqnarray*}

(102) Computing $C_{1212120121021}C_{21210121212}$
\begin{eqnarray*}
C_{1212120121021}C_{21210121212}&=&C_{12121201212}C_{1201210121212}=[2]^{2}C_{x_{\alpha}x_{\beta}w_{0}}+([2]^{4}+[2]^{2})C_{x_{\alpha}^{2}w_{0}}\\
& &+(2[2]^{4}-2[2]^{2})C_{w_{0}}+(2[2]^{4}+[2]^{2})C_{x_{\alpha}w_{0}}+[2]^{4}C_{x_{\beta}w_{0}}.
\end{eqnarray*}

(103) Computing $C_{1212120121021}C_{01210121212}$
\begin{eqnarray*}
C_{1212120121021}C_{01210121212}&=&C_{12121201210}C_{1201210121212}=[2]^{2}C_{x_{\alpha}x_{\beta}w_{0}}+2([2]^{4}-[2]^{2})C_{w_{0}}\\
& &+(2[2]^{4}+[2]^{2})C_{x_{\alpha}w_{0}}+[2]^{4}C_{x_{\beta}w_{0}}+([2]^{4}+[2]^{2})C_{x_{\alpha}^{2}w_{0}}.
\end{eqnarray*}

(104) Computing $C_{1212120121021}C_{201210121212}$
\begin{eqnarray*}
C_{1212120121021}C_{201210121212}&=&C_{121212012102}C_{1201210121212}=[2]C_{x_{\alpha}^{3}w_{0}}+(8[2]^{3}+[2])C_{x_{\alpha}w_{0}}\\
& &+(3[2]^{3}+2[2])C_{x_{\beta}w_{0}}+4[2]C_{x_{\alpha}x_{\beta}w_{0}}+(4[2]^{3}+3[2])C_{x_{\alpha}^{2}w_{0}}\\
& &+([2]^{5}+2[2]^{3}-[2])C_{w_{0}}.
\end{eqnarray*}

(105) Computing $C_{1212120121021}C_{1201210121212}$
\begin{eqnarray*}
C_{1212120121021}C_{1201210121212}&=&(C_{121212012102}C_{1}-C_{12121201210})C_{1201210121212}\\
&=&C_{121212012102}C_{1}C_{1201210121212}-C_{12121201210}C_{1201210121212}\\
&=&[2]C_{121212012102}C_{1201210121212}-C_{12121201210}C_{1201210121212}\\
&=&[2]^{2}C_{x_{\alpha}^{3}w_{0}}+(8[2]^{4}+[2]^{2})C_{x_{\alpha}w_{0}}+(3[2]^{4}+2[2]^{2})C_{x_{\beta}w_{0}}\\
& &+4[2]^{2}C_{x_{\alpha}x_{\beta}w_{0}}+(4[2]^{4}+3[2]^{2})C_{x_{\alpha}^{2}w_{0}}+([2]^{6}+2[2]^{4}-[2]^{2})C_{w_{0}}\\
& &-\{[2]^{2}C_{x_{\alpha}x_{\beta}w_{0}}+2([2]^{4}-[2]^{2})C_{w_{0}}+(2[2]^{4}+[2]^{2})C_{x_{\alpha}w_{0}}\\
& &+[2]^{4}C_{x_{\beta}w_{0}}+([2]^{4}+[2]^{2})C_{x_{\alpha}^{2}w_{0}}\}\\
&=&[2]^{2}C_{x_{\alpha}^{3}w_{0}}+6[2]^{4}C_{x_{\alpha}w_{0}}+(2[2]^{4}+2[2]^{2})C_{x_{\beta}w_{0}}\\
& &+3[2]^{2}C_{x_{\alpha}x_{\beta}w_{0}}+(3[2]^{4}+2[2]^{2})C_{x_{\alpha}^{2}w_{0}}+([2]^{6}+[2]^{2})C_{w_{0}}.\\
\end{eqnarray*}

(106) Computing $C_{1212120121021}C_{21201210121212}$
\begin{eqnarray*}
C_{1212120121021}C_{21201210121212}&=&(C_{121212012102}C_{1}-C_{12121201210})C_{21201210121212}\\
&=&C_{121212012102}C_{1}C_{21201210121212}-C_{12121201210}C_{21201210121212}\\
&=&C_{121212012102}C_{121201210121212}+C_{121212012102}C_{1201210121212}\\
& &-C_{12121201210}C_{21201210121212}\\
&=&[2]C_{x_{\beta}^{2}w_{0}}+([2]^{3}+[2])C_{x_{\alpha}x_{\beta}w_{0}}+(2[2]^{3}+3[2])C_{x_{\alpha}^{2}w_{0}}+2[2]C_{x_{\alpha}^{3}w_{0}}\\
& &+(6[2]^{3}-2[2])C_{x_{\alpha}w_{0}}+([2]^{3}+6[2])C_{x_{\beta}w_{0}}+(2[2]^{3}+2[2])C_{w_{0}}\\
& &[2]C_{x_{\alpha}^{3}w_{0}}+(8[2]^{3}+[2])C_{x_{\alpha}w_{0}}+(3[2]^{3}+2[2])C_{x_{\beta}w_{0}}\\
& &+4[2]C_{x_{\alpha}x_{\beta}w_{0}}+(4[2]^{3}+3[2])C_{x_{\alpha}^{2}w_{0}}+([2]^{5}+2[2]^{3}-[2])C_{w_{0}}\\
& &-\{[2]C_{x_{\alpha}^{3}w_{0}}+4[2]^{3}C_{x_{\alpha}w_{0}}+([2]^{3}+3[2])C_{x_{\beta}w_{0}}\\
& &+2[2]C_{x_{\alpha}x_{\beta}w_{0}}+([2]^{3}+3[2])C_{x_{\alpha}^{2}w_{0}}+2[2]^{3}C_{w_{0}}\}\\
&=&[2]C_{x_{\beta}^{2}w_{0}}+([2]^{3}+3[2])C_{x_{\alpha}x_{\beta}w_{0}}+(5[2]^{3}+3[2])C_{x_{\alpha}^{2}w_{0}}\\
& &+(10[2]^{3}-[2])C_{x_{\alpha}w_{0}}+(3[2]^{3}+5[2])C_{x_{\beta}w_{0}}\\
& &+([2]^{5}+2[2]^{3}+[2])C_{w_{0}}+2[2]C_{x_{\alpha}^{3}w_{0}}.\\
\end{eqnarray*}

(107) Computing $C_{1212120121021}C_{121201210121212}$
\begin{eqnarray*}
C_{1212120121021}C_{121201210121212}&=&(C_{121212012102}C_{1}-C_{12121201210})C_{121201210121212}\\
&=&C_{121212012102}C_{1}C_{121201210121212}-C_{12121201210}C_{121201210121212}\\
&=&[2]C_{121212012102}C_{121201210121212}-C_{12121201210}C_{121201210121212}\\
&=&[2]^{2}C_{x_{\beta}^{2}w_{0}}+([2]^{4}+[2]^{2})C_{x_{\alpha}x_{\beta}w_{0}}+(2[2]^{4}+3[2]^{2})C_{x_{\alpha}^{2}w_{0}}+2[2]^{2}C_{x_{\alpha}^{3}w_{0}}\\
& &+(6[2]^{4}-2[2]^{2})C_{x_{\alpha}w_{0}}+([2]^{4}+6[2]^{2})C_{x_{\beta}w_{0}}+(2[2]^{4}+2[2]^{2})C_{w_{0}}\\
& &-\{[2]^{2}C_{x_{\alpha}^{3}w_{0}}+2([2]^{4}-[2]^{2})C_{x_{\alpha}w_{0}}+3[2]^{2}C_{x_{\beta}w_{0}}\\
& &+[2]^{2}C_{x_{\alpha}x_{\beta}w_{0}}+2[2]^{2}C_{x_{\alpha}^{2}w_{0}}+2[2]^{2}C_{w_{0}}\}\\
&=&[2]^{2}C_{x_{\beta}^{2}w_{0}}+[2]^{4}C_{x_{\alpha}x_{\beta}w_{0}}+(2[2]^{4}+[2]^{2})C_{x_{\alpha}^{2}w_{0}}\\
& &+4[2]^{4}C_{x_{\alpha}w_{0}}+([2]^{4}+3[2]^{2})C_{x_{\beta}w_{0}}+2[2]^{4}C_{w_{0}}+[2]^{2}C_{x_{\alpha}^{3}w_{0}}.\\
\end{eqnarray*}

(108) Computing $C_{1212120121021}C_{0121201210121212}$
\begin{eqnarray*}
C_{121212012102}C_{1}&=&C_{1212120121021}+C_{12121201210},\\
C_{1212120121021}&=&C_{121212012102}C_{1}-C_{12121201210},\\
C_{1}C_{0121201210121212}&=&C_{10121201210121212}+C_{121201210121212}\\
&=&C_{01210121210121212}+C_{121201210121212}.\\
C_{1212120121021}C_{0121201210121212}&=&(C_{121212012102}C_{1}-C_{12121201210})C_{0121201210121212}\\
&=&C_{121212012102}C_{1}C_{0121201210121212}-C_{12121201210}C_{0121201210121212}\\
&=&C_{121212012102}C_{01210121210121212}+C_{121212012102}C_{121201210121212}\\
& &-C_{12121201210}C_{0121201210121212}\\
&=&[2]C_{x_{\alpha}^{2}x_{\beta}w_{0}}+[2]C_{x_{\beta}^{2}w_{0}}+3[2]^{3}C_{w_{0}}+6[2]^{3}C_{x_{\beta}w_{0}}\\
& &+(5[2]^{3}+2[2])C_{x_{\alpha}^{2}w_{0}}+(3[2]^{3}+[2])C_{x_{\alpha}x_{\beta}w_{0}}\\
& &+([2]^{5}+3[2]^{3}+[2])C_{x_{\alpha}w_{0}}+(2[2]^{3}+[2])C_{x_{\alpha}^{3}w_{0}}+[2]C_{x_{\beta}^{2}w_{0}}\\
& &+([2]^{3}+[2])C_{x_{\alpha}x_{\beta}w_{0}}+(2[2]^{3}+3[2])C_{x_{\alpha}^{2}w_{0}}+2[2]C_{x_{\alpha}^{3}w_{0}}\\
& &+(6[2]^{3}-2[2])C_{x_{\alpha}w_{0}}+([2]^{3}+6[2])C_{x_{\beta}w_{0}}+(2[2]^{3}+2[2])C_{w_{0}}\\
& &-\{[2]^{3}C_{x_{\alpha}^{3}w_{0}}+[2]^{3}C_{x_{\alpha}x_{\beta}w_{0}}+(2[2]^{3}+[2])C_{x_{\beta}w_{0}}+([2]^{5}-[2])C_{x_{\alpha}w_{0}}\\
& &+2[2]^{3}C_{x_{\alpha}^{2}w_{0}}+([2]^{3}+[2])C_{w_{0}}\}\\
&=&[2]C_{x_{\alpha}^{2}x_{\beta}w_{0}}+2[2]C_{x_{\beta}^{2}w_{0}}+(4[2]^{3}+[2])C_{w_{0}}+(5[2]^{3}+5[2])C_{x_{\beta}w_{0}}\\
& &+(5[2]^{3}+5[2])C_{x_{\alpha}^{2}w_{0}}+(3[2]^{3}+2[2])C_{x_{\alpha}x_{\beta}w_{0}}\\
& &+9[2]^{3}C_{x_{\alpha}w_{0}}+([2]^{3}+3[2])C_{x_{\alpha}^{3}w_{0}}\\
\end{eqnarray*}

\begin{eqnarray*}
C_{121212012102}C_{01210121210121212}&=&C_{12121201210}C_{2}C_{01210121210121212}\\
&=&C_{12121201210}C_{021210121210121212}\\
&=&(C_{1212120121}C_{0}-C_{121212012})C_{021210121210121212}\\
&=&[2]C_{1212120121}C_{021210121210121212}-[2]C_{121212012}C_{01210121210121212}\\
&=&[2](C_{121212012}C_{1}-C_{12121201})C_{021210121210121212}\\
& &-[2]C_{121212012}C_{01210121210121212}\\
&=&[2]\{C_{121212012}C_{1021210121210121212}+C_{121212012}C_{01210121210121212}\\
& &-C_{12121201}C_{021210121210121212}\}-[2]C_{121212012}C_{01210121210121212}\\
&=&[2]\{C_{121212012}C_{1021210121210121212}-C_{12121201}C_{021210121210121212}\}\\
&=&[2]\{(C_{12121201}C_{2}-C_{1212120})C_{1021210121210121212}\\
& &-C_{12121201}C_{021210121210121212}\}\\
&=&[2]\{C_{12121201}C_{2}C_{1021210121210121212}-C_{1212120}C_{1021210121210121212}\\
& &-C_{12121201}C_{021210121210121212}\}\\
&=&[2]\{C_{12121201}(C_{21021210121210121212}+C_{w_{0}}+C_{x_{\alpha}w_{0}}\\
& &+C_{x_{\alpha}^{2}w_{0}}+C_{x_{\beta}w_{0}}+C_{021210121210121212})\\
& &-C_{1212120}C_{1021210121210121212}-C_{12121201}C_{021210121210121212}\}\\
&=&[2]\{C_{12121201}(C_{21021210121210121212}+C_{w_{0}}+C_{x_{\alpha}w_{0}}\\
& &+C_{x_{\alpha}^{2}w_{0}}+C_{x_{\beta}w_{0}})-C_{1212120}C_{1021210121210121212}\}\\
&=&[2]\{(C_{1212120}C_{1}-C_{121212})C_{21021210121210121212}\\
& &+C_{12121201}(C_{w_{0}}+C_{x_{\alpha}w_{0}}+C_{x_{\alpha}^{2}w_{0}}+C_{x_{\beta}w_{0}})\\
& &-C_{1212120}C_{1021210121210121212}\}\\
&=&[2]\{C_{1212120}C_{121021210121210121212}+C_{1212120}C_{1021210121210121212}\\
& &-C_{121212}C_{21021210121210121212}+C_{12121201}(C_{w_{0}}+C_{x_{\alpha}w_{0}}\\
& &+C_{x_{\alpha}^{2}w_{0}}+C_{x_{\beta}w_{0}})-C_{1212120}C_{1021210121210121212}\}\\
&=&[2]\{C_{1212120}C_{121021210121210121212}-C_{121212}C_{21021210121210121212}\\
& &+C_{12121201}(C_{w_{0}}+C_{x_{\alpha}w_{0}}+C_{x_{\alpha}^{2}w_{0}}+C_{x_{\beta}w_{0}})\}\\
&=&[2]\{C_{121212}C_{0}C_{121021210121210121212}-C_{121212}C_{21021210121210121212}\\
& &+C_{12121201}(C_{w_{0}}+C_{x_{\alpha}w_{0}}+C_{x_{\alpha}^{2}w_{0}}+C_{x_{\beta}w_{0}})\}\\
&=&[2]\{C_{121212}C_{0121021210121210121212}+C_{121212}C_{10121210121212}\\
& &-C_{121212}C_{21021210121210121212}+C_{12121201}(C_{w_{0}}+C_{x_{\alpha}w_{0}}\\
& &+C_{x_{\alpha}^{2}w_{0}}+C_{x_{\beta}w_{0}})\}\\
&=&[2]\{C_{x_{\alpha}^{2}x_{\beta}w_{0}}+C_{x_{\beta}^{2}w_{0}}+2[2]^{2}C_{w_{0}}+4[2]^{2}C_{x_{\beta}w_{0}}\\
& &+([2]^{4}+2)C_{x_{\alpha}^{2}w_{0}}+([2]^{4}-[2]^{2}+1)C_{x_{\alpha}x_{\beta}w_{0}}\\
& &+(2[2]^{4}-[2]^{2}+1)C_{x_{\alpha}w_{0}}+([2]^{2}+1)C_{x_{\alpha}^{3}w_{0}}\\
& &+[2]^{2}C_{x_{\alpha}^{2}w_{0}}+([2]^{4}-[2]^{2})C_{x_{\alpha}w_{0}}+[2]^{2}C_{x_{\beta}w_{0}}\\
& &+[2]^{2}C_{w_{0}}+[2]^{2}C_{x_{\alpha}w_{0}}+([2]^{4}-2[2]^{2})C_{w_{0}}+[2]^{2}C_{x_{\alpha}x_{\beta}w_{0}}\\
& &+([2]^{4}-2[2]^{2})C_{x_{\beta}w_{0}}+[2]^{2}C_{x_{\alpha}w_{0}}+[2]^{2}C_{x_{\alpha}^{2}w_{0}}\\
& &+[2]^{2}C_{x_{\alpha}^{3}w_{0}}+[2]^{2}C_{x_{\alpha}w_{0}}+[2]^{2}C_{x_{\beta}w_{0}}+[2]^{2}C_{x_{\alpha}x_{\beta}w_{0}}\\
& &+([2]^{4}-[2]^{2})C_{x_{\alpha}^{2}w_{0}}+[2]^{2}C_{x_{\alpha}^{2}w_{0}}+([2]^{4}-[2]^{2})C_{x_{\alpha}w_{0}}\\
& &+[2]^{2}C_{x_{\beta}w_{0}}+[2]^{2}C_{w_{0}}-\{([2]^{4}-2[2]^{2})C_{x_{\alpha}x_{\beta}w_{0}}\\
& &+(3[2]^{4}-5[2]^{2})C_{x_{\alpha}w_{0}}+([2]^{4}-[2]^{2})C_{w_{0}}\\
& &+(2[2]^{4}-3[2]^{2})C_{x_{\alpha}^{2}w_{0}}+([2]^{4}-[2]^{2})C_{x_{\beta}w_{0}}\}\}\\
&=&[2]C_{x_{\alpha}^{2}x_{\beta}w_{0}}+[2]C_{x_{\beta}^{2}w_{0}}+3[2]^{3}C_{w_{0}}+6[2]^{3}C_{x_{\beta}w_{0}}\\
& &+(5[2]^{3}+2[2])C_{x_{\alpha}^{2}w_{0}}+(3[2]^{3}+[2])C_{x_{\alpha}x_{\beta}w_{0}}\\
& &+([2]^{5}+3[2]^{3}+[2])C_{x_{\alpha}w_{0}}+(2[2]^{3}+[2])C_{x_{\alpha}^{3}w_{0}}.\\
\end{eqnarray*}

\begin{eqnarray*}
C_{12121201}C_{x_{\alpha}^{2}w_{0}}&=&(C_{1212120}C_{1}-C_{121212})C_{x_{\alpha}^{2}w_{0}}\\
&=&[2]^{2}C_{x_{\alpha}^{3}w_{0}}+[2]^{2}C_{x_{\alpha}w_{0}}+[2]^{2}C_{x_{\beta}w_{0}}+[2]^{2}C_{x_{\alpha}x_{\beta}w_{0}}\\
& &+([2]^{6}-3[2]^{4}+2[2]^{2})C_{x_{\alpha}^{2}w_{0}}-([2]^{6}-4[2]^{4}+3[2]^{2})C_{x_{\alpha}^{2}w_{0}}\\
&=&[2]^{2}C_{x_{\alpha}^{3}w_{0}}+[2]^{2}C_{x_{\alpha}w_{0}}+[2]^{2}C_{x_{\beta}w_{0}}+[2]^{2}C_{x_{\alpha}x_{\beta}w_{0}}\\
& &+([2]^{4}-[2]^{2})C_{x_{\alpha}^{2}w_{0}}.\\
\end{eqnarray*}

Computing $C_{121212}C_{0121021210121210121212}$
\begin{eqnarray*}
C_{2}C_{0121021210121210121212}&=&C_{02121021210121210121212}+C_{210121212}+2C_{0121210121212}\\
& &+C_{210121210121212}+C_{0121210121210121212}+C_{2102121021210121212}\\
& &+C_{210121210121210121212},\\
C_{1}C_{02121021210121210121212}&=&C_{102121021210121210121212}+C_{x_{\alpha}x_{\beta}w_{0}},\\
C_{2}C_{102121021210121210121212}&=&C_{2102121021210121210121212}+C_{02121021210121210121212},\\
C_{1}C_{2102121021210121210121212}&=&C_{12102121021210121210121212}+C_{102121021210121210121212},\\
C_{2}C_{12102121021210121210121212}&=&C_{212102121021210121210121212}+C_{2102121021210121210121212},\\
C_{1}C_{212102121021210121210121212}&=&C_{x_{\alpha}^{2}x_{\beta}w_{0}}+C_{x_{\alpha}w_{0}}+C_{x_{\alpha}^{2}w_{0}}+C_{x_{\alpha}x_{\beta}w_{0}}+C_{x_{\beta}w_{0}}\\
& &+C_{x_{\beta}^{2}w_{0}}+C_{x_{\alpha}^{3}w_{0}}+C_{12102121021210121210121212}.\\
\end{eqnarray*}
\begin{eqnarray*}
C_{2}C_{0121021210121210121212}&=&C_{02121021210121210121212}+C_{210121212}+2C_{0121210121212}\\
& &+C_{210121210121212}+C_{0121210121210121212}\\
& &+C_{2102121021210121212}+C_{210121210121210121212},\\
C_{1}C_{2}C_{0121021210121210121212}&=&C_{102121021210121210121212}+C_{x_{\alpha}x_{\beta}w_{0}}+C_{1210121212}\\
& &+C_{10121212}+2C_{x_{\alpha}w_{0}}+C_{1210121210121212}\\
& &+3C_{10121210121212}+C_{1210121210121210121212}\\
& &+2C_{10121210121210121212}+C_{x_{\alpha}^{2}w_{0}}\\
& &+C_{12102121021210121212}+C_{102121021210121212}\\
C_{2}C_{1}C_{2}C_{0121021210121210121212}&=&C_{2102121021210121210121212}+C_{02121021210121210121212}\\
& &+[2]C_{x_{\alpha}x_{\beta}w_{0}}+C_{21210121212}+C_{210121212}+C_{210121212}\\
& &+C_{0121212}+2[2]C_{x_{\alpha}w_{0}}+C_{21210121210121212}\\
& &+C_{210121210121212}+3(C_{210121210121212}\\
& &+C_{0121210121212})+C_{21210121210121210121212}\\
& &+C_{210121210121210121212}+2(C_{210121210121210121212}\\
& &+C_{0121210121210121212})+[2]C_{x_{\alpha}^{2}w_{0}}\\
& &+C_{212102121021210121212}+C_{2102121021210121212}\\
& &+C_{2102121021210121212}+C_{02121021210121212}\\
&=&C_{2102121021210121210121212}+C_{02121021210121210121212}\\
& &+[2]C_{x_{\alpha}x_{\beta}w_{0}}+C_{21210121212}+2C_{210121212}+C_{0121212}\\
& &+2[2]C_{x_{\alpha}w_{0}}+C_{21210121210121212}+4C_{210121210121212}\\
& &+3C_{0121210121212}+C_{21210121210121210121212}\\
& &+3C_{210121210121210121212}+2C_{0121210121210121212}\\
& &+[2]C_{x_{\alpha}^{2}w_{0}}+C_{212102121021210121212}\\
& &+2C_{2102121021210121212}+C_{02121021210121212}\\
C_{1}C_{2}C_{1}C_{2}C_{0121021210121210121212}&=&C_{12102121021210121210121212}+C_{102121021210121210121212}\\
& &+C_{102121021210121210121212}+([2]^{2}+1)C_{x_{\alpha}x_{\beta}w_{0}}+C_{x_{\alpha}w_{0}}\\
& &+C_{1210121212}+2(C_{1210121212}+C_{10121212})+C_{10121212}+C_{w_{0}}\\
& &+2[2]^{2}C_{x_{\alpha}w_{0}}+(C_{x_{\alpha}^{2}w_{0}}+C_{x_{\alpha}w_{0}}+C_{w_{0}}+C_{x_{\beta}w_{0}}\\
& &+C_{1210121210121212})+4(C_{1210121210121212}+C_{10121210121212})\\
& &+3(C_{10121210121212}+C_{x_{\alpha}w_{0}})+(C_{x_{\alpha}^{3}w_{0}}+C_{x_{\alpha}w_{0}}+C_{x_{\alpha}^{2}w_{0}}\\
& &+C_{x_{\alpha}x_{\beta}w_{0}}+C_{x_{\beta}w_{0}}+C_{1210121210121210121212})\\
& &+3(C_{1210121210121210121212}+C_{10121210121210121212})\\
& &+2(C_{10121210121210121212}+C_{x_{\alpha}^{2}w_{0}})\\
& &+[2]^{2}C_{x_{\alpha}^{2}w_{0}}+(C_{x_{\alpha}x_{\beta}w_{0}}+C_{x_{\beta}w_{0}}\\
& &+C_{x_{\alpha}w_{0}}+C_{12102121021210121212})\\
& &+2(C_{12102121021210121212}+C_{102121021210121212})\\
& &+C_{102121021210121212}+C_{x_{\beta}w_{0}}\\
&=&C_{12102121021210121210121212}+2C_{102121021210121210121212}\\
& &([2]^{2}+1)C_{x_{\alpha}x_{\beta}w_{0}}+C_{x_{\alpha}w_{0}}+3C_{1210121212}\\
& &+3C_{10121212}+C_{w_{0}}+2[2]^{2}C_{x_{\alpha}w_{0}}+C_{x_{\alpha}^{2}w_{0}}+C_{x_{\alpha}w_{0}}+C_{w_{0}}\\
& &+C_{x_{\beta}w_{0}}+5C_{1210121210121212}+7C_{10121210121212}\\
& &+3C_{x_{\alpha}w_{0}}+(C_{x_{\alpha}^{3}w_{0}}+C_{x_{\alpha}w_{0}}+C_{x_{\alpha}^{2}w_{0}}\\
& &+C_{x_{\alpha}x_{\beta}w_{0}}+C_{x_{\beta}w_{0}})\\
& &+4C_{1210121210121210121212}+5C_{10121210121210121212}+2C_{x_{\alpha}^{2}w_{0}}\\
& &+[2]^{2}C_{x_{\alpha}^{2}w_{0}}+(C_{x_{\alpha}x_{\beta}w_{0}}+C_{x_{\beta}w_{0}}\\
& &+C_{x_{\alpha}w_{0}})+3C_{12102121021210121212}\\
& &+3C_{102121021210121212}+C_{x_{\beta}w_{0}}\\
&=&C_{12102121021210121210121212}+2C_{102121021210121210121212}\\
& &+3C_{1210121212}+3C_{10121212}+5C_{1210121210121212}\\
& &+7C_{10121210121212}+4C_{1210121210121210121212}\\
& &+5C_{10121210121210121212}+3C_{12102121021210121212}\\
& &+3C_{102121021210121212}+2C_{w_{0}}+([2]^{2}+4)C_{x_{\alpha}^{2}w_{0}}\\
& &+4C_{x_{\beta}w_{0}}+([2]^{2}+3)C_{x_{\alpha}x_{\beta}w_{0}}\\
& &+(2[2]^{2}+7)C_{x_{\alpha}w_{0}}+C_{x_{\alpha}^{3}w_{0}},\\
C_{2}C_{1}C_{2}C_{1}C_{2}C_{0121021210121210121212}&=&C_{212102121021210121210121212}+C_{2102121021210121210121212}\\
& &+2(C_{2102121021210121210121212}+C_{02121021210121210121212})\\
& &+3(C_{21210121212}+C_{210121212})+3(C_{210121212}+C_{0121212})\\
& &+5(C_{21210121210121212}+C_{210121210121212})\\
& &+7(C_{210121210121212}+C_{0121210121212})\\
& &+4(C_{21210121210121210121212}+C_{210121210121210121212})\\
& &+5(C_{210121210121210121212}+C_{0121210121210121212})\\
& &+3(C_{212102121021210121212}+C_{2102121021210121212})\\
& &+3(C_{2102121021210121212}+C_{02121021210121212})\\
& &+2[2]C_{w_{0}}+([2]^{3}+4[2])C_{x_{\alpha}^{2}w_{0}}+4[2]C_{x_{\beta}w_{0}}\\
& &+([2]^{3}+3[2])C_{x_{\alpha}x_{\beta}w_{0}}+(2[2]^{3}+7[2])C_{x_{\alpha}w_{0}}+[2]C_{x_{\alpha}^{3}w_{0}}\\
&=&C_{212102121021210121210121212}+3C_{2102121021210121210121212}\\
& &+2C_{02121021210121210121212}+3C_{21210121212}+6C_{210121212}\\
& &+3C_{0121212}+5C_{21210121210121212}+12C_{210121210121212}\\
& &+7C_{0121210121212}+4C_{21210121210121210121212}\\
& &+9C_{210121210121210121212}+5C_{0121210121210121212}\\
& &+3C_{212102121021210121212}+6C_{2102121021210121212}\\
& &+3C_{02121021210121212}+2[2]C_{w_{0}}+([2]^{3}+4[2])C_{x_{\alpha}^{2}w_{0}}\\
& &+4[2]C_{x_{\beta}w_{0}}+([2]^{3}+3[2])C_{x_{\alpha}x_{\beta}w_{0}}\\
& &+(2[2]^{3}+7[2])C_{x_{\alpha}w_{0}}+[2]C_{x_{\alpha}^{3}w_{0}},\\
C_{1}C_{2}C_{1}C_{2}C_{1}C_{2}C_{0121021210121210121212}&=&(C_{x_{\alpha}^{2}x_{\beta}w_{0}}+C_{x_{\alpha}w_{0}}+C_{x_{\alpha}^{2}w_{0}}+C_{x_{\alpha}x_{\beta}w_{0}}+C_{x_{\beta}w_{0}}\\
& &+C_{x_{\beta}^{2}w_{0}}+C_{x_{\alpha}^{3}w_{0}}+C_{12102121021210121210121212})\\
& &+3(C_{12102121021210121210121212}+C_{102121021210121210121212})\\
& &+2(C_{02121021210121210121212}+C_{x_{\alpha}x_{\beta}w_{0}})\\
& &+3(C_{x_{\alpha}w_{0}}+C_{1210121212})+6(C_{1210121212}+C_{10121212})\\
& &+3(C_{10121212}+C_{w_{0}})\\
& &+5(C_{x_{\alpha}^{2}w_{0}}+C_{x_{\alpha}w_{0}}+C_{w_{0}}\\
& &+C_{x_{\beta}w_{0}}+C_{1210121210121212})\\
& &+12(C_{1210121210121212}+C_{10121210121212})\\
& &+7(C_{10121210121212}+C_{x_{\alpha}w_{0}})\\
& &+4(C_{x_{\alpha}^{3}w_{0}}+C_{x_{\alpha}w_{0}}+C_{x_{\alpha}^{2}w_{0}}\\
& &+C_{x_{\alpha}x_{\beta}w_{0}}+C_{x_{\beta}w_{0}}+C_{1210121210121210121212})\\
& &+9(C_{1210121210121210121212}+C_{10121210121210121212})\\
& &+5(C_{10121210121210121212}+C_{x_{\alpha}^{2}w_{0}})\\
& &+3(C_{x_{\alpha}x_{\beta}w_{0}}+C_{x_{\beta}w_{0}}\\
& &+C_{x_{\alpha}w_{0}}+C_{12102121021210121212})\\
& &+6(C_{12102121021210121212}+C_{102121021210121212})\\
& &+3(C_{102121021210121212}+C_{x_{\beta}w_{0}})+2[2]^{2}C_{w_{0}}\\
& &+([2]^{4}+4[2]^{2})C_{x_{\alpha}^{2}w_{0}}+4[2]^{2}C_{x_{\beta}w_{0}}+([2]^{4}+3[2]^{2})C_{x_{\alpha}x_{\beta}w_{0}}\\
& &+(2[2]^{4}+7[2]^{2})C_{x_{\alpha}w_{0}}+[2]^{2}C_{x_{\alpha}^{3}w_{0}}\\
&=&C_{x_{\alpha}^{2}x_{\beta}w_{0}}+C_{x_{\beta}^{2}w_{0}}+4C_{12102121021210121210121212}\\
& &+5C_{102121021210121210121212}+9C_{1210121212}+9C_{10121212}\\
& &+17C_{1210121210121212}+19C_{10121210121212}\\
& &+13C_{1210121210121210121212}+14C_{10121210121210121212}\\
& &+9C_{12102121021210121212}+9C_{102121021210121212}\\
& &+(2[2]^{2}+5)C_{w_{0}}+([2]^{4}+4[2]^{2}+15)C_{x_{\alpha}^{2}w_{0}}\\
& &+(4[2]^{2}+16)C_{x_{\beta}w_{0}}+([2]^{4}+3[2]^{2}+10)C_{x_{\alpha}x_{\beta}w_{0}}\\
& &+(2[2]^{4}+7[2]^{2}+23)C_{x_{\alpha}w_{0}}+([2]^{2}+8)C_{x_{\alpha}^{3}w_{0}}\\
C_{121212}C_{0121021210121210121212}&=&C_{x_{\alpha}^{2}x_{\beta}w_{0}}+C_{x_{\beta}^{2}w_{0}}+4C_{12102121021210121210121212}\\
& &+5C_{102121021210121210121212}+9C_{1210121212}+9C_{10121212}\\
& &+17C_{1210121210121212}+19C_{10121210121212}\\
& &+13C_{1210121210121210121212}+14C_{10121210121210121212}\\
& &+9C_{12102121021210121212}+9C_{102121021210121212}\\
& &+(2[2]^{2}+8)C_{w_{0}}+([2]^{4}+4[2]^{2}+15)C_{x_{\alpha}^{2}w_{0}}\\
& &+(4[2]^{2}+16)C_{x_{\beta}w_{0}}+([2]^{4}+3[2]^{2}+10)C_{x_{\alpha}x_{\beta}w_{0}}\\
& &+(2[2]^{4}+7[2]^{2}+23)C_{x_{\alpha}w_{0}}+([2]^{2}+5)C_{x_{\alpha}^{3}w_{0}}\\
& &-4\{C_{12102121021210121210121212}+2C_{102121021210121210121212}\\
& &+3C_{1210121212}+3C_{10121212}+5C_{1210121210121212}\\
& &+7C_{10121210121212}+4C_{1210121210121210121212}\\
& &+5C_{10121210121210121212}+3C_{12102121021210121212}\\
& &+3C_{102121021210121212}+2C_{w_{0}}+([2]^{2}+4)C_{x_{\alpha}^{2}w_{0}}\\
& &+4C_{x_{\beta}w_{0}}+([2]^{2}+3)C_{x_{\alpha}x_{\beta}w_{0}}\\
& &+(2[2]^{2}+7)C_{x_{\alpha}w_{0}}+C_{x_{\alpha}^{3}w_{0}}\}\\
& &+3\{C_{102121021210121210121212}+C_{x_{\alpha}x_{\beta}w_{0}}+C_{1210121212}\\
& &+C_{10121212}+2C_{x_{\alpha}w_{0}}+C_{1210121210121212}+3C_{10121210121212}\\
& &+C_{1210121210121210121212}+2C_{10121210121210121212}+C_{x_{\alpha}^{2}w_{0}}\\
& &+C_{12102121021210121212}+C_{102121021210121212}\}\\
&=&C_{x_{\alpha}^{2}x_{\beta}w_{0}}+C_{x_{\beta}^{2}w_{0}}+2[2]^{2}C_{w_{0}}+4[2]^{2}C_{x_{\beta}w_{0}}\\
& &+([2]^{4}+2)C_{x_{\alpha}^{2}w_{0}}+([2]^{4}-[2]^{2}+1)C_{x_{\alpha}x_{\beta}w_{0}}\\
& &+(2[2]^{4}-[2]^{2}+1)C_{x_{\alpha}w_{0}}+([2]^{2}+1)C_{x_{\alpha}^{3}w_{0}}\\
\end{eqnarray*}

Computing $C_{121212}C_{21021210121210121212}$
\begin{eqnarray*}
C_{2}C_{21021210121210121212}&=&[2]C_{21021210121210121212},\\
C_{1}C_{21021210121210121212}&=&C_{121021210121210121212}+C_{1021210121210121212},\\
C_{1}C_{2}C_{21021210121210121212}&=&[2](C_{121021210121210121212}+C_{1021210121210121212}),\\
C_{2}C_{121021210121210121212}&=&C_{x_{\alpha}x_{\beta}w_{0}}+2C_{x_{\alpha}w_{0}}+C_{x_{\alpha}^{2}w_{0}}\\
& &+C_{21021210121210121212},\\
C_{2}C_{1021210121210121212}&=&C_{21021210121210121212}+C_{x_{\alpha}w_{0}}+C_{x_{\beta}w_{0}}\\
& &+C_{w_{0}}+C_{x_{\alpha}^{2}w_{0}}+C_{021210121210121212},\\
C_{2}C_{1}C_{2}C_{21021210121210121212}&=&[2](C_{x_{\alpha}x_{\beta}w_{0}}+3C_{x_{\alpha}w_{0}}+2C_{x_{\alpha}^{2}w_{0}}+C_{w_{0}}+C_{x_{\beta}w_{0}}\\
& &+2C_{21021210121210121212}+C_{021210121210121212}),\\
C_{1}C_{2}C_{1}C_{2}C_{21021210121210121212}&=&[2]\{[2]C_{x_{\alpha}x_{\beta}w_{0}}+3[2]C_{x_{\alpha}w_{0}}+2[2]C_{x_{\alpha}^{2}w_{0}}+[2]C_{w_{0}}\\
& &+[2]C_{x_{\beta}w_{0}}+2(C_{121021210121210121212}+C_{1021210121210121212})\\
& &+C_{1021210121210121212}+C_{01210121210121212}\}\\
&=&[2]\{[2]C_{x_{\alpha}x_{\beta}w_{0}}+3[2]C_{x_{\alpha}w_{0}}+2[2]C_{x_{\alpha}^{2}w_{0}}+[2]C_{w_{0}}\\
& &+[2]C_{x_{\beta}w_{0}}+2C_{121021210121210121212}+3C_{1021210121210121212}\\
& &+C_{01210121210121212}\},\\
C_{2}C_{1}C_{2}C_{1}C_{2}C_{21021210121210121212}&=&[2]\{[2]^{2}C_{x_{\alpha}x_{\beta}w_{0}}+3[2]^{2}C_{x_{\alpha}w_{0}}+2[2]^{2}C_{x_{\alpha}^{2}w_{0}}+[2]^{2}C_{w_{0}}\\
& &+[2]^{2}C_{x_{\beta}w_{0}}+2(C_{x_{\alpha}x_{\beta}w_{0}}+2C_{x_{\alpha}w_{0}}+C_{x_{\alpha}^{2}w_{0}}\\
& &+C_{21021210121210121212})+3(C_{21021210121210121212}+C_{x_{\alpha}w_{0}}+C_{x_{\beta}w_{0}}\\
& &+C_{w_{0}}+C_{x_{\alpha}^{2}w_{0}}+C_{021210121210121212})\\
& &+C_{021210121210121212}\}\\
&=&[2]\{([2]^{2}+2)C_{x_{\alpha}x_{\beta}w_{0}}+(3[2]^{2}+7)C_{x_{\alpha}w_{0}}+([2]^{2}+3)C_{w_{0}}\\
& &+(2[2]^{2}+5)C_{x_{\alpha}^{2}w_{0}}+([2]^{2}+3)C_{x_{\beta}w_{0}}\\
& &+5C_{21021210121210121212}+4C_{021210121210121212}\},\\
C_{1}C_{2}C_{1}C_{2}C_{1}C_{2}C_{21021210121210121212}&=&[2]\{([2]^{3}+2[2])C_{x_{\alpha}x_{\beta}w_{0}}+(3[2]^{3}+7[2])C_{x_{\alpha}w_{0}}\\
& &+([2]^{3}+3[2])C_{w_{0}}+(2[2]^{3}+5[2])C_{x_{\alpha}^{2}w_{0}}+([2]^{3}\\
& &+3[2])C_{x_{\beta}w_{0}}+5C_{121021210121210121212}\\
& &+9C_{1021210121210121212}+4C_{01210121210121212}\},\\
C_{121212}C_{21021210121210121212}&=&[2]\{([2]^{3}+2[2])C_{x_{\alpha}x_{\beta}w_{0}}+(3[2]^{3}+7[2])C_{x_{\alpha}w_{0}}\\
& &+([2]^{3}+3[2])C_{w_{0}}+(2[2]^{3}+5[2])C_{x_{\alpha}^{2}w_{0}}+([2]^{3}+3[2])C_{x_{\beta}w_{0}}\\
& &+5C_{121021210121210121212}+9C_{1021210121210121212}\\
& &+4C_{01210121210121212}\}-4[2]\{[2]C_{x_{\alpha}x_{\beta}w_{0}}\\
& &+3[2]C_{x_{\alpha}w_{0}}+2[2]C_{x_{\alpha}^{2}w_{0}}+[2]C_{w_{0}}+[2]C_{x_{\beta}w_{0}}\\
& &+2C_{121021210121210121212}+3C_{1021210121210121212}\\
& &+C_{01210121210121212}\}+3[2](C_{121021210121210121212}\\
& &+C_{1021210121210121212})\\
&=&([2]^{4}-2[2]^{2})C_{x_{\alpha}x_{\beta}w_{0}}+(3[2]^{4}-5[2]^{2})C_{x_{\alpha}w_{0}}\\
& &+([2]^{4}-[2]^{2})C_{w_{0}}+(2[2]^{4}-3[2]^{2})C_{x_{\alpha}^{2}w_{0}}\\
& &+([2]^{4}-[2]^{2})C_{x_{\beta}w_{0}}.\\
\end{eqnarray*}

(109) Computing $C_{12121201210212}C_{121212}$
\begin{eqnarray*}
C_{12121201210212}C_{121212}&=&C_{121212}C_{21201210121212}([2]^{4}-2[2]^{2})C_{x_{\beta}w_{0}}+([2]^{4}-[2]^{2})C_{x_{\alpha}w_{0}}\\
& &+([2]^{4}-2[2]^{2})C_{w_{0}}.
\end{eqnarray*}

(110) Computing $C_{12121201210212}C_{0121212}$
\begin{eqnarray*}
C_{12121201210212}C_{0121212}&=&C_{1212120}C_{21201210121212}==([2]^{3}-[2])C_{x_{\alpha}^{2}w_{0}}+(3[2]^{3}-2[2])C_{x_{\alpha}w_{0}}\\
& &+([2]^{3}-[2])C_{x_{\beta}w_{0}}+([2]^{5}-[2]^{3}-[2])C_{w_{0}}.
\end{eqnarray*}

(111) Computing $C_{12121201210212}C_{10121212}$
\begin{eqnarray*}
C_{12121201210212}C_{10121212}&=&C_{12121201}C_{21201210121212}C_{x_{\alpha}x_{\beta}w_{0}}+2[2]^{2}C_{x_{\alpha}^{2}w_{0}}\\
& &+2([2]^{4}-[2]^{2})C_{w_{0}}+([2]^{4}+3[2]^{2})C_{x_{\alpha}w_{0}}+2[2]^{2}C_{x_{\beta}w_{0}}.
\end{eqnarray*}

(112) Computing $C_{12121201210212}C_{210121212}$
\begin{eqnarray*}
C_{12121201210212}C_{210121212}&=&C_{121212012}C_{21201210121212}=[2]C_{x_{\alpha}x_{\beta}w_{0}}+([2]^{3}+[2])C_{x_{\alpha}^{2}w_{0}}\\
& &+([2]^{5}+2[2])C_{x_{\alpha}w_{0}}+([2]^{3}+[2])C_{x_{\beta}w_{0}}+([2]^{5}-[2]^{3}+[2])C_{w_{0}}.
\end{eqnarray*}

(113) Computing $C_{12121201210212}C_{1210121212}$
\begin{eqnarray*}
C_{12121201210212}C_{1210121212}&=&C_{1212120121}C_{21201210121212}=[2]^{2}C_{x_{\alpha}x_{\beta}w_{0}}+3[2]^{2}C_{x_{\alpha}^{2}w_{0}}\\
& &+([2]^{4}+[2]^{2})C_{w_{0}}+(2[2]^{4}+[2]^{2})C_{x_{\alpha}w_{0}}+3[2]^{2}C_{x_{\beta}w_{0}}.
\end{eqnarray*}

(114) Computing $C_{12121201210212}C_{21210121212}$
\begin{eqnarray*}
C_{12121201210212}C_{21210121212}&=&C_{12121201212}C_{21201210121212}=([2]^{3}-[2])C_{x_{\alpha}x_{\beta}w_{0}}\\
& &+(2[2]^{3}-[2])C_{x_{\alpha}^{2}w_{0}}+(2[2]^{3}-[2])C_{w_{0}}.
\end{eqnarray*}

(115) Computing $C_{12121201210212}C_{01210121212}$
\begin{eqnarray*}
C_{12121201210212}C_{01210121212}&=&C_{12121201210}C_{21201210121212}=[2]C_{x_{\alpha}^{3}w_{0}}+4[2]^{3}C_{x_{\alpha}w_{0}}\\
& &+([2]^{3}+3[2])C_{x_{\beta}w_{0}}+2[2]C_{x_{\alpha}x_{\beta}w_{0}}\\
& &+([2]^{3}+3[2])C_{x_{\alpha}^{2}w_{0}}+2[2]^{3}C_{w_{0}}.
\end{eqnarray*}

(116) Computing $C_{12121201210212}C_{201210121212}$
\begin{eqnarray*}
C_{12121201210212}C_{201210121212}&=&C_{121212012102}C_{21201210121212}=[2]^{2}C_{x_{\alpha}^{3}w_{0}}+4[2]^{4}C_{x_{\alpha}w_{0}}\\
& &+([2]^{4}+3[2]^{2})C_{x_{\beta}w_{0}}+2[2]^{2}C_{x_{\alpha}x_{\beta}w_{0}}\\
& &+([2]^{4}+3[2]^{2})C_{x_{\alpha}^{2}w_{0}}+2[2]^{4}C_{w_{0}}.
\end{eqnarray*}

(117) Computing $C_{12121201210212}C_{1201210121212}$
\begin{eqnarray*}
C_{12121201210212}C_{1201210121212}&=&C_{1212120121021}C_{21201210121212}=[2]C_{x_{\beta}^{2}w_{0}}+([2]^{3}+3[2])C_{x_{\alpha}x_{\beta}w_{0}}\\
& &+(5[2]^{3}+3[2])C_{x_{\alpha}^{2}w_{0}}+(10[2]^{3}-[2])C_{x_{\alpha}w_{0}}+(3[2]^{3}+5[2])C_{x_{\beta}w_{0}}\\
& &+([2]^{5}+2[2]^{3}+[2])C_{w_{0}}+2[2]C_{x_{\alpha}^{3}w_{0}}.
\end{eqnarray*}

(118) Computing $C_{12121201210212}C_{21201210121212}$
\begin{eqnarray*}
C_{12121201210212}&=&C_{1212120121021}C_{2}-C_{121212012102}-C_{121212012121},\\
C_{12121201210212}C_{21201210121212}&=&(C_{1212120121021}C_{2}-C_{121212012102}-C_{121212012121})C_{21201210121212}\\
&=&C_{1212120121021}C_{2}C_{21201210121212}-C_{121212012102}C_{21201210121212}\\
& &-C_{121212012121}C_{21201210121212}\\
&=&[2]C_{1212120121021}C_{21201210121212}-C_{121212012102}C_{21201210121212}\\
& &-C_{121212012121}C_{21201210121212}\\
&=&[2]^{2}C_{x_{\beta}^{2}w_{0}}+([2]^{4}+3[2]^{2})C_{x_{\alpha}x_{\beta}w_{0}}+(5[2]^{4}+3[2]^{2})C_{x_{\alpha}^{2}w_{0}}\\
& &+(10[2]^{4}-[2]^{2})C_{x_{\alpha}w_{0}}+(3[2]^{4}+5[2]^{2})C_{x_{\beta}w_{0}}\\
& &+([2]^{6}+2[2]^{4}+[2]^{2})C_{w_{0}}+2[2]^{2}C_{x_{\alpha}^{3}w_{0}}\\
& &-\{[2]^{2}C_{x_{\alpha}^{3}w_{0}}+4[2]^{4}C_{x_{\alpha}w_{0}}+([2]^{4}+3[2]^{2})C_{x_{\beta}w_{0}}\\
& &+2[2]^{2}C_{x_{\alpha}x_{\beta}w_{0}}+([2]^{4}+3[2]^{2})C_{x_{\alpha}^{2}w_{0}}+2[2]^{4}C_{w_{0}}\}\\
& &-\{([2]^{4}-2[2]^{2})C_{x_{\alpha}x_{\beta}w_{0}}+([2]^{4}-[2]^{2})C_{x_{\beta}w_{0}}\\
& &+(3[2]^{4}-5[2]^{2})C_{x_{\alpha}w_{0}}+(2[2]^{4}-3[2]^{2})C_{x_{\alpha}^{2}w_{0}}\\
& &+([2]^{4}-[2]^{2})C_{w_{0}}\}\\
&=&[2]^{2}C_{x_{\beta}^{2}w_{0}}+3[2]^{2}C_{x_{\alpha}x_{\beta}w_{0}}+(2[2]^{4}+3[2]^{2})C_{x_{\alpha}^{2}w_{0}}\\
& &+(3[2]^{4}+4[2]^{2})C_{x_{\alpha}w_{0}}+([2]^{4}+3[2]^{2})C_{x_{\beta}w_{0}}\\
& &+([2]^{6}-3[2]^{4}+2[2]^{2})C_{w_{0}}.\\
\end{eqnarray*}

\begin{eqnarray*}
C_{121212012121}C_{21201210121212}&=&(C_{12121201212}C_{1}-C_{1212120121})C_{21201210121212}\\
&=&C_{12121201212}C_{121201210121212}+C_{12121201212}C_{1201210121212}\\
& &-C_{1212120121}C_{21201210121212}\\
&=&([2]^{4}-2[2]^{2})C_{x_{\alpha}x_{\beta}w_{0}}+2[2]^{2}C_{x_{\beta}w_{0}}\\
& &+(3[2]^{4}-5[2]^{2})C_{x_{\alpha}w_{0}}+([2]^{4}-[2]^{2})C_{x_{\alpha}^{2}w_{0}}+2[2]^{2}C_{w_{0}}\\
& &+[2]^{2}C_{x_{\alpha}x_{\beta}w_{0}}+([2]^{4}+[2]^{2})C_{x_{\alpha}^{2}w_{0}}+(2[2]^{4}-2[2]^{2})C_{w_{0}}\\
& &+(2[2]^{4}+[2]^{2})C_{x_{\alpha}w_{0}}+[2]^{4}C_{x_{\beta}w_{0}}\\
& &-\{[2]^{2}C_{x_{\alpha}x_{\beta}w_{0}}+3[2]^{2}C_{x_{\alpha}^{2}w_{0}}+([2]^{4}+[2]^{2})C_{w_{0}}\\
& &+(2[2]^{4}+[2]^{2})C_{x_{\alpha}w_{0}}+3[2]^{2}C_{x_{\beta}w_{0}}\}\\
&=&([2]^{4}-2[2]^{2})C_{x_{\alpha}x_{\beta}w_{0}}+([2]^{4}-[2]^{2})C_{x_{\beta}w_{0}}\\
& &+(3[2]^{4}-5[2]^{2})C_{x_{\alpha}w_{0}}+(2[2]^{4}-3[2]^{2})C_{x_{\alpha}^{2}w_{0}}\\
& &+([2]^{4}-[2]^{2})C_{w_{0}}\\
\end{eqnarray*}

(119) Computing $C_{12121201210212}C_{121201210121212}$
\begin{eqnarray*}
C_{12121201210212}C_{121201210121212}&=&(C_{1212120121021}C_{2}-C_{121212012102}-C_{121212012121})C_{121201210121212}\\
&=&C_{1212120121021}C_{2}C_{121201210121212}-C_{121212012102}C_{121201210121212}\\
& &-C_{121212012121}C_{121201210121212}\\
&=&C_{1212120121021}C_{x_{\beta}w_{0}}+C_{1212120121021}C_{21201210121212}\\
& &+C_{1212120121021}C_{121212}-C_{121212012102}C_{121201210121212}\\
& &-C_{121212012121}C_{121201210121212}\\
&=&([2]^{5}-2[2]^{3}+[2])C_{x_{\alpha}x_{\beta}w_{0}}+([2]^{5}-[2]^{3})C_{x_{\alpha}^{2}w_{0}}\\
& &+([2]^{5}-2[2]^{3}+[2])C_{x_{\alpha}w_{0}}+2([2]^{3}-[2])C_{x_{\beta}w_{0}}\\
& &+([2]^{3}-[2])C_{w_{0}}+([2]^{3}-[2])C_{x_{\alpha}^{3}w_{0}}+([2]^{3}-[2])C_{x_{\beta}^{2}w_{0}}\\
& &+[2]C_{x_{\beta}^{2}w_{0}}+([2]^{3}+3[2])C_{x_{\alpha}x_{\beta}w_{0}}+(5[2]^{3}+3[2])C_{x_{\alpha}^{2}w_{0}}\\
& &+(10[2]^{3}-[2])C_{x_{\alpha}w_{0}}+(3[2]^{3}+5[2])C_{x_{\beta}w_{0}}\\
& &+([2]^{5}+2[2]^{3}+[2])C_{w_{0}}+2[2]C_{x_{\alpha}^{3}w_{0}}+([2]^{3}-[2])C_{x_{\beta}w_{0}}\\
& &+([2]^{5}-3[2]^{3}+[2])C_{x_{\alpha}w_{0}}+([2]^{3}-[2])C_{w_{0}}\\
& &-\{([2]^{5}-3[2]^{3}+[2])C_{x_{\alpha}x_{\beta}w_{0}}+([2]^{3}-[2])C_{x_{\beta}w_{0}}\\
& &+(2[2]^{5}-5[2]^{3}+[2])C_{x_{\alpha}w_{0}}+([2]^{5}-2[2]^{3})C_{x_{\alpha}^{2}w_{0}}\\
& &+([2]^{3}-[2])C_{w_{0}}\}-\{[2]C_{x_{\beta}^{2}w_{0}}+([2]^{3}+[2])C_{x_{\alpha}x_{\beta}w_{0}}\\
& &+(2[2]^{3}+3[2])C_{x_{\alpha}^{2}w_{0}}+2[2]C_{x_{\alpha}^{3}w_{0}}\\
& &+(6[2]^{3}-2[2])C_{x_{\alpha}w_{0}}+([2]^{3}+6[2])C_{x_{\beta}w_{0}}+(2[2]^{3}+2[2])C_{w_{0}}\}\\
&=&([2]^{3}+[2])C_{x_{\alpha}x_{\beta}w_{0}}+4[2]^{3}C_{x_{\alpha}^{2}w_{0}}+(4[2]^{3}+2[2])C_{x_{\alpha}w_{0}}\\
& &+(4[2]^{3}-3[2])C_{x_{\beta}w_{0}}+([2]^{5}+[2]^{3}-2[2])C_{w_{0}}\\
& &+([2]^{3}-[2])C_{x_{\alpha}^{3}w_{0}}+([2]^{3}-[2])C_{x_{\beta}^{2}w_{0}}\\
\end{eqnarray*}
We need to Computing $C_{1212120121021}C_{2121021210121212}$.

\begin{eqnarray*}
C_{1212120121021}C_{2121021210121212}&=&(C_{121212012102}C_{1}-C_{12121201210})C_{2121021210121212}\\
&=&([2]^{2}-1)C_{12121201210}C_{2121021210121212}\\
&=&([2]^{2}-1)\{([2]^{3}-[2])C_{x_{\alpha}x_{\beta}w_{0}}+[2]^{3}C_{x_{\alpha}^{2}w_{0}}+([2]^{3}-[2])C_{x_{\alpha}w_{0}}\\
& &+2[2]C_{x_{\beta}w_{0}}+[2]C_{w_{0}}+[2]C_{x_{\alpha}^{3}w_{0}}+[2]C_{x_{\beta}^{2}w_{0}}\}\\
&=&([2]^{5}-2[2]^{3}+[2])C_{x_{\alpha}x_{\beta}w_{0}}+([2]^{5}-[2]^{3})C_{x_{\alpha}^{2}w_{0}}\\
& &+([2]^{5}-2[2]^{3}+[2])C_{x_{\alpha}w_{0}}+2([2]^{3}-[2])C_{x_{\beta}w_{0}}\\
& &+([2]^{3}-[2])C_{w_{0}}+([2]^{3}-[2])C_{x_{\alpha}^{3}w_{0}}+([2]^{3}-[2])C_{x_{\beta}^{2}w_{0}}.\\
\end{eqnarray*}

\begin{eqnarray*}
C_{121212012121}C_{121201210121212}&=&(C_{12121201212}C_{1}-C_{1212120121})C_{121201210121212}\\
&=&[2]C_{12121201212}C_{121201210121212}-C_{1212120121}C_{121201210121212}\\
&=&([2]^{5}-2[2]^{3})C_{x_{\alpha}x_{\beta}w_{0}}+2[2]^{3}C_{x_{\beta}w_{0}}\\
& &+(3[2]^{5}-5[2]^{3})C_{x_{\alpha}w_{0}}+([2]^{5}-[2]^{3})C_{x_{\alpha}^{2}w_{0}}+2[2]^{3}C_{w_{0}}\\
& &-\{([2]^{3}-[2])C_{x_{\alpha}x_{\beta}w_{0}}+[2]^{3}C_{x_{\alpha}^{2}w_{0}}+([2]^{3}+[2])C_{w_{0}}\\
& &+([2]^{5}-[2])C_{x_{\alpha}w_{0}}+([2]^{3}+[2])C_{x_{\beta}w_{0}}\}\\
&=&([2]^{5}-3[2]^{3}+[2])C_{x_{\alpha}x_{\beta}w_{0}}+([2]^{3}-[2])C_{x_{\beta}w_{0}}\\
& &+(2[2]^{5}-5[2]^{3}+[2])C_{x_{\alpha}w_{0}}+([2]^{5}-2[2]^{3})C_{x_{\alpha}^{2}w_{0}}+([2]^{3}-[2])C_{w_{0}}\\
\end{eqnarray*}

(120) Computing $C_{12121201210212}C_{0121201210121212}$
\begin{eqnarray*}
C_{1212120121021}C_{2}&=&C_{12121201210212}+C_{121212012102}+C_{121212012121},\\
C_{12121201210212}&=&C_{1212120121021}C_{2}-C_{121212012102}-C_{121212012121},\\
C_{12121201210212}C_{0121201210121212}&=&(C_{1212120121021}C_{2}-C_{121212012102}-C_{121212012121})C_{0121201210121212}\\
&=&C_{1212120121021}C_{2}C_{0121201210121212}-C_{121212012102}C_{0121201210121212}\\
& &-C_{121212012121}C_{0121201210121212}\\
&=&C_{1212120121021}C_{20121201210121212}+C_{1212120121021}C_{210121210121212}\\
& &+C_{1212120121021}C_{0121212}-C_{121212012102}C_{0121201210121212}\\
& &-C_{121212012121}C_{0121201210121212}\\
&=&[2]^{2}C_{x_{\alpha}^{2}x_{\beta}w_{0}}+([2]^{4}+3[2]^{2})C_{x_{\beta}w_{0}}+([2]^{4}+[2]^{2})C_{x_{\alpha}^{2}w_{0}}\\
& &+3[2]^{2}C_{x_{\alpha}^{3}w_{0}}+2[2]^{2}C_{x_{\beta}^{2}w_{0}}+2[2]^{4}C_{x_{\alpha}x_{\beta}w_{0}}+2[2]^{2}C_{w_{0}}\\
& &+2[2]^{4}C_{x_{\alpha}w_{0}}+C_{x_{\alpha}^{2}x_{\beta}w_{0}}+C_{x_{\beta}^{2}w_{0}}+(2[2]^{4}+2[2]^{2})C_{w_{0}}\\
& &+(2[2]^{4}+6[2]^{2})C_{x_{\beta}w_{0}}+(2[2]^{4}+7[2]^{2}+2)C_{x_{\alpha}^{2}w_{0}}\\
& &+(6[2]^{2}+1)C_{x_{\alpha}x_{\beta}w_{0}}+(4[2]^{4}+5[2]^{2}+1)C_{x_{\alpha}w_{0}}\\
& &+(3[2]^{2}+1)C_{x_{\alpha}^{3}w_{0}}+[2]^{2}C_{x_{\alpha}^{2}w_{0}}\\
& &+(2[2]^{4}-[2]^{2})C_{x_{\alpha}121212}+2[2]^{2}C_{x_{\beta}121212}+[2]^{4}C_{w_{0}}\\
& &-\{[2]^{2}C_{x_{\beta}^{2}w_{0}}+([2]^{4}+[2]^{2})C_{x_{\alpha}x_{\beta}w_{0}}+(2[2]^{4}+2[2]^{2})C_{x_{\alpha}^{2}w_{0}}\\
& &+(4[2]^{4}-[2]^{2})C_{x_{\alpha}w_{0}}+([2]^{4}+4[2]^{2})C_{x_{\beta}w_{0}}+([2]^{4}+[2]^{2})C_{w_{0}}\\
& &+2[2]^{2}C_{x_{\alpha}^{3}w_{0}}\}\\
& &-\{C_{x_{\alpha}^{2}x_{\beta}w_{0}}+C_{x_{\beta}^{2}w_{0}}+2[2]^{2}C_{w_{0}}+4[2]^{2}C_{x_{\beta}w_{0}}\\
& &+([2]^{4}+[2]^{2}+1)C_{x_{\alpha}x_{\beta}w_{0}}+([2]^{4}+[2]^{2}+2)C_{x_{\alpha}^{2}w_{0}}\\
& &+(2[2]^{4}-2[2]^{2}+1)C_{x_{\alpha}w_{0}}+([2]^{2}+1)C_{x_{\alpha}^{3}w_{0}}\}\\
&=&[2]^{2}C_{x_{\alpha}^{2}x_{\beta}w_{0}}+(2[2]^{4}+3[2]^{2})C_{x_{\beta}w_{0}}+6[2]^{2}C_{x_{\alpha}^{2}w_{0}}\\
& &+3[2]^{2}C_{x_{\alpha}^{3}w_{0}}+[2]^{2}C_{x_{\beta}^{2}w_{0}}+4[2]^{2}C_{x_{\alpha}x_{\beta}w_{0}}\\
& &+(2[2]^{4}-[2]^{2})C_{w_{0}}+(3[2]^{4}+6[2]^{2})C_{x_{\alpha}w_{0}}.\\
\end{eqnarray*}

Computing $C_{1212120121021}C_{20121201210121212}$

\begin{eqnarray*}
C_{1212120121021}C_{20121201210121212}&=&(C_{121212012102}C_{1}-C_{12121201210})C_{20121201210121212}\\
&=&C_{121212012102}C_{120121201210121212}+C_{121212012102}C_{2121021210121212}\\
& &-C_{12121201210}C_{20121201210121212}\\
&=&C_{12121201210}C_{2120121201210121212}+C_{12121201210}C_{20121201210121212}\\
& &+C_{121212012102}C_{2121021210121212}-C_{12121201210}C_{20121201210121212}\\
&=&C_{12121201210}C_{2120121201210121212}+C_{121212012102}C_{2121021210121212}\\
&=&C_{12121201210}C_{2120121201210121212}+[2]C_{12121201210}C_{2121021210121212}\\
&=&[2]^{2}C_{x_{\alpha}^{2}x_{\beta}w_{0}}+([2]^{4}+[2]^{2})C_{x_{\beta}w_{0}}+[2]^{2}C_{x_{\alpha}^{2}w_{0}}+2[2]^{2}C_{x_{\alpha}^{3}w_{0}}\\
& &+[2]^{2}C_{x_{\beta}^{2}w_{0}}+([2]^{4}+[2]^{2})C_{x_{\alpha}x_{\beta}w_{0}}+[2]^{2}C_{w_{0}}+([2]^{4}+[2]^{2})C_{x_{\alpha}w_{0}}\\
& &+[2]^{2}C_{x_{\beta}^{2}w_{0}}+([2]^{4}-[2]^{2})C_{x_{\alpha}x_{\beta}w_{0}}\\
& &+[2]^{4}C_{x_{\alpha}^{2}w_{0}}+([2]^{4}-[2]^{2})C_{x_{\alpha}w_{0}}\\
& &+2[2]^{2}C_{x_{\beta}w_{0}}+[2]^{2}C_{w_{0}}+[2]^{2}C_{x_{\alpha}^{3}w_{0}}\\
&=&[2]^{2}C_{x_{\alpha}^{2}x_{\beta}w_{0}}+([2]^{4}+3[2]^{2})C_{x_{\beta}w_{0}}+([2]^{4}+[2]^{2})C_{x_{\alpha}^{2}w_{0}}\\
& &+3[2]^{2}C_{x_{\alpha}^{3}w_{0}}+2[2]^{2}C_{x_{\beta}^{2}w_{0}}\\
& &+2[2]^{4}C_{x_{\alpha}x_{\beta}w_{0}}+2[2]^{2}C_{w_{0}}+2[2]^{4}C_{x_{\alpha}w_{0}}\\
\end{eqnarray*}

\begin{eqnarray*}
C_{1212120121}C_{2102121021210121212}&=&(C_{121212012}C_{1}-C_{12121201})C_{2102121021210121212}\\
&=&C_{121212012}C_{12102121021210121212}+C_{121212012}C_{102121021210121212}\\
& &-C_{12121201}C_{2102121021210121212}\\
&=&C_{12121201}C_{212102121021210121212}+C_{12121201}C_{2102121021210121212}\\
& &-C_{1212120}C_{12102121021210121212}+C_{121212012}C_{102121021210121212}\\
& &-C_{12121201}C_{2102121021210121212}\\
&=&C_{12121201}C_{212102121021210121212}-C_{1212120}C_{12102121021210121212}\\
& &+C_{12121201}C_{2102121021210121212}+C_{12121201}C_{02121021210121212}\\
& &-C_{1212120}C_{102121021210121212}\\
&=&C_{12121201}C_{212102121021210121212}-C_{121212}C_{012102121021210121212}\\
& &-C_{121212}C_{2102121021210121212}+C_{12121201}C_{2102121021210121212}\\
& &+C_{12121201}C_{02121021210121212}-C_{1212120}C_{102121021210121212}\\
&=&C_{12121201}C_{212102121021210121212}-C_{121212}C_{012102121021210121212}\\
& &-C_{121212}C_{2102121021210121212}+C_{12121201}C_{2102121021210121212}\\
\end{eqnarray*}

\begin{eqnarray*}
C_{12121201}C_{02121021210121212}&=&(C_{1212120}C_{1}-C_{121212})C_{02121021210121212}\\
&=&C_{1212120}C_{102121021210121212}+C_{1212120}C_{2121021210121212}\\
& &-C_{121212}C_{02121021210121212}\\
&=&C_{1212120}C_{102121021210121212}+C_{121212}C_{02121021210121212}\\
& &-C_{121212}C_{02121021210121212}\\
&=&C_{1212120}C_{102121021210121212}
\end{eqnarray*}

\begin{eqnarray*}
C_{12121201210}C_{2120121201210121212}&=&[2]C_{1212120121}C_{2120121201210121212}-C_{121212012}C_{2120121201210121212}\\
&=&[2]C_{1212120121}C_{2120121201210121212}-\{[2]C_{12121201}C_{2120121201210121212}\\\
& &-[2]C_{121212}C_{2120121201210121212}\}\\
&=&[2]\{C_{12121201}C_{212102121021210121212}-C_{121212}C_{012102121021210121212}\\
& &-C_{121212}C_{2102121021210121212}+C_{12121201}C_{2102121021210121212}\}\\
& &-\{[2]C_{12121201}C_{2120121201210121212}-[2]C_{121212}C_{2120121201210121212}\}\\
&=&[2]\{C_{12121201}C_{212102121021210121212}-C_{121212}C_{012102121021210121212}\}\\
&=&[2]\{C_{1212120}(C_{x_{\alpha}x_{\beta}w_{0}}+C_{x_{\alpha}^{2}w_{0}}+C_{x_{\alpha}w_{0}})\\
& &+C_{121212}C_{012102121021210121212}+C_{121212}C_{2102121021210121212})\\
& &-C_{121212}C_{212102121021210121212}-C_{121212}C_{012102121021210121212}\}\\
&=&[2]\{C_{1212120}(C_{x_{\alpha}x_{\beta}w_{0}}+C_{x_{\alpha}^{2}w_{0}}+C_{x_{\alpha}w_{0}})\\
& &-C_{121212}C_{212102121021210121212}+C_{121212}C_{2102121021210121212})\}\\
&=&[2]\{[2]C_{x_{\alpha}^{2}x_{\beta}w_{0}}+[2]C_{x_{\beta}w_{0}}+[2]C_{x_{\alpha}^{2}w_{0}}+[2]C_{x_{\alpha}^{3}w_{0}}+[2]C_{x_{\beta}^{2}w_{0}}\\
& &+([2]^{5}-3[2]^{3}+2[2])C_{x_{\alpha}x_{\beta}w_{0}}\\
& &+[2]C_{x_{\alpha}^{3}w_{0}}+[2]C_{x_{\alpha}w_{0}}+[2]C_{x_{\beta}w_{0}}+[2]C_{x_{\alpha}x_{\beta}w_{0}}\\
& &+([2]^{5}-3[2]^{3}+2[2])C_{x_{\alpha}^{2}w_{0}}\\
& &+[2]C_{x_{\alpha}^{2}w_{0}}+([2]^{5}-3[2]^{3}+2[2])C_{x_{\alpha}w_{0}}\\
& &+[2]C_{x_{\beta}w_{0}}+[2]C_{w_{0}}\\
& &+([2]^{3}-[2])C_{x_{\alpha}x_{\beta}w_{0}}+([2]^{3}-[2])C_{x_{\alpha}^{2}w_{0}}\\
& &+([2]^{3}-[2])C_{x_{\alpha}w_{0}}+([2]^{3}-[2])C_{x_{\beta}w_{0}}\\
& &-\{([2]^{5}-3[2]^{3}+[2])C_{x_{\alpha}x_{\beta}w_{0}}+([2]^{5}-3[2]^{3}+[2])C_{x_{\alpha}^{2}w_{0}}\\
& &+([2]^{5}-3[2]^{3}+[2])C_{x_{\alpha}w_{0}}+[2]C_{x_{\beta}w_{0}}\}\}\\
&=&[2]^{2}C_{x_{\alpha}^{2}x_{\beta}w_{0}}+([2]^{4}+[2]^{2})C_{x_{\beta}w_{0}}+[2]^{2}C_{x_{\alpha}^{2}w_{0}}+2[2]^{2}C_{x_{\alpha}^{3}w_{0}}\\
& &+[2]^{2}C_{x_{\beta}^{2}w_{0}}+([2]^{4}+[2]^{2})C_{x_{\alpha}x_{\beta}w_{0}}+[2]^{2}C_{w_{0}}+([2]^{4}+[2]^{2})C_{x_{\alpha}w_{0}}.\\
\end{eqnarray*}

\begin{eqnarray*}
C_{121212}C_{212102121021210121212}&=&([2]^{5}-3[2]^{3}+[2])C_{x_{\alpha}x_{\beta}w_{0}}+([2]^{5}-3[2]^{3}+[2])C_{x_{\alpha}^{2}w_{0}}\\
& &+([2]^{5}-3[2]^{3}+[2])C_{x_{\alpha}w_{0}}+[2]C_{x_{\beta}w_{0}}.\\
\end{eqnarray*}

\begin{eqnarray*}
C_{121212}C_{2102121021210121212}&=&([2]^{3}-[2])C_{x_{\alpha}x_{\beta}w_{0}}+([2]^{3}-[2])C_{x_{\alpha}^{2}w_{0}}\\
& &+([2]^{3}-[2])C_{x_{\alpha}w_{0}}+([2]^{3}-[2])C_{x_{\beta}w_{0}}.\\
\end{eqnarray*}

\begin{eqnarray*}
C_{12121201}C_{212102121021210121212}&=&C_{1212120}(C_{x_{\alpha}x_{\beta}w_{0}}+C_{x_{\alpha}^{2}w_{0}}+C_{x_{\alpha}w_{0}}\\
& &+C_{12102121021210121212})-C_{121212}C_{212102121021210121212}\\
&=&C_{1212120}(C_{x_{\alpha}x_{\beta}w_{0}}+C_{x_{\alpha}^{2}w_{0}}+C_{x_{\alpha}w_{0}})+C_{121212}C_{012102121021210121212}\\
& &+C_{121212}C_{2102121021210121212}-C_{121212}C_{212102121021210121212}\\
\end{eqnarray*}

Computing $C_{1212120}C_{x_{\alpha}x_{\beta}w_{0}}=C_{121212}C_{01212102121021210121212}$.
\begin{eqnarray*}
C_{2}C_{01212102121021210121212}&=&[2]C_{01212102121021210121212},\\
C_{1}C_{01212102121021210121212}&=&C_{101212102121021210121212}+C_{x_{\alpha}x_{\beta}w_{0}},\\
C_{1}C_{2}C_{01212102121021210121212}&=&[2](C_{101212102121021210121212}+C_{x_{\alpha}x_{\beta}w_{0}}),\\
C_{2}C_{101212102121021210121212}&=&C_{2101212102121021210121212}+C_{01212102121021210121212},\\
C_{2}C_{1}C_{2}C_{01212102121021210121212}&=&[2](C_{2101212102121021210121212}+C_{01212102121021210121212}\\
& &+[2]C_{x_{\alpha}x_{\beta}w_{0}}),\\
C_{1}C_{2101212102121021210121212}&=&C_{12101212102121021210121212}+C_{101212102121021210121212},\\
C_{1}C_{2}C_{1}C_{2}C_{01212102121021210121212}&=&[2]\{C_{12101212102121021210121212}+2C_{101212102121021210121212}\\
& &+([2]^{2}+1)C_{x_{\alpha}x_{\beta}w_{0}}\},\\
C_{2}C_{1}C_{2}C_{1}C_{2}C_{01212102121021210121212}&=&[2]\{C_{212101212102121021210121212}+3C_{2101212102121021210121212}\\
& &+2C_{01212102121021210121212}+([2]^{3}+[2])C_{x_{\alpha}x_{\beta}w_{0}}\},\\
C_{1}C_{212101212102121021210121212}&=&C_{x_{\alpha}^{2}x_{\beta}w_{0}}+C_{x_{\beta}w_{0}}+C_{x_{\alpha}^{2}w_{0}}\\
& &+C_{x_{\alpha}x_{\beta}w_{0}}+C_{x_{\alpha}^{3}w_{0}}+C_{x_{\beta}^{2}w_{0}}+C_{12101212102121021210121212},\\
C_{1}C_{2}C_{1}C_{2}C_{1}C_{2}C_{01212102121021210121212}&=&[2]\{C_{x_{\alpha}^{2}x_{\beta}w_{0}}+C_{x_{\beta}w_{0}}+C_{x_{\alpha}^{2}w_{0}}+C_{x_{\alpha}^{3}w_{0}}+C_{x_{\beta}^{2}w_{0}}\\
& &+([2]^{4}+[2]^{2}+3)C_{x_{\alpha}x_{\beta}w_{0}}+C_{12101212102121021210121212}\\
& &+3C_{12101212102121021210121212}+5C_{101212102121021210121212}\},\\
C_{1212120}C_{x_{\alpha}x_{\beta}w_{0}}&=&C_{121212}C_{01212102121021210121212}\\
&=&[2]\{C_{x_{\alpha}^{2}x_{\beta}w_{0}}+C_{x_{\beta}w_{0}}+C_{x_{\alpha}^{2}w_{0}}+C_{x_{\alpha}^{3}w_{0}}+C_{x_{\beta}^{2}w_{0}}\\
& &+([2]^{4}+[2]^{2}+3)C_{x_{\alpha}x_{\beta}w_{0}}+C_{12101212102121021210121212}\\
& &+3C_{12101212102121021210121212}+5C_{101212102121021210121212}\}\\
& &-4[2]\{C_{12101212102121021210121212}+2C_{101212102121021210121212}\\
& &+([2]^{2}+1)C_{x_{\alpha}x_{\beta}w_{0}}\}+3[2](C_{101212102121021210121212}\\
& &+C_{x_{\alpha}x_{\beta}w_{0}})\\
&=&[2]C_{x_{\alpha}^{2}x_{\beta}w_{0}}+[2]C_{x_{\beta}w_{0}}+[2]C_{x_{\alpha}^{2}w_{0}}+[2]C_{x_{\alpha}^{3}w_{0}}+[2]C_{x_{\beta}^{2}w_{0}}\\
& &+([2]^{5}-3[2]^{3}+2[2])C_{x_{\alpha}x_{\beta}w_{0}}\\
\end{eqnarray*}

Computing $C_{1212120121021}C_{210121210121212}$
\begin{eqnarray*}
C_{1212120121021}C_{210121210121212}&=&(C_{121212012102}C_{1}-C_{12121201210})C_{210121210121212}\\
&=&C_{121212012102}C_{1210121210121212}+C_{121212012102}C_{10121210121212}\\
& &-C_{12121201210}C_{210121210121212}\\
&=&C_{12121201210}C_{21210121210121212}+C_{12121201210}C_{210121210121212}\\
& &+C_{12121201210}C_{210121210121212}+C_{12121201210}C_{210121210121212}\\
& &-C_{12121201210}C_{210121210121212}\\
&=&C_{12121201210}C_{21210121210121212}+C_{12121201210}C_{210121210121212}\\
& &+C_{12121201210}C_{0121210121212}\\
&=&C_{x_{\alpha}^{2}x_{\beta}w_{0}}+C_{x_{\beta}^{2}w_{0}}+3[2]^{2}C_{w_{0}}+6[2]^{2}C_{x_{\beta}w_{0}}\\
& &+(5[2]^{2}+2)C_{x_{\alpha}^{2}w_{0}}+(3[2]^{2}+1)C_{x_{\alpha}x_{\beta}w_{0}}\\
& &+([2]^{4}+3[2]^{2}+1)C_{x_{\alpha}w_{0}}+(2[2]^{2}+1)C_{x_{\alpha}^{3}w_{0}}\\
& &+[2]^{2}C_{x_{\alpha}^{3}w_{0}}+(2[2]^{4}+[2]^{2})C_{x_{\alpha}w_{0}}+([2]^{4}+[2]^{2})C_{x_{\beta}w_{0}}\\
& &+2[2]^{2}C_{x_{\alpha}x_{\beta}w_{0}}+([2]^{4}+2[2]^{2})C_{x_{\alpha}^{2}w_{0}}+[2]^{4}C_{w_{0}}\\
& &+[2]^{2}C_{x_{\alpha}x_{\beta}w_{0}}+[2]^{4}C_{x_{\alpha}^{2}w_{0}}+([2]^{4}+[2]^{2})C_{x_{\alpha}w_{0}}\\
& &+([2]^{4}-[2]^{2})C_{x_{\beta}w_{0}}+([2]^{4}-[2]^{2})C_{w_{0}}\\
&=&C_{x_{\alpha}^{2}x_{\beta}w_{0}}+C_{x_{\beta}^{2}w_{0}}+(2[2]^{4}+2[2]^{2})C_{w_{0}}+(2[2]^{4}+6[2]^{2})C_{x_{\beta}w_{0}}\\
& &+(2[2]^{4}+7[2]^{2}+2)C_{x_{\alpha}^{2}w_{0}}+(6[2]^{2}+1)C_{x_{\alpha}x_{\beta}w_{0}}\\
& &+(4[2]^{4}+5[2]^{2}+1)C_{x_{\alpha}w_{0}}+(3[2]^{2}+1)C_{x_{\alpha}^{3}w_{0}}.\\
\end{eqnarray*}

\begin{eqnarray*}
C_{12121201210}C_{0121210121212}&=&C_{12121201210}C_{0}C_{121210121212}\\
&=&[2]C_{12121201210}C_{121210121212}\\
&=&[2]^{2}C_{x_{\alpha}x_{\beta}w_{0}}+[2]^{4}C_{x_{\alpha}^{2}w_{0}}+([2]^{4}+[2]^{2})C_{x_{\alpha}w_{0}}\\\
& &+([2]^{4}-[2]^{2})C_{x_{\beta}w_{0}}+([2]^{4}-[2]^{2})C_{w_{0}}.\\
\end{eqnarray*}

\begin{eqnarray*}
C_{12121201210}C_{210121210121212}&=&[2]^{2}C_{x_{\alpha}^{3}w_{0}}+5[2]^{4}C_{x_{\alpha}w_{0}}+(3[2]^{4}-[2]^{2})C_{x_{\beta}w_{0}}\\
& &+4[2]^{2}C_{x_{\alpha}x_{\beta}w_{0}}+(3[2]^{4}+[2]^{2})C_{x_{\alpha}^{2}w_{0}}+(3[2]^{4}-2[2]^{2})C_{w_{0}}.\\
\end{eqnarray*}

\begin{eqnarray*}
C_{121212012102}C_{01210121210121212}&=&C_{12121201210}C_{2}C_{01210121210121212}\\
&=&C_{12121201210}C_{021210121210121212}\\
&=&C_{12121201210}C_{0}C_{21210121210121212}\\
&=&[2]C_{12121201210}C_{21210121210121212}\\
&=&[2]C_{x_{\alpha}^{2}x_{\beta}w_{0}}+[2]C_{x_{\beta}^{2}w_{0}}+3[2]^{3}C_{w_{0}}+6[2]^{3}C_{x_{\beta}w_{0}}\\
& &+(5[2]^{3}+2[2])C_{x_{\alpha}^{2}w_{0}}+(3[2]^{3}+[2])C_{x_{\alpha}x_{\beta}w_{0}}\\
& &+([2]^{5}+3[2]^{3}+[2])C_{x_{\alpha}w_{0}}+(2[2]^{3}+[2])C_{x_{\alpha}^{3}w_{0}}.\\
C_{12121201210}C_{21210121210121212}&=&C_{x_{\alpha}^{2}x_{\beta}w_{0}}+C_{x_{\beta}^{2}w_{0}}+3[2]^{2}C_{w_{0}}+6[2]^{2}C_{x_{\beta}w_{0}}\\
& &+(5[2]^{2}+2)C_{x_{\alpha}^{2}w_{0}}+(3[2]^{2}+1)C_{x_{\alpha}x_{\beta}w_{0}}\\
& &+([2]^{4}+3[2]^{2}+1)C_{x_{\alpha}w_{0}}+(2[2]^{2}+1)C_{x_{\alpha}^{3}w_{0}}\\
\end{eqnarray*}

Computing $C_{121212012121}C_{0121201210121212}$

\begin{eqnarray*}
C_{121212012121}C_{0121201210121212}&=&(C_{12121201212}C_{1}-C_{1212120121})C_{0121201210121212}\\
&=&C_{12121201212}C_{01210121210121212}+C_{12121201212}C_{121201210121212}\\
& &-C_{1212120121}C_{0121201210121212}\\
&=&C_{1212120121}C_{201210121210121212}-C_{121212012}C_{01210121210121212}\\
& &+C_{12121201212}C_{121201210121212}-C_{1212120121}C_{0121201210121212}\\
&=&C_{1212120121}C_{201210121210121212}-C_{121212012}C_{01210121210121212}\\
& &+C_{12121201212}C_{121201210121212}-C_{1212120121}C_{0121201210121212}\\
&=&C_{x_{\alpha}^{2}x_{\beta}w_{0}}+C_{x_{\beta}^{2}w_{0}}+3[2]^{2}C_{w_{0}}+6[2]^{2}C_{x_{\beta}w_{0}}\\
& &+(5[2]^{2}+2)C_{x_{\alpha}^{2}w_{0}}+(3[2]^{2}+1)C_{x_{\alpha}x_{\beta}w_{0}}\\
& &+([2]^{4}+3[2]^{2}+1)C_{x_{\alpha}w_{0}}+(2[2]^{2}+1)C_{x_{\alpha}^{3}w_{0}}\\
& &+([2]^{4}-2[2]^{2})C_{x_{\alpha}x_{\beta}w_{0}}+2[2]^{2}C_{x_{\beta}w_{0}}\\
& &+(3[2]^{4}-5[2]^{2})C_{x_{\alpha}w_{0}}+([2]^{4}-[2]^{2})C_{x_{\alpha}^{2}w_{0}}+2[2]^{2}C_{w_{0}}\\
& &-\{[2]^{2}C_{x_{\alpha}^{3}w_{0}}+2[2]^{4}C_{x_{\alpha}w_{0}}+4[2]^{2}C_{x_{\beta}w_{0}}\\
& &+2[2]^{2}C_{x_{\alpha}x_{\beta}w_{0}}+3[2]^{2}C_{x_{\alpha}^{2}w_{0}}+3[2]^{2}C_{w_{0}}\}\\
&=&C_{x_{\alpha}^{2}x_{\beta}w_{0}}+C_{x_{\beta}^{2}w_{0}}+2[2]^{2}C_{w_{0}}+4[2]^{2}C_{x_{\beta}w_{0}}\\
& &+([2]^{4}+[2]^{2}+1)C_{x_{\alpha}x_{\beta}w_{0}}+([2]^{4}+[2]^{2}+2)C_{x_{\alpha}^{2}w_{0}}\\
& &+(2[2]^{4}-2[2]^{2}+1)C_{x_{\alpha}w_{0}}+([2]^{2}+1)C_{x_{\alpha}^{3}w_{0}}.\\
\end{eqnarray*}

Computing $C_{1212120121}C_{021210121210121212}-C_{121212012}C_{01210121210121212}$
\begin{eqnarray*}
&&C_{1212120121}C_{021210121210121212}[2]-C_{121212012}C_{01210121210121212}[2]\\
&&=[2]C_{x_{\alpha}^{2}x_{\beta}w_{0}}+[2]C_{x_{\beta}^{2}w_{0}}+3[2]^{3}C_{w_{0}}+6[2]^{3}C_{x_{\beta}w_{0}}\\
&&+(5[2]^{3}+2[2])C_{x_{\alpha}^{2}w_{0}}+(3[2]^{3}+[2])C_{x_{\alpha}x_{\beta}w_{0}}\\
&&+([2]^{5}+3[2]^{3}+[2])C_{x_{\alpha}w_{0}}+(2[2]^{3}+[2])C_{x_{\alpha}^{3}w_{0}},\\
&&C_{1212120121}C_{021210121210121212}-C_{121212012}C_{01210121210121212}\\
&&=C_{x_{\alpha}^{2}x_{\beta}w_{0}}+C_{x_{\beta}^{2}w_{0}}+3[2]^{2}C_{w_{0}}+6[2]^{2}C_{x_{\beta}w_{0}}\\
&&+(5[2]^{2}+2)C_{x_{\alpha}^{2}w_{0}}+(3[2]^{2}+1)C_{x_{\alpha}x_{\beta}w_{0}}\\
&&+([2]^{4}+3[2]^{2}+1)C_{x_{\alpha}w_{0}}+(2[2]^{2}+1)C_{x_{\alpha}^{3}w_{0}}.\\
\end{eqnarray*}

(121) Computing $C_{121212012102121}C_{121212}$
\begin{eqnarray*}
C_{121212012102121}C_{121212}&=&C_{121212}C_{121201210121212}=([2]^{5}-3[2]^{3}+[2])C_{x_{\beta}w_{0}}\\
& &+([2]^{3}-[2])C_{x_{\alpha}w_{0}}+([2]^{5}-3[2]^{3}+[2])C_{w_{0}}.
\end{eqnarray*}

(122) Computing $C_{121212012102121}C_{0121212}$
\begin{eqnarray*}
C_{121212012102121}C_{0121212}&=&C_{1212120}C_{121201210121212}=C_{x_{\alpha}x_{\beta}w_{0}}+[2]^{2}C_{x_{\alpha}^{2}w_{0}}+(3[2]^{2})C_{x_{\alpha}w_{0}}\\
& &+([2]^{4}-2[2]^{2})C_{x_{\beta}w_{0}}+(2[2]^{4}-4[2]^{2})C_{w_{0}}.
\end{eqnarray*}

(123) Computing $C_{121212012102121}C_{10121212}$
\begin{eqnarray*}
C_{12121201210212}C_{10121212}&=&C_{12121201}C_{21201210121212}=[2]C_{x_{\alpha}x_{\beta}121212}+[2]^{3}C_{x_{\alpha}^{2}121212}\\
& &+([2]^{5}-[2]^{3}-[2])C_{w_{0}}+(2[2]^{3}+[2])C_{x_{\alpha}w_{0}}\\
& &+([2]^{3}-[2])C_{x_{\beta}w_{0}}.
\end{eqnarray*}

(124) Computing $C_{121212012102121}C_{210121212}$
\begin{eqnarray*}
C_{121212012102121}C_{210121212}&=&C_{121212012}C_{121201210121212}=[2]^{2}C_{x_{\alpha}x_{\beta}w_{0}}+2[2]^{2}C_{x_{\beta}w_{0}}+([2]^{4}\\
& &+2[2]^{2})C_{x_{\alpha}w_{0}}+2[2]^{2}C_{x_{\alpha}^{2}w_{0}}+[2]^{4}C_{w_{0}}.
\end{eqnarray*}

(125) Computing $C_{121212012102121}C_{1210121212}$
\begin{eqnarray*}
C_{121212012102121}C_{1210121212}&=&C_{1212120121}C_{121201210121212}=([2]^{3}-[2])C_{x_{\alpha}x_{\beta}w_{0}}+[2]^{3}C_{x_{\alpha}^{2}w_{0}}\\
& &+([2]^{3}+[2])C_{w_{0}}+([2]^{5}-[2])C_{x_{\alpha}w_{0}}+([2]^{3}+[2])C_{x_{\beta}w_{0}}.
\end{eqnarray*}

(126) Computing $C_{121212012102121}C_{21210121212}$
\begin{eqnarray*}
C_{12121201210212}C_{21210121212}&=&C_{12121201212}C_{121201210121212}=([2]^{4}-2[2]^{2})C_{x_{\alpha}x_{\beta}w_{0}}+2[2]^{2}C_{x_{\beta}w_{0}}\\
& &+(3[2]^{4}-5[2]^{2})C_{x_{\alpha}w_{0}}+([2]^{4}-[2]^{2})C_{x_{\alpha}^{2}w_{0}}+2[2]^{2}C_{w_{0}}.
\end{eqnarray*}

(127) Computing $C_{121212012102121}C_{01210121212}$
\begin{eqnarray*}
C_{12121201210212}C_{01210121212}&=&C_{12121201210}C_{121201210121212}=[2]^{2}C_{x_{\alpha}^{3}w_{0}}+2([2]^{4}-[2]^{2})C_{x_{\alpha}w_{0}}\\
& &+3[2]^{2}C_{x_{\beta}w_{0}}+[2]^{2}C_{x_{\alpha}x_{\beta}w_{0}}+2[2]^{2}C_{x_{\alpha}^{2}w_{0}}+2[2]^{2}C_{w_{0}}.
\end{eqnarray*}

(128) Computing $C_{121212012102121}C_{201210121212}$
\begin{eqnarray*}
C_{12121201210212}C_{201210121212}&=&C_{121212012102}C_{121201210121212}=[2]C_{x_{\beta}^{2}w_{0}}+([2]^{3}+[2])C_{x_{\alpha}x_{\beta}w_{0}}\\
& &+(2[2]^{3}+3[2])C_{x_{\alpha}^{2}w_{0}}+2[2]C_{x_{\alpha}^{3}w_{0}}\\
& &+(6[2]^{3}-2[2])C_{x_{\alpha}w_{0}}+([2]^{3}+6[2])C_{x_{\beta}w_{0}}+(2[2]^{3}+2[2])C_{w_{0}}.
\end{eqnarray*}

(129) Computing $C_{121212012102121}C_{1201210121212}$
\begin{eqnarray*}
C_{12121201210212}C_{1201210121212}&=&C_{1212120121021}C_{121201210121212}=[2]^{2}C_{x_{\beta}^{2}w_{0}}+[2]^{4}C_{x_{\alpha}x_{\beta}w_{0}}\\
& &+(2[2]^{4}+[2]^{2})C_{x_{\alpha}^{2}w_{0}}+4[2]^{4}C_{x_{\alpha}w_{0}}+([2]^{4}+3[2]^{2})C_{x_{\beta}w_{0}}\\
& &+2[2]^{4}C_{w_{0}}+[2]^{2}C_{x_{\alpha}^{3}w_{0}}.
\end{eqnarray*}

(130) Computing $C_{121212012102121}C_{21201210121212}$
\begin{eqnarray*}
C_{12121201210212}C_{21201210121212}&=&C_{12121201210212}C_{121201210121212}=([2]^{3}+[2])C_{x_{\alpha}x_{\beta}w_{0}}\\
& &+4[2]^{3}C_{x_{\alpha}^{2}w_{0}}+(4[2]^{3}+2[2])C_{x_{\alpha}w_{0}}\\
& &+(4[2]^{3}-3[2])C_{x_{\beta}w_{0}}+([2]^{5}+[2]^{3}-2[2])C_{w_{0}}\\
& &+([2]^{3}-[2])C_{x_{\alpha}^{3}w_{0}}+([2]^{3}-[2])C_{x_{\beta}^{2}w_{0}}.
\end{eqnarray*}

(131) Computing $C_{121212012102121}C_{121201210121212}$
\begin{eqnarray*}
C_{12121201210212}C_{121201210121212}&=&(C_{12121201210212}C_{1}-C_{1212120121021})C_{121201210121212}\\
&=&C_{12121201210212}C_{1}C_{121201210121212}-C_{1212120121021}C_{121201210121212}\\
&=&[2]C_{12121201210212}C_{121201210121212}-C_{1212120121021}C_{121201210121212}\\
&=&([2]^{4}+[2]^{2})C_{x_{\alpha}x_{\beta}w_{0}}+4[2]^{4}C_{x_{\alpha}^{2}w_{0}}+(4[2]^{4}+2[2]^{2})C_{x_{\alpha}w_{0}}\\
& &+(4[2]^{4}-3[2]^{2})C_{x_{\beta}w_{0}}+([2]^{6}+[2]^{4}-2[2]^{2})C_{w_{0}}\\
& &+([2]^{4}-2[2]^{2})C_{x_{\alpha}^{3}w_{0}}+([2]^{4}-2[2]^{2})C_{x_{\beta}^{2}w_{0}}\\
&&-\{[2]^{2}C_{x_{\beta}^{2}w_{0}}+[2]^{4}C_{x_{\alpha}x_{\beta}w_{0}}+(2[2]^{4}+[2]^{2})C_{x_{\alpha}^{2}w_{0}}\\
& &+4[2]^{4}C_{x_{\alpha}w_{0}}+([2]^{4}+3[2]^{2})C_{x_{\beta}w_{0}}+2[2]^{4}C_{w_{0}}+[2]^{2}C_{x_{\alpha}^{3}w_{0}}\}\\
&=&[2]^{4}C_{x_{\alpha}x_{\beta}w_{0}}+(2[2]^{4}-[2]^{2})C_{x_{\alpha}^{2}w_{0}}+2[2]^{2}C_{x_{\alpha}w_{0}}\\
& &+(3[2]^{4}-6[2]^{2})C_{x_{\beta}w_{0}}+([2]^{6}-[2]^{4}-2[2]^{2})C_{w_{0}}\\
& &+([2]^{4}-3[2]^{2})C_{x_{\alpha}^{3}w_{0}}+([2]^{4}-3[2]^{2})C_{x_{\beta}^{2}w_{0}}\\
\end{eqnarray*}

(132) Computing $C_{121212012102121}C_{0121201210121212}$
\begin{eqnarray*}
C_{12121201210212}C_{0121201210121212}&=&(C_{12121201210212}C_{1}-C_{1212120121021})C_{0121201210121212}\\
&=&C_{12121201210212}C_{1}C_{0121201210121212}-C_{1212120121021}C_{0121201210121212}\\
&=&C_{12121201210212}C_{01210121210121212}+C_{12121201210212}C_{121201210121212}\\
& &-C_{1212120121021}C_{0121201210121212}\\
&=&[2]C_{x_{\alpha}^{4}w_{0}}+5[2]C_{x_{\alpha}^{2}x_{\beta}w_{0}}+(8[2]^{3}+12[2])C_{x_{\alpha}^{2}w_{0}}+(3[2]^{3}+9[2])C_{x_{\alpha}x_{\beta}w_{0}}\\
& &+([2]^{5}+5[2]^{3}+9[2])C_{x_{\beta}w_{0}}+([2]^{3}+8[2])C_{x_{\alpha}^{3}w_{0}}\\
& &+([2]^{5}+7[2]^{3}+9[2])C_{x_{\alpha}w_{0}}+(4[2]^{3}+3[2])C_{w_{0}}+3[2]C_{x_{\beta}^{2}w_{0}}\\
& &+([2]^{3}+[2])C_{x_{\alpha}x_{\beta}w_{0}}+4[2]^{3}C_{x_{\alpha}^{2}w_{0}}+(4[2]^{3}+2[2])C_{x_{\alpha}w_{0}}\\
& &+(4[2]^{3}-3[2])C_{x_{\beta}w_{0}}+([2]^{5}+[2]^{3}-2[2])C_{w_{0}}\\
& &+([2]^{3}-[2])C_{x_{\alpha}^{3}w_{0}}+([2]^{3}-[2])C_{x_{\beta}^{2}w_{0}}\\
& &-\{[2]C_{x_{\alpha}^{2}x_{\beta}w_{0}}+2[2]C_{x_{\beta}^{2}w_{0}}+(4[2]^{3}+[2])C_{w_{0}}+(5[2]^{3}+5[2])C_{x_{\beta}w_{0}}\\
& &+(5[2]^{3}+5[2])C_{x_{\alpha}^{2}w_{0}}+(3[2]^{3}+2[2])C_{x_{\alpha}x_{\beta}w_{0}}\\
& &+9[2]^{3}C_{x_{\alpha}w_{0}}+([2]^{3}+3[2])C_{x_{\alpha}^{3}w_{0}}\}\\
&=&[2]C_{x_{\alpha}^{4}w_{0}}+4[2]C_{x_{\alpha}^{2}x_{\beta}w_{0}}+(7[2]^{3}+7[2])C_{x_{\alpha}^{2}w_{0}}+([2]^{3}+8[2])C_{x_{\alpha}x_{\beta}w_{0}}\\
& &+([2]^{5}+4[2]^{3}+[2])C_{x_{\beta}w_{0}}+([2]^{3}+4[2])C_{x_{\alpha}^{3}w_{0}}\\
& &+([2]^{5}+2[2]^{3}+11[2])C_{x_{\alpha}w_{0}}+([2]^{5}+[2]^{3}+2[2])C_{w_{0}}+[2]^{3}C_{x_{\beta}^{2}w_{0}}\\
\end{eqnarray*}

\begin{eqnarray*}
C_{12121201210212}C_{01210121210121212}&=&C_{1212120121021}C_{2}C_{01210121210121212}-C_{121212012102}C_{01210121210121212}\\
&=&C_{1212120121021}C_{021210121210121212}-C_{12121201210}C_{021210121210121212}\\
&=&C_{12121201210}(C_{21021210121210121212}+C_{x_{\beta}w_{0}}+C_{x_{\alpha}^{2}w_{0}}\\
& &+C_{x_{\alpha}w_{0}}+C_{w_{0}}+C_{021210121210121212})-C_{12121201210}C_{021210121210121212}\\
&=&C_{12121201210}(C_{21021210121210121212}+C_{x_{\beta}w_{0}}+C_{x_{\alpha}^{2}w_{0}}\\
& &+C_{x_{\alpha}w_{0}}+C_{w_{0}})\\
&=&[2]C_{x_{\alpha}^{4}w_{0}}+4[2]C_{x_{\alpha}^{2}x_{\beta}w_{0}}+(5[2]^{3}+11[2])C_{x_{\alpha}^{2}w_{0}}\\
& &+([2]^{3}+9[2])C_{x_{\alpha}x_{\beta}w_{0}}+([2]^{5}+3[2]^{3}+7[2])C_{x_{\beta}w_{0}}+7[2]C_{x_{\alpha}^{3}w_{0}}\\
& &+([2]^{5}+3[2]^{3}+10[2])C_{x_{\alpha}w_{0}}+(3[2]^{3}+2[2])C_{w_{0}}+2[2]C_{x_{\beta}^{2}w_{0}}\\
& &+[2]C_{x_{\beta}^{2}w_{0}}+([2]^{3}-[2])C_{x_{\alpha}x_{\beta}w_{0}}+[2]^{3}C_{x_{\alpha}^{2}w_{0}}\\
& &+([2]^{3}-[2])C_{x_{\alpha}w_{0}}+2[2]C_{x_{\beta}w_{0}}+[2]C_{w_{0}}+[2]C_{x_{\alpha}^{3}w_{0}}\\
& &+[2]C_{x_{\alpha}^{2}x_{\beta}w_{0}}+([2]^{3}+[2])C_{x_{\alpha}^{2}w_{0}}\\
& &+[2]^{3}C_{x_{\alpha}x_{\beta}w_{0}}+[2]^{3}C_{x_{\beta}w_{0}}+[2]^{3}C_{x_{\alpha}^{3}w_{0}}+[2]^{3}C_{x_{\alpha}w_{0}}\\
& &[2]C_{x_{\alpha}x_{\beta}w_{0}}+[2]^{3}C_{x_{\alpha}^{2}w_{0}}+([2]^{3}+[2])C_{x_{\alpha}w_{0}}\\
& &+([2]^{3}-[2])C_{x_{\beta}w_{0}}+([2]^{3}-[2])C_{w_{0}}\\
& &+[2]C_{x_{\beta}w_{0}}+([2]^{3}-[2])C_{x_{\alpha}w_{0}}+[2]C_{w_{0}}\\
&=&[2]C_{x_{\alpha}^{4}w_{0}}+5[2]C_{x_{\alpha}^{2}x_{\beta}w_{0}}+(8[2]^{3}+12[2])C_{x_{\alpha}^{2}w_{0}}\\
& &+(3[2]^{3}+9[2])C_{x_{\alpha}x_{\beta}w_{0}}+([2]^{5}+5[2]^{3}+9[2])C_{x_{\beta}w_{0}}\\
& &+([2]^{3}+8[2])C_{x_{\alpha}^{3}w_{0}}+([2]^{5}+7[2]^{3}+9[2])C_{x_{\alpha}w_{0}}\\
& &+(4[2]^{3}+3[2])C_{w_{0}}+3[2]C_{x_{\beta}^{2}w_{0}}\\
\end{eqnarray*}

\begin{eqnarray*}
C_{12121201210}C_{x_{\beta}w_{0}}&=&[2]C_{x_{\beta}^{2}w_{0}}+2([2]^{3}-[2])C_{x_{\alpha}x_{\beta}w_{0}}+(2[2]^{3}-[2])C_{x_{\alpha}^{2}w_{0}}\\
& &+2([2]^{3}-[2])C_{x_{\alpha}w_{0}}+([2]^{3}+[2])C_{x_{\beta}w_{0}}+[2]C_{w_{0}}+[2]C_{x_{\alpha}^{3}w_{0}}\\
& &-\{([2]^{3}-[2])C_{x_{\alpha}x_{\beta}w_{0}}+([2]^{3}-[2])C_{x_{\alpha}^{2}w_{0}}+([2]^{3}-[2])C_{x_{\alpha}w_{0}}\\
& &+([2]^{3}-[2])C_{x_{\beta}w_{0}}\}\\
&=&[2]C_{x_{\beta}^{2}w_{0}}+([2]^{3}-[2])C_{x_{\alpha}x_{\beta}w_{0}}+[2]^{3}C_{x_{\alpha}^{2}w_{0}}\\
& &+([2]^{3}-[2])C_{x_{\alpha}w_{0}}+2[2]C_{x_{\beta}w_{0}}+[2]C_{w_{0}}+[2]C_{x_{\alpha}^{3}w_{0}}\\
\end{eqnarray*}

\begin{eqnarray*}
C_{1212120121021}C_{021210121210121212}&=&(C_{121212012102}C_{1}-C_{12121201210})C_{021210121210121212}\\
&=&C_{121212012102}C_{1021210121210121212}+C_{121212012102}C_{01210121210121212}\\
& &-C_{12121201210}C_{021210121210121212}\\
&=&C_{121212012102}C_{1021210121210121212}+C_{12121201210}C_{021210121210121212}\\
& &-C_{12121201210}C_{021210121210121212}\\
&=&C_{121212012102}C_{1021210121210121212}\\
&=&C_{12121201210}C_{2}C_{1021210121210121212}\\
&=&C_{12121201210}(C_{21021210121210121212}+C_{x_{\beta}w_{0}}+C_{x_{\alpha}^{2}w_{0}}\\
& &+C_{x_{\alpha}w_{0}}+C_{w_{0}}+C_{021210121210121212})\\
\end{eqnarray*}

\begin{eqnarray*}
C_{12121201210}(C_{21021210121210121212}&=&(C_{1212120121}C_{0}-C_{121212012})C_{21021210121210121212}\\
&=&C_{1212120121}C_{0}C_{21021210121210121212}-C_{121212012}C_{21021210121210121212}\\
&=&C_{1212120121}(C_{210121210121210121212}+C_{210121212}+C_{0121210121212}\\
& &+2C_{210121210121212}+C_{2102121021210121212})\\
& &-C_{121212012}C_{21021210121210121212}\\
&=&[2]C_{x_{\alpha}^{4}w_{0}}+3[2]C_{x_{\alpha}^{2}x_{\beta}w_{0}}+([2]^{5}+5[2])C_{x_{\alpha}^{2}w_{0}}\\
& &+(2[2]^{3}+2[2]^{3})C_{x_{\alpha}x_{\beta}w_{0}}+(2[2]^{3}+[2])C_{x_{\beta}w_{0}}+(2[2]^{3}+3[2])C_{x_{\alpha}^{3}w_{0}}\\
& &+(2[2]^{3}+3[2]^{3})C_{x_{\alpha}w_{0}}+[2]C_{w_{0}}+[2]C_{x_{\beta}^{2}w_{0}}\\
& &+[2]C_{x_{\alpha}^{2}w_{0}}+2[2]C_{x_{\beta}w_{0}}\\
& &+2[2]^{3}C_{x_{\alpha}w_{0}}+([2]^{5}-2[2]^{3}+[2])C_{w_{0}}\\
& &+(2[2]^{3}-[2])C_{x_{\alpha}^{2}w_{0}}+[2]C_{x_{\alpha}x_{\beta}w_{0}}+2([2]^{3}-[2])C_{x_{\beta}w_{0}}\\
& &+(3[2]^{3}-[2]^{3})C_{x_{\alpha}w_{0}}+2([2]^{3}-[2])C_{w_{0}}\\
& &+2\{[2]C_{x_{\alpha}^{3}w_{0}}+([2]^{5}+4[2])C_{x_{\alpha}w_{0}}+(2[2]^{3}+[2])C_{x_{\beta}w_{0}}\\
& &+3[2]C_{x_{\alpha}x_{\beta}w_{0}}+(2[2]^{3}+3[2])C_{x_{\alpha}^{2}w_{0}}+2[2]^{3}C_{w_{0}}\}\\
& &+2[2]C_{x_{\alpha}^{2}x_{\beta}w_{0}}+3[2]C_{x_{\alpha}^{3}w_{0}}+(2[2]^{3}+[2])C_{x_{\alpha}w_{0}}\\
& &+([2]^{5}+4[2])C_{x_{\beta}w_{0}}+(2[2]^{3}+[2])C_{x_{\alpha}x_{\beta}w_{0}}\\
& &+([2]^{3}+2[2])C_{x_{\alpha}^{2}w_{0}}+2[2]C_{w_{0}}+2[2]C_{x_{\beta}^{2}w_{0}}\\
& &-\{[2]C_{x_{\alpha}^{2}x_{\beta}w_{0}}+(2[2]^{3}+[2])C_{x_{\alpha}^{3}w_{0}}+([2]^{5}+6[2]^{3}+[2])C_{x_{\alpha}w_{0}}\\
& &+5[2]^{3}C_{x_{\beta}w_{0}}+(3[2]^{3}+[2])C_{x_{\alpha}x_{\beta}w_{0}}\\
& &+([2]^{5}+2[2]^{3}+2[2])C_{x_{\alpha}^{2}w_{0}}+([2]^{5}+2[2]^{3})C_{w_{0}}+[2]C_{x_{\beta}^{2}w_{0}}\}\\
&=&[2]C_{x_{\alpha}^{4}w_{0}}+4[2]C_{x_{\alpha}^{2}x_{\beta}w_{0}}+(5[2]^{3}+11[2])C_{x_{\alpha}^{2}w_{0}}\\
& &+([2]^{3}+9[2])C_{x_{\alpha}x_{\beta}w_{0}}+([2]^{5}+3[2]^{3}+7[2])C_{x_{\beta}w_{0}}+7[2]C_{x_{\alpha}^{3}w_{0}}\\
& &+([2]^{5}+3[2]^{3}+10[2])C_{x_{\alpha}w_{0}}+(3[2]^{3}+2[2])C_{w_{0}}+2[2]C_{x_{\beta}^{2}w_{0}}\\
\end{eqnarray*}

\begin{eqnarray*}
C_{1212120121}C_{0121210121212}&=&C_{12121201210}C_{121210121212}+C_{121212012}C_{121210121212}\\
&=&(2[2]^{3}-[2])C_{x_{\alpha}^{2}w_{0}}+[2]C_{x_{\alpha}x_{\beta}w_{0}}+2([2]^{3}-[2])C_{x_{\beta}w_{0}}\\
& &+(3[2]^{3}-[2]^{3})C_{x_{\alpha}w_{0}}+2([2]^{3}-[2])C_{w_{0}}.\\
\end{eqnarray*}

\begin{eqnarray*}
C_{1212120121}C_{210121210121210121212}&=&(C_{121212012}C_{1}-C_{12121201})C_{210121210121210121212}\\
&=&C_{121212012}C_{1210121210121210121212}+C_{121212012}C_{10121210121210121212}\\
& &-C_{12121201}C_{210121210121210121212}\\
&=&C_{121212012}C_{1210121210121210121212}+C_{12121201}C_{210121210121210121212}\\
& &+C_{12121201}C_{0121210121210121212}-C_{1212120}C_{10121210121210121212}\\
& &-C_{12121201}C_{210121210121210121212}\\
&=&C_{121212012}C_{1210121210121210121212}+C_{12121201}C_{0121210121210121212}\\
& &-C_{1212120}C_{10121210121210121212}\\
&=&C_{121212012}C_{1210121210121210121212}\\
\end{eqnarray*}

\begin{eqnarray*}
C_{12121201}C_{0121210121210121212}&=&(C_{1212120}C_{1}-C_{121212})C_{0121210121210121212}\\
&=&C_{1212120}C_{1}C_{0121210121210121212}-C_{121212}C_{0121210121210121212}\\
&=&C_{1212120}C_{10121210121210121212}+C_{1212120}C_{121210121210121212}\\
& &-C_{121212}C_{0121210121210121212}\\
&=&C_{1212120}C_{10121210121210121212}\\
\end{eqnarray*}

\begin{eqnarray*}
C_{121212012}C_{1210121210121210121212}&=&(C_{12121201}C_{2}-C_{1212120})C_{1210121210121210121212}\\
&=&C_{12121201}C_{21210121210121210121212}+C_{12121201}C_{210121210121210121212}\\
& &-C_{1212120}C_{1210121210121210121212}\\
&=&(C_{1212120}C_{1}-C_{121212})C_{21210121210121210121212}\\
& &+C_{12121201}C_{210121210121210121212}-C_{1212120}C_{1210121210121210121212}\\
&=&C_{1212120}(C_{x_{\alpha}^{3}w_{0}}+C_{x_{\alpha}w_{0}}+C_{x_{\alpha}^{2}w_{0}}+C_{x_{\alpha}x_{\beta}w_{0}}\\
& &+C_{x_{\beta}w_{0}}+C_{1210121210121210121212})-C_{121212}C_{21210121210121210121212}\\
& &+C_{12121201}C_{210121210121210121212}-C_{1212120}C_{1210121210121210121212}\\
&=&C_{1212120}(C_{x_{\alpha}^{3}w_{0}}+C_{x_{\alpha}w_{0}}+C_{x_{\alpha}^{2}w_{0}}+C_{x_{\alpha}x_{\beta}w_{0}}+C_{x_{\beta}w_{0}})\\
& &-C_{121212}C_{21210121210121210121212}+C_{12121201}C_{210121210121210121212}\\
&=&[2]C_{x_{\alpha}^{4}w_{0}}+[2]C_{x_{\alpha}^{2}w_{0}}+[2]C_{x_{\alpha}^{2}x_{\beta}w_{0}}+[2]C_{x_{\beta}w_{0}}\\
& &+([2]^{5}-3[2]^{3}+2[2])C_{x_{\alpha}^{3}w_{0}}+[2]C_{x_{\alpha}^{2}w_{0}}+([2]^{5}-3[2]^{3}+2[2])C_{x_{\alpha}w_{0}}\\
& &+[2]C_{x_{\beta}w_{0}}+[2]C_{w_{0}}+[2]C_{x_{\alpha}^{3}w_{0}}+[2]C_{x_{\alpha}w_{0}}+[2]C_{x_{\beta}w_{0}}\\
& &+[2]C_{x_{\alpha}x_{\beta}w_{0}}+([2]^{5}-3[2]^{3}+2[2])C_{x_{\alpha}^{2}w_{0}}+[2]C_{x_{\alpha}^{2}x_{\beta}w_{0}}+[2]C_{x_{\beta}w_{0}}\\
& &+[2]C_{x_{\alpha}^{2}w_{0}}+[2]C_{x_{\alpha}^{3}w_{0}}+[2]C_{x_{\beta}^{2}w_{0}}+([2]^{5}-3[2]^{3}+2[2])C_{x_{\alpha}x_{\beta}w_{0}}\\
& &+[2]C_{x_{\alpha}x_{\beta}w_{0}}+([2]^{5}-3[2]^{3}+[2])C_{x_{\beta}w_{0}}+[2]C_{x_{\alpha}w_{0}}\\
& &+[2]C_{x_{\alpha}^{2}w_{0}}-\{([2]^{5}-3[2]^{3}+[2])C_{x_{\alpha}^{3}w_{0}}\\
& &+([2]^{5}-3[2]^{3}+[2])C_{x_{\alpha}w_{0}}+([2]^{5}-3[2]^{3}+[2])C_{x_{\beta}w_{0}}\\
& &+([2]^{5}-3[2]^{3}+[2])C_{x_{\alpha}x_{\beta}w_{0}}+([2]^{5}-3[2]^{3}+2[2])C_{x_{\alpha}^{2}w_{0}}\}\\
& &+[2]C_{x_{\alpha}^{2}x_{\beta}w_{0}}+([2]^{5}+[2])C_{x_{\alpha}^{2}w_{0}}+2[2]^{3}C_{x_{\alpha}x_{\beta}w_{0}}\\
& &+2[2]^{3}C_{x_{\beta}w_{0}}+2[2]^{3}C_{x_{\alpha}^{3}w_{0}}+2[2]^{3}C_{x_{\alpha}w_{0}}\\
&=&[2]C_{x_{\alpha}^{4}w_{0}}+3[2]C_{x_{\alpha}^{2}x_{\beta}w_{0}}+([2]^{5}+5[2])C_{x_{\alpha}^{2}w_{0}}+(2[2]^{3}+2[2]^{3})C_{x_{\alpha}x_{\beta}w_{0}}\\
& &+(2[2]^{3}+[2])C_{x_{\beta}w_{0}}+(2[2]^{3}+3[2])C_{x_{\alpha}^{3}w_{0}}\\
& &+(2[2]^{3}+3[2]^{3})C_{x_{\alpha}w_{0}}+[2]C_{w_{0}}+[2]C_{x_{\beta}^{2}w_{0}}.\\
\end{eqnarray*}

\begin{eqnarray*}
C_{121212}C_{21210121210121210121212}&=&C_{121212}(C_{2}C_{1210121210121210121212}-C_{210121210121210121212})\\
&=&[2]C_{121212}C_{1210121210121210121212}-C_{121212}C_{210121210121210121212}\\
&=&[2]C_{121212}(C_{1}C_{210121210121210121212}-C_{10121210121210121212}\\
& &-C_{121212}C_{210121210121210121212}\\
&=&([2]^{2}-1)C_{121212}C_{210121210121210121212}-[2]C_{121212}C_{10121210121210121212}\\
&=&([2]^{2}-1)\{([2]^{3}-[2])C_{x_{\alpha}^{3}w_{0}}+([2]^{3}-[2])C_{x_{\alpha}w_{0}}+([2]^{3}-[2])C_{x_{\beta}w_{0}}\\
& &+([2]^{3}-[2])C_{x_{\alpha}x_{\beta}w_{0}}+2([2]^{3}-[2])C_{x_{\alpha}^{2}w_{0}}\}-\{[2]^{3}C_{x_{\alpha}^{3}w_{0}}\\
& &+[2]^{3}C_{x_{\alpha}w_{0}}+[2]^{3}C_{x_{\beta}w_{0}}+[2]^{3}C_{x_{\alpha}x_{\beta}w_{0}}\\
& &+([2]^{5}-[2]^{3})C_{x_{\alpha}^{2}w_{0}}\}\\
&=&([2]^{5}-3[2]^{3}+[2])C_{x_{\alpha}^{3}w_{0}}+([2]^{5}-3[2]^{3}+[2])C_{x_{\alpha}w_{0}}\\
& &+([2]^{5}-3[2]^{3}+[2])C_{x_{\beta}w_{0}}+([2]^{5}-3[2]^{3}+[2])C_{x_{\alpha}x_{\beta}w_{0}}\\
& &+([2]^{5}-3[2]^{3}+2[2])C_{x_{\alpha}^{2}w_{0}}\\
\end{eqnarray*}

\begin{eqnarray*}
C_{121212}C_{210121210121210121212}&=&C_{121212}(C_{2}C_{10121210121210121212}-C_{0121210121210121212})\\
&=&[2]C_{121212}C_{10121210121210121212}-C_{121212}C_{0121210121210121212}\\
&=&[2]^{3}C_{x_{\alpha}^{3}w_{0}}+[2]^{3}C_{x_{\alpha}w_{0}}+[2]^{3}C_{x_{\beta}w_{0}}+[2]^{3}C_{x_{\alpha}x_{\beta}w_{0}}\\
& &+([2]^{5}-[2]^{3})C_{x_{\alpha}^{2}w_{0}}\\
& &-\{[2]C_{x_{\alpha}^{3}w_{0}}+[2]C_{x_{\alpha}w_{0}}+[2]C_{x_{\beta}w_{0}}+[2]C_{x_{\alpha}x_{\beta}w_{0}}\\
& &+([2]^{5}-3[2]^{3}+2[2])C_{x_{\alpha}^{2}w_{0}}\}\\
&=&([2]^{3}-[2])C_{x_{\alpha}^{3}w_{0}}+([2]^{3}-[2])C_{x_{\alpha}w_{0}}+([2]^{3}-[2])C_{x_{\beta}w_{0}}\\
& &+([2]^{3}-[2])C_{x_{\alpha}x_{\beta}w_{0}}+2([2]^{3}-[2])C_{x_{\alpha}^{2}w_{0}}\\
\end{eqnarray*}

\begin{eqnarray*}
C_{121212}C_{10121210121210121212}&=&C_{121212}(C_{1}C_{0121210121210121212}-C_{121210121210121212})\\
&=&[2]C_{121212}C_{0121210121210121212}-C_{121212}C_{121210121210121212}\\
C_{121212}C_{0121210121210121212}[2]&=&C_{121212}C_{10121210121210121212}+C_{121212}C_{121210121210121212}\\
&=&[2]^{2}C_{x_{\alpha}^{3}w_{0}}+[2]^{2}C_{x_{\alpha}w_{0}}+[2]^{2}C_{x_{\beta}w_{0}}+[2]^{2}C_{x_{\alpha}x_{\beta}w_{0}}\\
& &+([2]^{4}-[2]^{2})C_{x_{\alpha}^{2}w_{0}}+([2]^{6}-4[2]^{4}+3[2]^{2})C_{x_{\alpha}^{2}w_{0}}\\
&=&[2]^{2}C_{x_{\alpha}^{3}w_{0}}+[2]^{2}C_{x_{\alpha}w_{0}}+[2]^{2}C_{x_{\beta}w_{0}}+[2]^{2}C_{x_{\alpha}x_{\beta}w_{0}}\\
& &+([2]^{6}-3[2]^{4}+2[2]^{2})C_{x_{\alpha}^{2}w_{0}}\\
C_{121212}C_{0121210121210121212}&=&[2]C_{x_{\alpha}^{3}w_{0}}+[2]C_{x_{\alpha}w_{0}}+[2]C_{x_{\beta}w_{0}}+[2]C_{x_{\alpha}x_{\beta}w_{0}}\\
& &+([2]^{5}-3[2]^{3}+2[2])C_{x_{\alpha}^{2}w_{0}}.\\
\end{eqnarray*}

\begin{eqnarray*}
C_{12121201}C_{210121210121210121212}&=&(C_{1212120}C_{1}-C_{121212})C_{210121210121210121212}\\
&=&C_{1212120}C_{1}C_{210121210121210121212}-C_{121212}C_{210121210121210121212}\\
&=&C_{1212120}C_{1210121210121210121212}+C_{1212120}C_{10121210121210121212}\\
& &-C_{121212}C_{210121210121210121212}\\
&=&C_{1212120}C_{1210121210121210121212}+C_{1212120}C_{10121210121210121212}\\
& &-C_{121212}C_{210121210121210121212}\\
&=&C_{121212}C_{01210121210121210121212}+C_{121212}C_{210121210121210121212}\\
& &+C_{1212120}C_{10121210121210121212}-C_{121212}C_{210121210121210121212}\\
&=&C_{121212}C_{01210121210121210121212}+[2]C_{121212}C_{10121210121210121212}\\
&=&[2]C_{x_{\alpha}^{2}x_{\beta}w_{0}}+([2]^{3}+[2])C_{x_{\alpha}^{2}w_{0}}\\
& &+[2]^{3}C_{x_{\alpha}x_{\beta}w_{0}}+[2]^{3}C_{x_{\beta}w_{0}}+[2]^{3}C_{x_{\alpha}^{3}w_{0}}+[2]^{3}C_{x_{\alpha}w_{0}}\\
& &+[2]^{3}C_{x_{\alpha}^{3}w_{0}}+[2]^{3}C_{x_{\alpha}w_{0}}+[2]^{3}C_{x_{\beta}w_{0}}+[2]^{3}C_{x_{\alpha}x_{\beta}w_{0}}\\
& &+([2]^{5}-[2]^{3})C_{x_{\alpha}^{2}w_{0}}\\
&=&[2]C_{x_{\alpha}^{2}x_{\beta}w_{0}}+([2]^{5}+[2])C_{x_{\alpha}^{2}w_{0}}\\
& &+2[2]^{3}C_{x_{\alpha}x_{\beta}w_{0}}+2[2]^{3}C_{x_{\beta}w_{0}}+2[2]^{3}C_{x_{\alpha}^{3}w_{0}}+2[2]^{3}C_{x_{\alpha}w_{0}}\\
\end{eqnarray*}

\begin{eqnarray*}
C_{121212}C_{01210121210121210121212}&=&[2]C_{x_{\alpha}^{2}x_{\beta}w_{0}}+([2]^{3}+[2])C_{x_{\alpha}^{2}w_{0}}\\
& &+[2]^{3}C_{x_{\alpha}x_{\beta}w_{0}}+[2]^{3}C_{x_{\beta}w_{0}}+[2]^{3}C_{x_{\alpha}^{3}w_{0}}+[2]^{3}C_{x_{\alpha}w_{0}}.\\
\end{eqnarray*}

\begin{eqnarray*}
C_{1212120}C_{121210121210121210121212}&=&C_{121212}C_{0121210121210121210121212}\\
&=&[2]C_{x_{\alpha}^{4}w_{0}}+[2]C_{x_{\alpha}^{2}w_{0}}+[2]C_{x_{\alpha}^{2}x_{\beta}w_{0}}\\
& &+[2]C_{x_{\beta}w_{0}}+([2]^{5}-3[2]^{3}+2[2])C_{x_{\alpha}^{3}w_{0}}.\\
\end{eqnarray*}

(133) Computing $C_{1212120121021210}C_{121212}$
\begin{eqnarray*}
C_{1212120121021210}C_{121212}&=&C_{121212}C_{0121201210121212}=C_{x_{\alpha}x_{\beta}w_{0}}+[2]^{2}C_{x_{\alpha}^{2}w_{0}}+(2[2]^{2})C_{x_{\alpha}w_{0}}\\
& &+([2]^{4}-2[2]^{2})C_{x_{\beta}w_{0}}+([2]^{4}-2[2]^{2})C_{w_{0}}.
\end{eqnarray*}

(134) Computing $C_{1212120121021210}C_{0121212}$
\begin{eqnarray*}
C_{1212120121021210}C_{0121212}&=&C_{1212120}C_{0121201210121212}=[2]C_{121212}C_{0121201210121212}.
\end{eqnarray*}

(135) Computing $C_{1212120121021210}C_{10121212}$
\begin{eqnarray*}
C_{1212120121021210}C_{10121212}&=&C_{12121201}C_{0121201210121212}=[2]^{2}C_{x_{\alpha}x_{\beta}w_{0}}+(2[2]^{4}-3[2]^{2})C_{w_{0}}\\
& &+([2]^{4}+2[2]^{2})C_{x_{\alpha}w_{0}}+([2]^{4}-[2]^{2})C_{x_{\beta}w_{0}}+[2]^{4}C_{x_{\alpha}^{2}w_{0}}.
\end{eqnarray*}

(136) Computing $C_{1212120121021210}C_{210121212}$
\begin{eqnarray*}
C_{1212120121021210}C_{210121212}&=&C_{121212012}C_{0121201210121212}=[2]^{3}C_{x_{\alpha}x_{\beta}w_{0}}+(2[2]^{3}-[2])C_{x_{\beta}w_{0}}\\
& &+([2]^{5}+[2])C_{x_{\alpha}w_{0}}+2[2]^{3}C_{x_{\alpha}^{2}w_{0}}+(2[2]^{3}-[2])C_{w_{0}}.
\end{eqnarray*}

(137) Computing $C_{1212120121021210}C_{1210121212}$
\begin{eqnarray*}
C_{1212120121021210}C_{1210121212}&=&C_{1212120121}C_{0121201210121212}=[2]^{2}C_{x_{\alpha}^{3}w_{0}}+2[2]^{4}C_{x_{\alpha}w_{0}}\\
& &+4[2]^{2}C_{x_{\beta}w_{0}}+2[2]^{2}C_{x_{\alpha}x_{\beta}w_{0}}+4[2]^{2}C_{x_{\alpha}^{2}w_{0}}+3[2]^{2}C_{w_{0}}.
\end{eqnarray*}

(138) Computing $C_{1212120121021210}C_{21210121212}$
\begin{eqnarray*}
C_{1212120121021210}C_{21210121212}&=&C_{12121201212}C_{0121201210121212}=([2]^{3}+[2])C_{x_{\alpha}x_{\beta}w_{0}}\\
& &+(2[2]^{3}+[2])C_{x_{\alpha}^{2}w_{0}}+(4[2]^{3}-[2])C_{x_{\alpha}w_{0}}\\
& &+([2]^{3}+4[2])C_{x_{\beta}w_{0}}+([2]^{3}+2[2])C_{w_{0}}+2[2]C_{x_{\alpha}^{3}w_{0}}.
\end{eqnarray*}

(139) Computing $C_{1212120121021210}C_{01210121212}$
\begin{eqnarray*}
C_{1212120121021210}C_{01210121212}&=&C_{12121201210}C_{0121201210121212}=[2]^{3}C_{x_{\alpha}^{3}w_{0}}+[2]^{3}C_{x_{\alpha}x_{\beta}w_{0}}\\
& &+(2[2]^{3}+[2])C_{x_{\beta}w_{0}}+([2]^{5}-[2])C_{x_{\alpha}w_{0}}+2[2]^{3}C_{x_{\alpha}^{2}w_{0}}\\
& &+([2]^{3}+[2])C_{w_{0}}.
\end{eqnarray*}

(140) Computing $C_{1212120121021210}C_{201210121212}$
\begin{eqnarray*}
C_{1212120121021210}C_{201210121212}&=&C_{121212012102}C_{0121201210121212}=[2]^{2}C_{x_{\beta}^{2}w_{0}}+([2]^{4}+[2]^{2})C_{x_{\alpha}x_{\beta}w_{0}}\\
& &+(2[2]^{4}+2[2]^{2})C_{x_{\alpha}^{2}w_{0}}+(4[2]^{4}-[2]^{2})C_{x_{\alpha}w_{0}}+([2]^{4}+4[2]^{2})C_{x_{\beta}w_{0}}\\
& &+([2]^{4}+[2]^{2})C_{w_{0}}+2[2]^{2}C_{x_{\alpha}^{3}w_{0}}.
\end{eqnarray*}

(141) Computing $C_{1212120121021210}C_{1201210121212}$
\begin{eqnarray*}
C_{1212120121021210}C_{1201210121212}&=&C_{1212120121021}C_{0121201210121212}=[2]C_{x_{\alpha}^{2}x_{\beta}w_{0}}+2[2]C_{x_{\beta}^{2}w_{0}}\\
& &+(4[2]^{3}+[2])C_{w_{0}}+(5[2]^{3}+5[2])C_{x_{\beta}w_{0}}\\
& &+(5[2]^{3}+5[2])C_{x_{\alpha}^{2}w_{0}}+(3[2]^{3}+2[2])C_{x_{\alpha}x_{\beta}w_{0}}\\
& &+9[2]^{3}C_{x_{\alpha}w_{0}}+([2]^{3}+3[2])C_{x_{\alpha}^{3}w_{0}}.
\end{eqnarray*}

(142) Computing $C_{1212120121021210}C_{21201210121212}$
\begin{eqnarray*}
C_{1212120121021210}C_{21201210121212}&=&C_{12121201210212}C_{0121201210121212}=[2]^{2}C_{x_{\alpha}^{2}x_{\beta}w_{0}}\\
& &+(2[2]^{4}+3[2]^{2})C_{x_{\beta}w_{0}}+6[2]^{2}C_{x_{\alpha}^{2}w_{0}}\\
& &+3[2]^{2}C_{x_{\alpha}^{3}w_{0}}+[2]^{2}C_{x_{\beta}^{2}w_{0}}+4[2]^{2}C_{x_{\alpha}x_{\beta}w_{0}}\\
& &+(2[2]^{4}-[2]^{2})C_{w_{0}}+(3[2]^{4}+6[2]^{2})C_{x_{\alpha}w_{0}}.
\end{eqnarray*}

(143) Computing $C_{1212120121021210}C_{121201210121212}$
\begin{eqnarray*}
C_{1212120121021210}C_{121201210121212}&=&C_{121212012102121}C_{0121201210121212}=[2]C_{x_{\alpha}^{4}w_{0}}+4[2]C_{x_{\alpha}^{2}x_{\beta}w_{0}}\\
& &+(7[2]^{3}+7[2])C_{x_{\alpha}^{2}w_{0}}+([2]^{3}+8[2])C_{x_{\alpha}x_{\beta}w_{0}}\\
& &+([2]^{5}+4[2]^{3}+[2])C_{x_{\beta}w_{0}}+([2]^{3}+4[2])C_{x_{\alpha}^{3}w_{0}}\\
& &+([2]^{5}+2[2]^{3}+11[2])C_{x_{\alpha}w_{0}}+([2]^{5}+[2]^{3}+2[2])C_{w_{0}}\\
& &+[2]^{3}C_{x_{\beta}^{2}w_{0}}\\.
\end{eqnarray*}

(144) Computing $C_{1212120121021210}C_{0121201210121212}$
\begin{eqnarray*}
C_{121212012102121}C_{0}&=&C_{1212120121021210}+C_{12121201}.\\
C_{1212120121021210}C_{0121201210121212}&=&(C_{121212012102121}C_{0}-C_{12121201})C_{0121201210121212}\\
&=&[2]C_{121212012102121}C_{0121201210121212}-C_{12121201}C_{0121201210121212}\\
&=&[2]^{2}C_{x_{\alpha}^{4}w_{0}}+4[2]^{2}C_{x_{\alpha}^{2}x_{\beta}w_{0}}+(7[2]^{4}+7[2]^{2})C_{x_{\alpha}^{2}w_{0}}+([2]^{4}\\
& &+8[2]^{2})C_{x_{\alpha}x_{\beta}w_{0}}+([2]^{6}+4[2]^{4}+[2]^{2})C_{x_{\beta}w_{0}}+([2]^{4}+4[2]^{2})C_{x_{\alpha}^{3}w_{0}}\\
& &+([2]^{6}+2[2]^{4}+11[2]^{2})C_{x_{\alpha}w_{0}}+([2]^{6}+[2]^{4}+2[2]^{2})C_{w_{0}}+[2]^{4}C_{x_{\beta}^{2}w_{0}}\\
& &-\{[2]^{2}C_{x_{\alpha}x_{\beta}w_{0}}+(2[2]^{4}-3[2]^{2})C_{w_{0}}+([2]^{4}+2[2]^{2})C_{x_{\alpha}w_{0}}\\
& &+([2]^{4}-[2]^{2})C_{x_{\beta}w_{0}}+[2]^{4}C_{x_{\alpha}^{2}w_{0}}\}\\
&=&[2]^{2}C_{x_{\alpha}^{4}w_{0}}+4[2]^{2}C_{x_{\alpha}^{2}x_{\beta}w_{0}}+(6[2]^{4}+7[2]^{2})C_{x_{\alpha}^{2}w_{0}}+([2]^{4}\\
& &+7[2]^{2})C_{x_{\alpha}x_{\beta}w_{0}}+([2]^{6}+3[2]^{4}+2[2]^{2})C_{x_{\beta}w_{0}}+([2]^{4}+4[2]^{2})C_{x_{\alpha}^{3}w_{0}}\\
& &+([2]^{6}+[2]^{4}+9[2]^{2})C_{x_{\alpha}w_{0}}+([2]^{6}-[2]^{4}+5[2]^{2})C_{w_{0}}+[2]^{4}C_{x_{\beta}^{2}w_{0}}\\
\end{eqnarray*}
\quad\quad We know that $\delta_{w_{0}d_{u}^{-1},d_{u^{'}}w_{0},z_{1}w_{0}}$ is the coefficient of $q^{\frac{5}{2}}$ of
$h_{w_{0}d_{u}^{-1},d_{u^{'}}w_{0},z_{1}w_{0}}$. By the above process of calculation, we get the conclusion.
\endproof

Define

$\mathcal {U}_{1}=\{(d_{u},d_{u^{\prime}},z)|u, u^{\prime} \in W_{0}$, there exists unique $z\in
\{0,x_{\alpha},x_{\beta}\},\delta_{w_{0}d_{u}^{-1},d_{u^{\prime}}w_{0},zw_{0}}=1\}$,

\begin{eqnarray*}
\mathcal {U}_{2}&=&\{(d_{u},d_{u^{\prime}})|((d_{u}^{-1},d_{u^{\prime}},z_{1}),(d_{u}^{-1},d_{u^{\prime}},z_{2})),u, u^{\prime} \in
W_{0},(z_{1},z_{2})=(0,x_{\alpha}),\\
& &\delta_{w_{0}d_{u}^{-1},d_{u^{\prime}}w_{0},z_{1}w_{0}}=\delta_{w_{0}d_{u}^{-1},d_{u^{\prime}}w_{0},z_{2}w_{0}}=1\},
\end{eqnarray*}
\begin{eqnarray*}
\mathcal {U}_{3}&=&\{(d_{u},d_{u^{\prime}})|((d_{u}^{-1},d_{u^{\prime}},z_{1}),(d_{u}^{-1},d_{u^{\prime}},z_{2})), u, u^{\prime} \in
W_{0},(z_{1},z_{2})=(0,x_{\beta}),\\
& &\delta_{w_{0}d_{u}^{-1},d_{u^{\prime}}w_{0},z_{1}w_{0}}=\delta_{w_{0}d_{u}^{-1},d_{u^{\prime}}w_{0},z_{2}w_{0}}=1\},
\end{eqnarray*}
\begin{eqnarray*}
\mathcal{U}_{4}&=&\{(d_{u},d_{u^{\prime}})|((d_{u}^{-1},d_{u^{\prime}},z_{1}),(d_{u}^{-1},d_{u^{\prime}},z_{2}),(d_{u}^{-1},d_{u^{\prime}},z_{3})),
u,u^{\prime} \in W_{0}, (z_{1},z_{2},z_{3})=(0,x_{\alpha},x_{\beta}),\\
&
&\delta_{w_{0}d_{u}^{-1},d_{u^{\prime}}w_{0},z_{1}w_{0}}=\delta_{w_{0}d_{u}^{-1},d_{u^{\prime}}w_{0},z_{2}w_{0}}=\delta_{w_{0}d_{u}^{-1},d_{u^{\prime}}w_{0},z_{3}w_{0}}=1\}.
\end{eqnarray*}
Then by Proposition 4.1, we get
\begin{eqnarray*}
\mathcal {U}_{1}&=&\{(e,r,0),(r,e,0),(e,tstsr,x_{\alpha}),(tstsr,e,x_{\alpha}),(e,strstsr,x_{\alpha}),(strstsr,e,x_{\alpha}),\\
& &(sr,r,0),(r,trstsr,x_{\alpha}),(trstsr,r,x_{\alpha}),(r,tstrstsr,0),(tstrstsr,r,0),(r,trstsr,x_{\alpha}),\\
& &(trstsr,r,x_{\alpha}),(r,sr,0),(sr,tsr,0),(tsr,sr,0),(sr,rstsr,x_{\alpha}),(rstsr,sr,x_{\alpha}),\\
& &(sr,ststrstsr,0),(ststrstsr,sr,0),(tsr,stsr,0),(stsr,tsr,0),(tsr,rststrstsr,x_{\alpha}),\\
& &(rststrstsr,tsr,x_{\alpha}),(stsr,tstsr,0),(tstsr,stsr,0),(stsr,rstsr,0),(rstsr,stsr,0),\\
& &(stsr,ststrstsr,x_{\alpha}),(ststrstsr,stsr,x_{\alpha}),(tstsr,trstsr,0),(trstsr,tstsr,0),\\
& &(tstsr,tstrstsr,x_{\alpha}),(tstrstsr,tstsr,x_{\alpha}),(rstsr,trstsr,0),(trstsr,rstsr,0),\\
& &(rstsr,rststrstsr,x_{\alpha}),(rststrstsr,rstsr,x_{\alpha}),(trstsr,strstsr,0),(strstsr,trstsr,0),\\
& &(rstsrts,tstrstsr,0),(tstrstsr,strstsr,0),(strstsr,tstrstsr,0),(tstrstsr,strstsr,0),\\
& &(tstrstsr,ststrstsr,0),(ststrstsr,tstrstsr,0)\}.
\end{eqnarray*}
\begin{eqnarray*}
\mathcal {U}_{2}&=&\{(sr,strstsr),(strstsr,sr),(tsr,trstsr),(trstsr,tsr),(tsr,tstrstsr),\\
& &(tstrstsr,tsr),(stsr,strstsr),(strstsr,stsr)\}.\\
\mathcal {U}_{3}&=&\{(e,ststrstsr),(ststrstsr,e)\}.\\
\mathcal {U}_{4}&=&\{(ststrstsr,rststrstsr),(rststrstsr,ststrstsr)\}.
\end{eqnarray*}

The following lemma is a well-known result in representation theory (see
[24]).

\begin{lemma} \label{lemma4.2}
Assume that $L$ is a semisimple complex Lie algebra, $\lambda$ and $\lambda^{\prime}$
are dominant weights. Let $K_{\lambda,\gamma}$ be the weight multiplicity of weight $\gamma$ in the
irreducible $L$-module $V(\lambda)$, and $m_{\lambda,\lambda^{\prime},\nu}$ be the tensor product multiplicity of
irreducible module $V(\nu)$ in $V(\lambda)\bigotimes V(\lambda^{\prime})$. Then we have that $0 \leqslant m_{\lambda,\lambda^{\prime},\nu} \leqslant
K_{\lambda,\nu-\lambda^{\prime}}$.
\end{lemma}

\begin{theorem} \label{Thm}
For any $y, w\in c_{0}$, if $y\thicksim_{L}w$,that is,$y=d_{u}xw_{0}d_{v}, w=d_{u^{\prime}}x^{\prime}w_{0}d_{v}$, where $u, u^{\prime}, v \in W_{0}$
and $x,x^{\prime} \in \Lambda^{+}$. Then we must have:

(i). If there is a $z\in \{0, x_{\alpha}, x_{\beta}\}$ such that $(d_{u},d_{u^{\prime}},z)\in \mathcal {U}_{1}$, then we have
$\mu(y,w)=m_{x,x^{\prime},z}$.

(ii).If there is a $(d_{u},d_{u^{\prime}})\in \mathcal {U}_{2}$, then we have $\mu(y,w)=m_{x,x^{\prime},0}+m_{x,x^{\prime},x_{\alpha}}$.

(iii).If there is a $(d_{u},d_{u^{\prime}})\in \mathcal {U}_{3}$, then we have $\mu(y,w)=m_{x,x^{\prime},0}+m_{x,x^{\prime},x_{\beta}}$.

(iv).If there is a $(d_{u},d_{u^{\prime}})\in \mathcal {U}_{4}$, then we have
$\mu(y,w)=m_{x,x^{\prime},0}+m_{x,x^{\prime},x_{\alpha}}+m_{x,x^{\prime},x_{\beta}}$.

Otherwise, $\mu(y,w)=0$. Moreover, for any $y, w\in c_{0}$, we get $\mu(y,w)\leqslant 3$.
\end{theorem}

\beginproof
For any $y, w\in c_{0}$, if $\mu(y,w)\neq0$, then we know that $y\thicksim_{L}w$ or $y\thicksim_{R}w$. If $y\thicksim_{R}w$, then
$y^{-1}\thicksim_{L}w^{-1}$. Since we have $\mu(y,w)=\mu(y^{-1},w^{-1})$, we can assume $y\thicksim_{L}w$.

If $y=d_{u}xw_{0}d_{v}, w=d_{u^{\prime}}x^{\prime}w_{0}d_{v}$, for some $u, u^{\prime}, v \in W_{0}$ and $x,x^{\prime} \in \Lambda^{+}$.
By II.B(c) in section II, we also can assume $y=d_{u}xw_{0}, w=d_{u^{\prime}}x^{\prime}w_{0}$. By II.II(d) in section II, we know that

$$\mu(y,w)=\sum_{z\in \Lambda^{+}} m_{x^{\ast},x^{\prime},z^{\ast}}\delta_{w_{0}d_{u}^{-1},d_{u^{\prime}}w_{0},zw_{0}},$$

where $x^{\ast}=w_{0}x^{-1}w_{0}, z^{\ast}=w_{0}z^{-1}w_{0} \in \Lambda^{+}$.

From this, we get:

(1). If there is a $z\in \{0, x_{\alpha}, x_{\beta}\}$ such that $(d_{u},d_{u^{\prime}},z)\in \mathcal {U}_{1}$, then
$\mu(y,w)=m_{x,x^{\prime},z}$.

(2). If there is a $(d_{u},d_{u^{\prime}})\in \mathcal {U}_{2}$, then $\mu(y,w)=m_{x,x^{\prime},0}+m_{x,x^{\prime},x_{\alpha}}$.

(3). If there is a $(d_{u},d_{u^{\prime}})\in \mathcal {U}_{3}$, then $\mu(y,w)=m_{x,x^{\prime},0}+m_{x,x_{\prime},x_{\beta}}$.

(4). If there is a $(d_{u},d_{u^{\prime}})\in \mathcal {U}_{4}$, then
$\mu(y,w)=m_{x,x^{\prime},0}+m_{x,x^{\prime},x_{\alpha}}+m_{x,x^{\prime},x_{\beta}}$.

Since we have
\begin{eqnarray*}
m_{x^{\ast},x^{\prime},z^{\ast}}&=&dim Hom_{G}(V(x^{\ast})\bigotimes V(x^{\prime}),V(z^{\ast}))\\
&=&dim Hom_{G}(V(z^{\ast})\bigotimes V(x^{\prime}),V(x^{\ast}))\\
&=&m_{z,x^{\prime},x}.\\
\end{eqnarray*}
By the Freudenthal formula in Lie algebra and lemma VI.2, we can get if $z\in \{0,x_{\alpha},x_{\beta}\}$, then we have $dimV(z)_{\lambda}$ for any weight of
$V(z)$,
the by the lemma 2.2(e), we get $m_{x^{\ast},x^{\prime},z^{\ast}}=m_{z,x^{\prime},x}\leqslant 1$ for $x,x^{\prime} \in \Lambda^{+}$ and $z\in
\{0,x_{\alpha},x_{\beta}\}$. Then we get the theorem.
\endproof

\textbf{Acknowledgments:}
I would like to thank Xun Xie, Yanmin Yang and Zongxing Xiong for helpful discussions. I am very grateful to the referee for carefully reading
and helpful comments.

\end{document}